\documentclass[a4paper, 10pt, reqno]{amsart}
\usepackage[T1]{fontenc}
\usepackage{lmodern}

\usepackage{amsmath,amssymb,graphicx, mathtools} 
\usepackage{amsthm,amsaddr}
\usepackage[a4paper,left=2.5cm,right=2.5cm,top=3cm,bottom=3cm]{geometry}
\usepackage{enumerate} 
\usepackage{enumitem}
\usepackage{tikz} 
\usetikzlibrary{arrows,decorations.markings}
\usepackage{bm}
\usepackage[colorlinks, linkcolor=red, urlcolor=blue, citecolor=blue]{hyperref}
\usepackage{url}    

\usepackage[giveninits=true,url=false,doi=false,isbn=false,eprint=true,datamodel=mrnumber,sorting=nyt,sortcites=false,maxcitenames=4,maxbibnames=99, maxsortnames=99, backref=false,block=space,backend=biber,style=phys, biblabel=brackets]{biblatex} 
\AtEveryBibitem{\clearfield{month}}
\AtEveryCitekey{\clearfield{month}} 
\renewbibmacro{in:}{}
\DeclareFieldFormat[article]{title}{\emph{#1}} 
\DeclareFieldFormat{mrnumber}{\ifhyperref{\href{http://www.ams.org/mathscinet-getitem?mr=#1}{\nolinkurl{MR#1}}}{\nolinkurl{#1}}}
\DeclareFieldFormat{pmid}{\ifhyperref{\href{https://www.ncbi.nlm.nih.gov/pubmed/#1}{\nolinkurl{PMID#1}}}{\nolinkurl{#1}}}
\DeclareFieldFormat{eprint}{\ifhyperref{\href{https://arxiv.org/abs/#1}{\nolinkurl{arXiv:#1}}}{\nolinkurl{#1}}}
\renewbibmacro*{doi+eprint+url}{%
  \iftoggle{bbx:doi}{\printfield{doi}}{}
  \newunit\newblock%
  \printfield{mrnumber}%
  \newunit\newblock%
  \printfield{pmid}%
  \newunit\newblock%
  \printfield{eprint}%
  \iftoggle{bbx:url}{\usebibmacro{url+urldate}}{}}
\bibliography{bibliography}
\usepackage{caption}
\usepackage{subcaption}

\DeclareMathOperator{\Tr}{Tr}

\newcommand{\R}{\mathbb{R}} 
\newcommand{\E}{\mathbb{E}}

\theoremstyle{plain}
\newtheorem{thm}{Theorem}[section]
\newtheorem{cor}[thm]{Corollary}
\newtheorem{lem}[thm]{Lemma}
\newtheorem{prop}[thm]{Proposition}

\theoremstyle{definition}
\newtheorem{defn}[thm]{Definition}

\theoremstyle{remark}
\newtheorem{rmk}[thm]{Remark}

\numberwithin{equation}{section}

\makeatletter
\renewcommand\subsection{\@startsection{subsection}{2}%
  \z@{-.5\linespacing\@plus-.7\linespacing}{.5\linespacing}%
  {\normalfont\scshape}}
\renewcommand\subsubsection{\@startsection{subsubsection}{3}%
  \z@{.5\linespacing\@plus.7\linespacing}{-.5em}%
  {\normalfont\scshape}}
\makeatother

\makeatletter
\@namedef{subjclassname@2020}{%
  \textup{2020} Mathematics Subject Classification}
\makeatother

\usepackage[format=plain,
            labelfont=sc]{caption}
\usepackage[toc,page]{appendix} 
\usepackage{scalerel,stackengine}
\stackMath%
\newcommand\reallywidehat[1]{%
\savestack{\tmpbox}{\stretchto{%
  \scaleto{%
    \scalerel*[\widthof{\ensuremath{#1}}]{\kern.1pt\mathchar"0362\kern.1pt}%
    {\rule{0ex}{\textheight}}
  }{\textheight}%
}{2.4ex}}%
\stackon[-6.9pt]{#1}{\tmpbox}%
}

\usepackage{fancyhdr}

\newcommand\shortitle{Complexity of spiked random polynomials and finite-rank spherical integrals}
\newcommand\name{Vanessa Piccolo}

\begin{document}

\title{Topological complexity of spiked random polynomials and finite-rank spherical integrals}
\author{Vanessa Piccolo}
\address{Unité de Mathématiques Pures et Appliquées (UMPA), ENS Lyon}
\email{vanessa.piccolo@ens-lyon.fr}

\subjclass[2020]{60B20, 60G15, 60F10} 
\keywords{Landscape complexity, Spherical spin glasses, Kac--Rice formula, Multi-spiked GOE matrix, Finite-rank spherical integrals}
\date{\today}

\begin{abstract}
We study the annealed complexity of Gaussian random homogeneous polynomials on the \((N-1)\)-dimensional unit sphere in the presence of deterministic perturbations depending on fixed orthonormal vectors and external parameters. We derive variational formulas for the exponential asymptotics of the average number of critical points and local maxima. Our approach combines the Kac--Rice formula with determinant asymptotics for finite-rank perturbations of Gaussian Wigner matrices. In particular, the determinant analysis builds on recent results by Guionnet and Husson~\cite{MR4436026} on finite-rank spherical integrals, which we use to establish large deviation estimates for the largest eigenvalue of finite-rank Gaussian Wigner matrices. The resulting variational problems reveal a topological phase transition: above an explicit threshold in the external parameters, new zero-complexity regions emerge, corresponding to critical points with large correlation with the perturbation vectors. We also identify regions associated with critical points having large correlations with several vectors simultaneously; numerical evidence suggests that these critical points are more likely to be saddles than local maxima.
\end{abstract}

\maketitle

{
\hypersetup{linkcolor=black}
\tableofcontents
}
\section{Introduction}

In this paper, we study the complexity of finite-rank spiked random polynomials on the \((N-1)\)-dimensional unit sphere \(\mathbb{S}^{N-1} = \{\boldsymbol{x} \in \R^N \colon \lVert \boldsymbol{x} \rVert_2=1 \}\), where \(\lVert \cdot \rVert_2\) is the standard \(2\)-norm. These landscapes consist of an isotropic Gaussian random field perturbed by a deterministic polynomial involving finitely many orthogonal spike directions. While a single spike introduces one preferred direction in the landscape, several spikes lead to competing directions, allowing critical points to correlate with one spike, with several spikes simultaneously, or with none of them. We characterize the exponential growth rates of the numbers of critical points and local maxima in the presence of this finite-rank deterministic structure. Such questions are central in the study of high-dimensional nonconvex landscapes arising in spin glasses, statistical inference, and high-dimensional optimization.

\subsection{Model and main results}

Our main object of study is the finite-rank spiked random polynomial
\begin{equation} \label{eq: function f}
f_N(\boldsymbol{\sigma})= \sum_{i=1}^r \lambda_i \langle \boldsymbol{u}_i, \boldsymbol{\sigma} \rangle ^{k_i} +  H_{N,p}(\boldsymbol{\sigma}),
\end{equation}
where \(\boldsymbol{\sigma} = (\sigma_1, \ldots, \sigma_N) \in \mathbb{S}^{N-1}\). Here, the rank \(r\geq 1\) and the degrees \(k_1,\ldots,k_r\geq 3\) are fixed integers, the vectors \(\boldsymbol{u}_1,\ldots,\boldsymbol{u}_r\in\mathbb{S}^{N-1}\) are deterministic and orthogonal, and \(\lambda_1, \ldots, \lambda_r>0\) are fixed signal strengths. The Gaussian component \(H_{N,p}\) is the random homogeneous polynomial of fixed degree \(p\geq 3\) defined by
\begin{equation} \label{eq: Hamiltonian H_{N,p}}
H_{N,p}(\boldsymbol{\sigma}) = \frac{1}{\sqrt{2N}} \sum_{1 \leq i_1, \ldots, i_p \leq N} W_{i_1, \ldots, i_p} \sigma_{i_1}\cdots \sigma_{i_p}.
\end{equation}
The couplings \(W_{i_1, \ldots, i_p}\) are obtained by symmetrizing i.i.d.\ standard Gaussian random variables:
\[
W_{i_1, \ldots, i_p} \coloneqq \frac{1}{p!} \sum_{\pi \in S_p} G_{i_{\pi(1)}, \ldots, i_{\pi(p)}}, 
\]
where \(S_p\) denotes the symmetric group on \(\{1, \ldots, p\}\) and \((G_{i_1,\ldots,i_p})_{1\leq i_1,\ldots,i_p\leq N}\) are i.i.d.\ standard Gaussian random variables. As a result, \(H_{N,p}\) is a centered Gaussian function on \(\mathbb{S}^{N-1}\) with covariance 
\[
\E \left [H_{N,p}(\boldsymbol{\sigma}) H_{N,p}(\boldsymbol{\tau}) \right] = \frac{\langle \boldsymbol{\sigma},\boldsymbol{\tau} \rangle^p}{2N}.
\]
Here, \(\langle \cdot, \cdot \rangle\) denotes the standard Euclidean inner product. In the physics literature, \(H_{N,p}\) corresponds to the Hamiltonian of the \emph{spherical pure \(p\)-spin model}. Thus, \(f_N\) is a finite-rank deterministic perturbation of the pure spherical \(p\)-spin Hamiltonian.

We characterize the \emph{annealed complexity} of the function \(f_N\), namely the exponential-scale asymptotics of the expected number of its critical points and local maxima. More precisely, for any Borel set \(B \subset \R\), let \(\text{Crt}_N^\textnormal{tot}(B)\) denote the total number of critical points of \(f_N\) whose values lie in \(B\), and let \(\text{Crt}_N^\text{max}(B)\) denote the number of such critical points that are local maxima. We study the large-\(N\) asymptotic behavior of
\[
\frac{1}{N} \log \E \left [\text{Crt}_N^\textnormal{tot}(B) \right ] \quad \text{and} \quad \frac{1}{N} \log \E \left [\text{Crt}_N^\text{max}(B) \right ].
\]
The study of complexity in the pure spherical case was initiated in physics by Crisanti and Sommers~\cite{CS95}, and was later made rigorous by Auffinger, Ben Arous, and \v{C}ern\'{y}~\cite{MR2999295}. Fyodorov~\cite{MR3469265} extended this analysis to include spherical pure spin glasses in the presence of a random magnetic field, which is related to the rank-one case of~\eqref{eq: function f} with \(k = 1\). For \(r=1\) and \(k=p\), the random function \(f_N\) reduces to the \emph{spiked tensor model}, introduced by Richard and Montanari~\cite{MR14} in the context of statistical inference. Its annealed complexity was subsequently studied by Ben Arous, Mei, Montanari, and Nica~\cite{MR4011861}. Further related literature is discussed in Subsection~\ref{related work}. Our work extends these analyses beyond the rank-one setting by allowing multiple orthogonal spikes with arbitrary polynomial degrees.

Our main result is as follows. For fixed parameters \(r,p,k_1,\ldots, k_r\), and \(\lambda_1, \ldots, \lambda_r\), there exists an upper semicontinuous function \(\Sigma^{\textnormal{tot}}\), referred to as the \emph{total complexity function}, such that for any Borel set \(B \subset \R\), 
\begin{equation} \label{eq: LDP}
\begin{split}
\sup_{x \in B^\circ} \Sigma^{\textnormal{tot}} (x) & \leq \liminf_{N \to \infty}\frac{1}{N} \log \E \left [\text{Crt}_N^\textnormal{tot}(B) \right ]  \leq \limsup_{N \to \infty} \frac{1}{N} \log \E \left [\text{Crt}_N^\textnormal{tot}(B) \right ]  \leq  \sup_{x \in \overline{B}} \Sigma^{\textnormal{tot}} (x),
\end{split}
\end{equation}
where $\overline{B}$ and $B^\circ$ denote the closure and interior of $B$, respectively. We obtain an analogous result for local maxima. More generally, our results provide a refined complexity function that also keeps track of the correlations \(\langle\boldsymbol{\sigma},\boldsymbol{u}_i\rangle\) between critical points and the spike directions. The proof is based on the Kac--Rice formula, reducing the expected count of critical points to an integral involving the determinant of a Gaussian Wigner matrix with a finite-rank perturbation. The asymptotic analysis of this determinant presents a significant technical challenge, which we address using recent advances in finite-dimensional spherical integrals due to Guionnet and Husson~\cite{MR4436026}. Further details of the proof strategy are given in Subsection~\ref{outline proofs}. 

A key consequence of the complexity characterization underlying~\eqref{eq: LDP} is the emergence of regions in which the total complexity function vanishes, so that the expected number of critical points grows at most sub-exponentially in \(N\). As the signal strengths \(\lambda_1,\ldots,\lambda_r\) cross sharp thresholds, new zero-complexity regions appear at nontrivial values of the correlations. These regimes include critical-point configurations correlated with a single spike, as well as configurations simultaneously correlated with several spikes. 

\subsection{Motivation}

The main motivation for studying finite-rank perturbations of spherical spin-glass landscapes comes from high-dimensional signal-recovery problems involving several latent signal directions. The deterministic part of \(f_N\) favors all vectors having a large correlation with the distinguished vectors \(\boldsymbol{u}_1, \ldots, \boldsymbol{u}_r \in \mathbb{S}^{N-1}\), while the Gaussian term introduces a highly nonconvex random landscape. The resulting landscape provides a tractable setting in which to investigate how several competing signal directions shape the number, location, and geometry of critical points. 

When \(r=1\) and \(k=p\), this model reduces to the landscape associated with the rank-one spiked tensor problem~\cite{MR4011861}. For general \(r\) and degrees \(k_1,\ldots,k_r\), however, the landscape \(f_N\) is not itself directly associated with an inference problem. Nevertheless, we believe that its complexity provides insight into estimation problems arising naturally in high-dimensional statistics, such as the multi-rank spiked tensor model described below. In this model, the goal is to recover the unknown signal vectors \(\boldsymbol{u}_1, \ldots, \boldsymbol{u}_r \in \mathbb{S}^{N-1}\) from a noisy observation of a tensor \(\boldsymbol{Y} = (Y_{i_1, \ldots, i_p}) \in (\R^N)^{\otimes p}\) given by 
\begin{equation*}
\boldsymbol{Y} = \sum_{i=1}^r \lambda_i \boldsymbol{u}_i^{\otimes p} +\boldsymbol{W},
\end{equation*}
where \(\boldsymbol{W}\in (\R^N)^{\otimes p}\) is a symmetric Gaussian tensor. A natural estimator for the spikes \(\boldsymbol{u}_1, \ldots, \boldsymbol{u}_r\) is the maximum likelihood estimator (MLE). In the rank-one case, the MLE corresponds to the argmax of the random function \(f_N\) with \(k=p\)~\cite{MR14,MR4011861}. For \(r>1\), the likelihood landscape is instead defined on the Stiefel manifold \(\textnormal{St}(N,r) \coloneqq \{ (\boldsymbol{\sigma}_1, \ldots, \boldsymbol{\sigma}_r) \in (\R^N)^r \colon \langle \boldsymbol{\sigma}_i, \boldsymbol{\sigma}_j \rangle = \delta_{ij}\} \), and the MLE is given by the maximizer of the random function \(F_N \colon \textnormal{St}(N,r) \to \R\) defined by
\begin{equation*}
F_N(\boldsymbol{\sigma}_1, \ldots, \boldsymbol{\sigma}_r) = \sum_{i=1}^r \lambda_i f_N(\boldsymbol{\sigma}_i) ,
\end{equation*}
where \(f_N\) corresponds to~\eqref{eq: function f} with \(k_j = p\) for all \(1 \leq j \leq r\). 

Although \(f_N\) is not the full likelihood landscape of the multi-rank problem, it captures some of its central geometric features. In particular, its critical points may correlate with one spike, with several spikes simultaneously, or with none of them. Understanding the complexity of these different configurations provides a first step toward describing the more intricate geometry of multi-component estimation problems. A direct complexity analysis of \(F_N\) is substantially more difficult, particularly for local maxima. Through the Kac--Rice formula, the problem leads to the exponential asymptotics of determinants of large block matrices with dependent Gaussian blocks and full-rank deformations. The analysis of local maxima additionally requires large deviation principles for the extreme eigenvalues of such matrices and involves challenging tools, including the asymptotics of spherical integrals. The complexity of the multi-rank spiked tensor model therefore requires further work and is left for future research.

\subsection{Related work} \label{related work}

As mentioned above, early work on the complexity of high-dimensional landscapes was carried out in statistics and physics in the context of mean-field glasses and spin glasses~\cite{BM80,MR1136952,CS95,CGP98}. In the mathematical literature, a major mathematical breakthrough was Fyodorov’s analysis~\cite{MR2115095} on the model known as the ``zero-dimensional elastic manifold''. Since then, the landscape properties have been studied for many different models. Of interest here are the works on spherical spin glasses~\cite{MR2999295,MR3161473,MR3706746,MR4346481} and in the presence of a rank-one perturbation~\cite{MR3469265,MR4011861,ros2019complex,SBCKZ19,SKUZ19, MR4354698,MR4404771}. The classical approach for counting the total number of such points is the Kac--Rice formula, which relates the average number of critical points to conditional averages of random matrix determinants. In most of the models studied so far, the random matrices appearing in the Kac--Rice formula are related to the Gaussian Orthogonal Ensemble (GOE). Recently,~\cite{MR4552227, MR4673883} studied landscape models with few distributional symmetries whose conditioned Hessian is a random matrix with a noninvariant distribution. We wish to emphasize that we address here the question of the average number of critical points, which in general may not be close to its typical value. A more representative quantity of the typical number of critical points is given by the \emph{quenched complexity}, which corresponds to the large-\(N\) asymptotics of the average of the logarithm of the number of critical points. In most of the models, the random variable \(\text{Crt}^{\textnormal{tot}}_N\) does not concentrate about its expectation, and quenched and annealed complexity differ from each other even at leading order, as in the case of the spiked tensor model~\cite{ros2019complex}. The computation of the quenched complexity is a challenging problem and is usually based on a nonrigorous but exact approach involving the Kac--Rice formula and the replica theory from statistical physics (see e.g.~\cite{ros2019complex,MBB20}). In the case of spherical pure spin glasses, it has been proved by means of a second moment analysis that the number of local minima and of critical points concentrates around its expectation~\cite{MR3706746,MR4346481}. 

In addition to the topological complexity question, models that depend on external parameters also provide a characterization of the landscape by identifying a \emph{topological trivialization transition} when tuning the parameters. Specifically, there is an exact threshold between regimes where the complexity is positive and regimes where the complexity vanishes. Understanding these phase transitions can be useful in predicting the dynamics of optimization on the landscape. A positive complexity indicates a landscape with exponentially many critical points and may therefore signal greater algorithmic difficulty, whereas a sub-exponential number of critical points is often associated with a simpler landscape.

 Optimization in high-dimensional landscapes is actually computationally hard. For the special case \(p = 2\), the model~\eqref{eq: function f} reduces to the \emph{spiked matrix model} introduced in~\cite{MR1863961}. This model exhibits a phase transition, known as the \emph{BBP transition}~\cite{MR2165575}, where there exists an order \(1\) critical threshold \(\lambda_{\textnormal{IT}}\) such that below \(\lambda_{\textnormal{IT}}\), it is information-theoretically impossible to detect the signal vectors, and above \(\lambda_{\textnormal{IT}}\), it is possible to detect the spikes by Principal Component Analysis. The same phenomenon has been observed for the spiked tensor model. In particular, in the high-dimensional asymptotic regime, there is an order \(1\) critical threshold \(\lambda_{\textnormal{IT}}\) (depending on the order \(p\)), below which it is information-theoretically impossible to detect the spikes, and above, the MLE is known to be a good, consistent estimator. Computing the MLE, however, is NP-hard and finding a good algorithm to compute it quickly remains a challenge (this phenomenon is known as a \emph{computational-to-statistical gap}). For the rank-one case, it was shown heuristically in~\cite{MR14} that for the tensor power iteration method, it is possible to recover the spike provided \(\lambda \gtrsim N^\frac{p-2}{2}\). This conjecture was rigorously proved in~\cite{JMLR:v23:21-1290}. The same threshold was obtained for Langevin dynamics and gradient flow in~\cite{MR4124533} and for stochastic gradient descent in~\cite{arous2021online}. Moreover,~\cite{MR4677744} investigated the tensor unfolding algorithm and proved that it is possible to successfully recover the spike provided \(\lambda \gtrsim N^\frac{p-2}{4}\), as also previously conjectured in~\cite{MR14}. The sharp threshold \(\lambda \gtrsim N^\frac{p-2}{4}\) was also achieved using Sum-of-Squares algorithms in~\cite{pmlr-v40-Hopkins15}. 

For a more complete review on the characterization of high-dimensional random landscapes, we direct the reader to~\cite{ros2022highd} and the references therein.

\subsection{Outline of proofs} \label{outline proofs}

In the first part of the paper, we derive the variational formula~\eqref{eq: LDP} by following the strategy developed by Auffinger, Ben Arous, and \v{C}ern\'{y}~\cite{MR2999295} and Auffinger and Ben Arous~\cite{MR3161473} for the Hamiltonians of spherical spin glasses, as well as by Ben Arous, Mei, Montanari, and Nica~\cite{MR4011861} for the spiked tensor model. A standard approach for computing \(\E \left [\text{Crt}^{\textnormal{tot}}_N(B) \right ]\) is through the Kac--Rice formula, which in our setting reads
\[
\E \left [\text{Crt}^{\textnormal{tot}}_N(B) \right] = \int_{\mathbb{S}^{N-1}} \E \left [ |\det \nabla^2  f_N(\boldsymbol{\sigma}) |  \mathbf{1}_{\{f_N(\boldsymbol{\sigma}) \in B\}} \Big |  \nabla f_N(\boldsymbol{\sigma}) = \boldsymbol{0} \right ] \varphi_{ \nabla f_N}(\boldsymbol{0}) \text{d} \boldsymbol{\sigma}.
\]
Here, \( \nabla f_N\) and \( \nabla^2 f_N\) denote the spherical gradient and Hessian of \(f_N\), respectively, and \(\varphi_{ \nabla  f_N}(\boldsymbol{0})\) is the density of \( \nabla f_N\) at \(\boldsymbol{0}\). Conditionally on the gradient, the Hessian appearing in the Kac--Rice formula is distributed as a finite-rank perturbation of a shifted Gaussian Orthogonal Ensemble (GOE) matrix. More precisely, it has the form
\[
\boldsymbol{H}_{N-1} (\boldsymbol{\gamma},t) \stackrel{d}{=} \boldsymbol{W}_{N-1}+ \sum_{i=1}^r \gamma_i \boldsymbol{e}_i \boldsymbol{e}_i^\top - t \boldsymbol{I}_{N-1},
\]
where \(\boldsymbol{W}_{N-1} \sim \textnormal{GOE}(N-1)\), \(\boldsymbol{e}_1, \ldots, \boldsymbol{e}_r\) are the standard basis vectors of \(\R^{N-1}\), and the parameters \(\boldsymbol{\gamma}=(\gamma_1,\ldots,\gamma_r)\) and \(t\) are deterministic functions of the conditioning variables appearing in the Kac--Rice formula. We are therefore led to study the exponential asymptotics of \(\E\left[ \left|\det \boldsymbol{H}_{N-1}(\boldsymbol{\gamma},t)\right| \right]\). For every compact product set \(\mathcal{U} = \mathcal{U}_1 \times \cdots \times \mathcal{U}_r \subset \R^r\) and every compact set \(\mathcal{T} \subset \R\), we prove that
\[
\lim_{N \to \infty} \sup_{\boldsymbol{\gamma} \in \mathcal{U}, t \in \, \mathcal{T}} \left | \frac{1}{N} \log \E \left [ |\det \boldsymbol{H}_{N-1} (\boldsymbol{\gamma},t) | \right ] - \int_\R \log |\lambda - t| \mu_{\textnormal{sc}}( \text{d}\lambda) \right |=0,
\]
where \(\mu_{\textnormal{sc}}\) denotes the semicircle law supported on \([-2,2]\). This conclusion is expected because a finite-rank perturbation does not change the limiting empirical spectral distribution, while the contribution of any outlier eigenvalues vanishes after normalization by \(N\). Making this heuristic rigorous requires additional estimates, which are obtained from Theorem~1.2 of Ben Arous, Bourgade, and McKenna~\cite{MR4552227}.

We also compute the annealed complexity of local maxima. In this case, the Kac--Rice formula contains the additional constraint that the conditioned Hessian be negative definite, and one has to study the large-\(N\) limit of \(\frac{1}{N} \log \E \left [| \det \boldsymbol{H}_{N-1} (\boldsymbol{\gamma},t) | \mathbf{1}_{\{\boldsymbol{H}_{N-1} (\boldsymbol{\gamma},t) \prec 0\}} \right ]\). Thus the main additional difficulty is to understand the asymptotic behavior of 
\[
\frac{1}{N} \log \mathbb{P} \left ( \lambda_{\textnormal{max}} \left (\boldsymbol{W}_{N-1} + \sum_{i=1}^r \gamma_i \boldsymbol{e}_i \boldsymbol{e}_i^\top \right ) \leq t \right).
\]
To this end, we require large deviation principles (LDPs) for the extreme eigenvalues of finite-rank spiked GOE matrices. In the case of a deterministic rank-one perturbation, the LDP for the largest eigenvalue was established by Ma\"{i}da~\cite{MR2336602} and subsequently used in~\cite{MR4011861} to compute the complexity of local maxima in the spiked tensor model. For finite-rank perturbations, the main difficulty is the asymptotics of \emph{finite-rank spherical integrals}, which take the form
\begin{equation} \label{HCIZ}
\E \left [ \exp \left \{\frac{N}{2} \sum_{i=1}^r \gamma_i \langle \boldsymbol{v}_i, \boldsymbol{D}_N \boldsymbol{v}_i \rangle \right \} \right ] ,
\end{equation}
where \(\boldsymbol{D}_N\) is an \(N \times N\) diagonal matrix and the expectation is taken over a Haar-distributed orthonormal family \(\boldsymbol{v}_1, \ldots, \boldsymbol{v}_r\). Equivalently, these vectors may be taken to be the first \(r\) columns of a Haar-distributed orthogonal matrix. The integral~\eqref{HCIZ} is a special case of the orthogonal Harish-Chandra/Itzykson/Zuber (HCIZ) integral, which is defined for two \(N \times N\) matrices \(\boldsymbol{A}_N\) and \(\boldsymbol{B}_N\) by 
\[
HCIZ (\boldsymbol{A}_N, \boldsymbol{B}_N) \coloneqq \int_{\mathcal{O}_N} \exp \left \{N\Tr (\boldsymbol{U} \boldsymbol{A}_N \boldsymbol{U}^\top \boldsymbol{B}_N) \right \} \textnormal{d} m_N (\boldsymbol{U}),
\]
where \(m_N\) denotes the Haar probability measure on the orthogonal group \(\mathcal{O}_N\) of \(N \times N\) orthogonal matrices. By choosing \(\boldsymbol{A}_N = \boldsymbol{D}_N\) and \(\boldsymbol{B}_N = \frac{1}{2} \textnormal{diag} (\gamma_1,\ldots, \gamma_r, 0, \ldots, 0)\), the HCIZ integral reduces to~\eqref{HCIZ}. Guionnet and Husson~\cite{MR4436026} showed that finite-rank spherical integrals are asymptotically equivalent to the product of rank-one spherical integrals. Using this result, they established a LDP for the joint law of the \(k\) largest eigenvalues and the \(\ell\) smallest eigenvalues of a Gaussian Wigner matrix perturbed by a finite-rank matrix having \(k\) nonnegative eigenvalues and \(\ell\) nonpositive eigenvalues (see~\cite[Proposition 2.7]{MR4436026}). By combining this LDP with classical techniques, we derive a LDP for the largest eigenvalue of the finite-rank spiked GOE, thereby extending the rank-one result of~\cite{MR2336602}. Our variational formula additionally records the scalar products with each of the distinguished directions \(\boldsymbol{u}_1,\ldots,\boldsymbol{u}_r\), and therefore yields a more refined description of the complexity landscape than one based only on the critical values.

In the second part of the paper, we analyze the variational problem~\eqref{eq: LDP} and identify the regions in which the complexity function vanishes and the expected number of critical points is therefore sub-exponential in \(N\). In particular, we identify a topological phase transition characterized by a sharp threshold in the parameters \(\lambda_1, \ldots, \lambda_r\). Beyond this threshold, new zero-complexity regions emerge in which critical points have large scalar products with \(\boldsymbol{u}_1, \ldots, \boldsymbol{u}_r\). Some of these regions contain critical points that are strongly correlated with more than one distinguished direction. However, numerical evidence for \(r=2\) suggests that local maxima with large correlations occur in regions where the critical points are closely aligned with a single spike. This extends the qualitative picture obtained in~\cite{MR4011861} for a single spike to the multi-spike setting.

\subsection{Overview}
The paper is organized as follows. In Section~\ref{main results}, we state our main results on the landscape complexity of the random function \(f_N\). Section~\ref{background} provides key intermediate results, including the computation of the random matrix arising from the Kac--Rice formula and the derivation of the large deviation principle for the largest eigenvalue of a finite-rank spiked GOE matrix. The proofs of the main theorems are presented in Section~\ref{proof main results}. Finally, in Section~\ref{analysis complexity function}, we analyze the variational problem associated with the total complexity function and determine the critical threshold that characterizes the topological phase transition.\\

\textbf{Acknowledgements.} I am deeply grateful to my PhD advisors, Prof.\ Alice Guionnet and Prof.\ Gérard Ben Arous, for proposing this project and for their invaluable guidance, encouragement, and insightful discussions throughout the research process. I also wish to thank Ben McKenna and Justin Ko for their helpful comments. This work was supported by the ERC Advanced Grant LDRaM (No.\ 884584). 
\section{Main results}\label{main results}
 
Our main results concern the exponential asymptotics of the expected number of critical points and local maxima of the function \(f_N\) introduced in~\eqref{eq: function f}. Throughout, let \(\nabla f_N(\boldsymbol{\sigma})\) and \(\nabla^2 f_N(\boldsymbol{\sigma})\) denote, respectively, the Riemannian gradient and Hessian at \(\boldsymbol{\sigma} \in \mathbb{S}^{N-1}\), computed with respect to the standard Riemannian metric on \(\mathbb{S}^{N-1}\). 

\begin{defn} \label{def: random variables}
Let \(M_1, \ldots, M_r \subset [-1,1]\) and \(B \subset \R\) be Borel sets, and write \(M=M_1 \times \cdots \times M_r\). We define the random variable \(\textnormal{Crt}_N^{\textnormal{tot}} (M, B)\) counting the critical points of \(f_N\) whose scalar products with \(\boldsymbol{u}_i\) lie in \(M_i\) and whose values lie in \(B\) by
\begin{equation} \label{def: tot crit points}
\textnormal{Crt}_N^{\textnormal{tot}} (M, B) \coloneqq \sum_{\boldsymbol{\sigma} \in \mathbb{S}^{N-1}\colon  \nabla f_N(\boldsymbol{\sigma})=0} \mathbf{1}_{\{(\langle \boldsymbol{\sigma} , \boldsymbol{u}_i \rangle)_{i=1}^r \in M\}} \mathbf{1}_{\{f_N(\boldsymbol{\sigma}) \in B\}}.
\end{equation}
We similarly define the random variable counting local maxima by
\begin{equation} \label{def: crit points index}
\textnormal{Crt}_N^{\textnormal{max}} (M, B) \coloneqq \sum_{\boldsymbol{\sigma} \in \mathbb{S}^{N-1} \colon \nabla f_N(\boldsymbol{\sigma})=0}  \mathbf{1}_{\{(\langle \boldsymbol{\sigma} , \boldsymbol{u}_i \rangle)_{i=1}^r \in M\}} \mathbf{1}_{\{f_N(\boldsymbol{\sigma}) \in B \}} \mathbf{1}_{\{ \nabla^2 f_N(\boldsymbol{\sigma}) \prec 0 \}},
\end{equation}
where \(\nabla^2 f_N(\boldsymbol{\sigma}) \prec 0\) denotes that the Riemannian Hessian is negative definite on the tangent space \(T_{\boldsymbol{\sigma}} \mathbb{S}^{N-1}\), or equivalently, that all of its \(N-1\) eigenvalues are strictly negative.
\end{defn}

\subsection{Annealed complexity of total critical points} \label{annealed complexity}

We first state the annealed complexity result for all critical points. To this end, we introduce the total complexity function. We begin with the logarithmic potential of the semicircle law.

\begin{defn} \label{def: log-potential}
The \emph{\(\log\)-potential of a measure \(\mu\)} is defined by
\[
\Omega_\mu(x) \coloneqq \int_\R \log |x-\lambda| \textnormal{d} \mu(\lambda).
\]
Let \(\mu_{\textnormal{sc}}\) denote the semicircle measure, whose density with respect to Lebesgue measure is 
\[
\mu_{\textnormal{sc}} (\textnormal{d} x)= \frac{1}{2 \pi} \sqrt{4-x^2} \mathbf{1}_{\{|x| \le 2\}} \textnormal{d}x.
\]
We denote its \(\log\)-potential by \(\Omega(x) \coloneqq \Omega_{\mu_{\textnormal{sc}}}(x)\). The explicit expression for \(\Omega(x)\) is given by (see e.g.~\cite[eq.\ (2.4)]{MR3706746})
\begin{equation} \label{def: function Phi star}
\Omega(x) = \begin{cases}
\frac{x^2}{4}-\frac{1}{2} & \text{if} \enspace 0 \le |x|\le 2,\\
\frac{x^2}{4}-\frac{1}{2} -\left[ \frac{|x|}{4} \sqrt{x^2-4} - \log \left ( \frac{|x| + \sqrt{x^2-4}}{2}\right )\right] & \text{if} \enspace |x| > 2.
\end{cases}
\end{equation}
\end{defn}

The total complexity function is then defined as follows.

\begin{defn} \label{def: function Sigma tot}
Define the admissible domain by
\[ 
\mathcal{D}_r \coloneqq \left\{ \boldsymbol{m} \in [-1,1]^r \colon \| \boldsymbol{m}\|_2 < 1 \right\}. 
\] 
The \emph{total complexity function} \(\Sigma^{\textnormal{tot}} \colon [-1, 1]^r \times \R \to \R \cup \{-\infty\}\) is defined by
\begin{equation*} 
\Sigma^{\textnormal{tot}}(\boldsymbol{m}, x) = 
\begin{cases}
\Sigma(\boldsymbol{m}, y(\boldsymbol{m},x)) & \text{if} \enspace \boldsymbol{m} \in \mathcal{D}_r ,  \\
-\infty & \text{if} \enspace \boldsymbol{m} \notin \mathcal{D}_r,
\end{cases}
\end{equation*}
where
\begin{equation*}
y(\boldsymbol{m},x) = x - \sum_{i=1}^r \lambda_i \left ( 1 - \frac{k_i}{p} \right ) m_i^{k_i},
\end{equation*}
and for \((\boldsymbol{m},y) \in \mathcal{D}_r \times \R\), \(\Sigma(\boldsymbol{m},y) \) is given by
\begin{equation*} 
\begin{split} 
\Sigma(\boldsymbol{m},y) & = \frac{1}{2} (\log(p-1) +1) + \frac{1}{2} \log \left (1- \| \boldsymbol{m}\|_2^2 \right) - \frac{1}{p} \sum_{i=1}^r \lambda_i^2 k_i^2 m_i^{2k_i-2}\\
& \quad +\frac{1}{p} \left (\sum_{i=1}^r \lambda_i k_i m_i^{k_i} \right)^2 - \left ( y- \frac{1}{p}\sum_{i=1}^r \lambda_i k_i m_i^{k_i}\right )^2 + \Omega \left ( \sqrt{\frac{2p}{p-1}} y \right ).
\end{split}
\end{equation*}
\end{defn}

The main result on the annealed complexity of total critical points is the following.

\begin{thm}\label{thm crit points}
Fix $p$, $r$, $k_1,\ldots,k_r$, and $\lambda_1,\ldots,\lambda_r$. Then, for all Borel sets \(M_1, \ldots, M_r \subset [-1,1]\) and \(B \subset \R\), it holds that
\begin{align}
\limsup_{N \to \infty} \frac{1}{N} \log \E \left [ \textnormal{Crt}_N^{\textnormal{tot}}(M,B) \right ]& \leq \sup_{(\boldsymbol{m},x) \in \overline{M} \times \overline{B} } \Sigma^{\textnormal{tot}}(\boldsymbol{m}, x), \label{sup crit points}\\
\liminf_{N \to \infty} \frac{1}{N} \log \E \left [ \textnormal{Crt}_N^{\textnormal{tot}}(M,B) \right ]& \geq \sup_{(\boldsymbol{m},x) \in M^\circ \times B^\circ } \Sigma^{\textnormal{tot}}(\boldsymbol{m}, x). \label{inf crit points}
\end{align}
Here and throughout, we use the convention \(\sup\emptyset=-\infty\).
\end{thm}

As mentioned in the introduction, when \(r=1\), the objective function \(f_N\) reduces to
\[
f_N(\boldsymbol{\sigma})= \lambda \langle \boldsymbol{\sigma},\boldsymbol{u} \rangle ^k + \frac{1}{\sqrt{2N}} \sum_{1 \leq i_1, \ldots, i_p \leq N} W_{i_1, \ldots, i_p} \sigma_{i_1}\cdots \sigma_{i_p}.
\]
This model is known as the \((p,k)\) spiked tensor model and was introduced in~\cite{MR4404771} to study the sharp asymptotics for the average number of local maxima. When \(k=p\), the argmax of \(f_N\) corresponds to the MLE for the signal vector \(\boldsymbol{u}\) and the complexity of this model was studied in~\cite{MR4011861}. In particular, Theorem~\ref{thm crit points} reduces to~\cite[Theorem 2.1]{MR4011861} in the special case \(r=1\) and \(k=p\).

\subsection{Annealed complexity of local maxima} 

We now state the annealed complexity result for local maxima. In contrast to the complexity of total critical points, the negative-definiteness constraint on the Hessian gives rise to an additional spectral contribution. To describe this contribution, we associate with each correlation vector \(\boldsymbol{m}\in\mathcal{D}_r\) an \(r\times r\) symmetric matrix and denote its ordered eigenvalues by \(\gamma_1(\boldsymbol{m}),\ldots,\gamma_r(\boldsymbol{m})\). In Section~\ref{background}, we show how these eigenvalues arise from the finite-rank deformation of the Hessian in the Kac--Rice formula.

\begin{defn}  \label{def: perturbation matrix}
For $\boldsymbol{m}\in\mathcal{D}_r$, let \(\boldsymbol{G}(\boldsymbol{m}) \in \R^{r \times r}\) denote the matrix with entries 
\begin{equation} \label{overlap v}
G_{ij}(\boldsymbol{m}) = \frac{\delta_{ij} - m_im_j}{\sqrt{1-m_i^2} \sqrt{1-m_j^2}} ,
\end{equation}
and define
\[
\boldsymbol{\Theta}(\boldsymbol{m}) = \operatorname{diag} \bigl( \theta_1(m_1),\ldots,\theta_r(m_r) \bigr),
\]
where each \(\theta_i \colon [-1,1] \to \R\) is given by
\begin{equation} \label{function theta}
\theta_i(m_i) = \sqrt{\frac{2}{p(p-1)}} k_i (k_i-1) \lambda_i m_i^{k_i-2}(1-m_i^2).
\end{equation}
Since \(\boldsymbol{G}(\boldsymbol{m})\) is symmetric positive definite on \(\mathcal{D}_r\), the matrix \(\boldsymbol{G}(\boldsymbol{m})^{1/2} \boldsymbol{\Theta}(\boldsymbol{m}) \boldsymbol{G}(\boldsymbol{m})^{1/2}\) is well-defined and symmetric. We denote its ordered eigenvalues by \(\gamma_1(\boldsymbol{m}) \geq\cdots\geq \gamma_r(\boldsymbol{m}) \).
\end{defn}

We next introduce the function that encodes the additional spectral contribution associated with the negative-definiteness constraint.

\begin{defn} \label{def: function L}
Let 
\[ 
\mathcal{W}_r \coloneqq \left\{ \boldsymbol{\gamma} = (\gamma_1,\ldots,\gamma_r)\in \mathbb{R}^r \colon  \gamma_1\geq\cdots\geq\gamma_r \right\} 
\] 
denote the set of nonincreasing vectors in $\mathbb{R}^r$. We define the function \( L\colon \mathcal{W}_r\times\mathbb{R} \to \mathbb{R}_+\cup\{+\infty\} \) by
\[
L(\boldsymbol{\gamma},t) =
\begin{cases}
\sum\limits_{\substack{1\le i\le r\colon\\ \gamma_i>1,t<\gamma_i+\gamma_i^{-1}}} I_{\gamma_i}(t) & \text{if} \enspace t \ge  2,\\
+ \infty & \text{if} \enspace t < 2,
\end{cases}
\]
where, as usual, an empty sum is equal to zero. Here, for \(\gamma >1\), \(I_\gamma \colon [2, + \infty) \to \R\) is given by 
\begin{equation} \label{eq: I gamma}
I_{\gamma}(x) = \frac{1}{4} \int^x_{\gamma + \frac{1}{\gamma}} \sqrt{y^2-4} \, \textnormal{d}y - \frac{1}{2} \gamma \left ( x - \left ( \gamma + \frac{1}{\gamma}\right ) \right ) + \frac{1}{8} \left ( x^2 - \left (\gamma + \frac{1}{\gamma} \right )^2 \right ).
\end{equation}
\end{defn}

As shown in Section~\ref{background}, \(L(\boldsymbol{\gamma},t)\) gives the exponential cost of forcing the largest eigenvalue of the corresponding finite-rank deformed GOE matrix below the threshold \(t\). The complexity of local maxima is therefore obtained by subtracting this spectral cost from the total complexity.

\begin{defn} \label{def: function Sigma max}
The local-maxima complexity function \(\Sigma^{\textnormal{max}}\colon [-1, 1]^r \times \R \to \R \cup \{-\infty\}\) is defined by
\begin{equation*}
\Sigma^{\textnormal{max}}(\boldsymbol{m}, x) =
\begin{cases}
 \Sigma^{\textnormal{tot}}(\boldsymbol{m}, x)- L \left (\boldsymbol{\gamma}(\boldsymbol{m}), t(\boldsymbol{m},x) \right ) & \text{if} \enspace \boldsymbol{m} \in \mathcal{D}_r ,  \\
-\infty & \text{if} \enspace \boldsymbol{m} \notin \mathcal{D}_r,
\end{cases}
\end{equation*}
where \(\Sigma^{\textnormal{tot}}\) and \(L\) are given in Definitions~\ref{def: function Sigma tot} and~\ref{def: function L}, respectively. Moreover, \(\boldsymbol{\gamma}(\boldsymbol{m}) = (\gamma_1(\boldsymbol{m}), \ldots, \gamma_r(\boldsymbol{m})) \in \mathcal{W}_r\) is given in Definition~\ref{def: perturbation matrix}, and \(t(\boldsymbol{m},x)\) is defined by
\begin{equation}\label{function t}
t(\boldsymbol{m},x) \coloneqq \sqrt{\frac{2p}{p-1}} y (\boldsymbol{m},x) =  \sqrt{\frac{2p}{p-1}} \left ( x - \sum_{i=1}^r \lambda_i \left (1-\frac{k_i}{p} \right )m_i^{k_i} \right).
\end{equation}
Here, we use the convention $a-(+\infty)=-\infty$ for every $a\in\mathbb{R}$. 
\end{defn}

We can now state the main result on the annealed complexity of local maxima.

\begin{thm} \label{thm local maxima}
Fix $p$, $r$, $k_1,\ldots,k_r$, and $\lambda_1,\ldots,\lambda_r$. Then, for all Borel sets \(M_1, \ldots, M_r \subset [-1,1]\) and \(B \subset \R\), it holds that
\begin{align}
\limsup_{N \to \infty} \frac{1}{N} \log \E \left [ \textnormal{Crt}_N^{\textnormal{max}}(M, B) \right ]& \leq \sup_{(\boldsymbol{m},x) \in \overline{M} \times \overline{B} } \Sigma^{\textnormal{max}}(\boldsymbol{m}, x), \label{sup loc maxima}\\
\liminf_{N \to \infty} \frac{1}{N} \log \E \left [ \textnormal{Crt}_N^{\textnormal{max}}(M, B) \right]& \geq \sup_{(\boldsymbol{m},x) \in M^\circ \times B^\circ } \Sigma^{\textnormal{max}}(\boldsymbol{m}, x). \label{inf loc maxima}
\end{align}
\end{thm}

As for total critical points, Theorem~\ref{thm local maxima} reduces to~\cite[Theorem 2.2]{MR4011861} in the special case \(r=1\) and \(k=p\).

\subsection{Analyzing the variational formulas}  \label{analyzing the variational formulas}

In this subsection, we study the variational problems arising in Theorems~\ref{thm crit points} and~\ref{thm local maxima}. Specifically, we analyze the complexity functions
\[
\Sigma^{\textnormal{tot}}(\boldsymbol{m}) \coloneqq \max_{x \in \R} \Sigma^{\textnormal{tot}}(\boldsymbol{m},x) \quad \textnormal{and} \quad \Sigma^{\textnormal{max}}(\boldsymbol{m}) \coloneqq \max_{x \in \R} \Sigma^{\textnormal{max}}(\boldsymbol{m},x),
\]
which describe the exponential growth rates of the number of critical points and local maxima, respectively, as functions of the scalar products \(\langle \boldsymbol{\sigma},\boldsymbol{u}_i \rangle = m_i\) for \(1 \leq i \leq r\). Throughout this subsection, we restrict our attention to the nonnegative admissible domain 
\[
\mathcal{D}_r^+ \coloneqq \mathcal{D}_r\cap[0,1]^r = \left\{ \boldsymbol{m}\in[0,1]^r \colon \|\boldsymbol{m}\|_2<1 \right\}.
\]

\begin{figure}[htbp] 
\begin{subfigure}[b]{0.8\textwidth}
\centering
\includegraphics[scale=0.75]{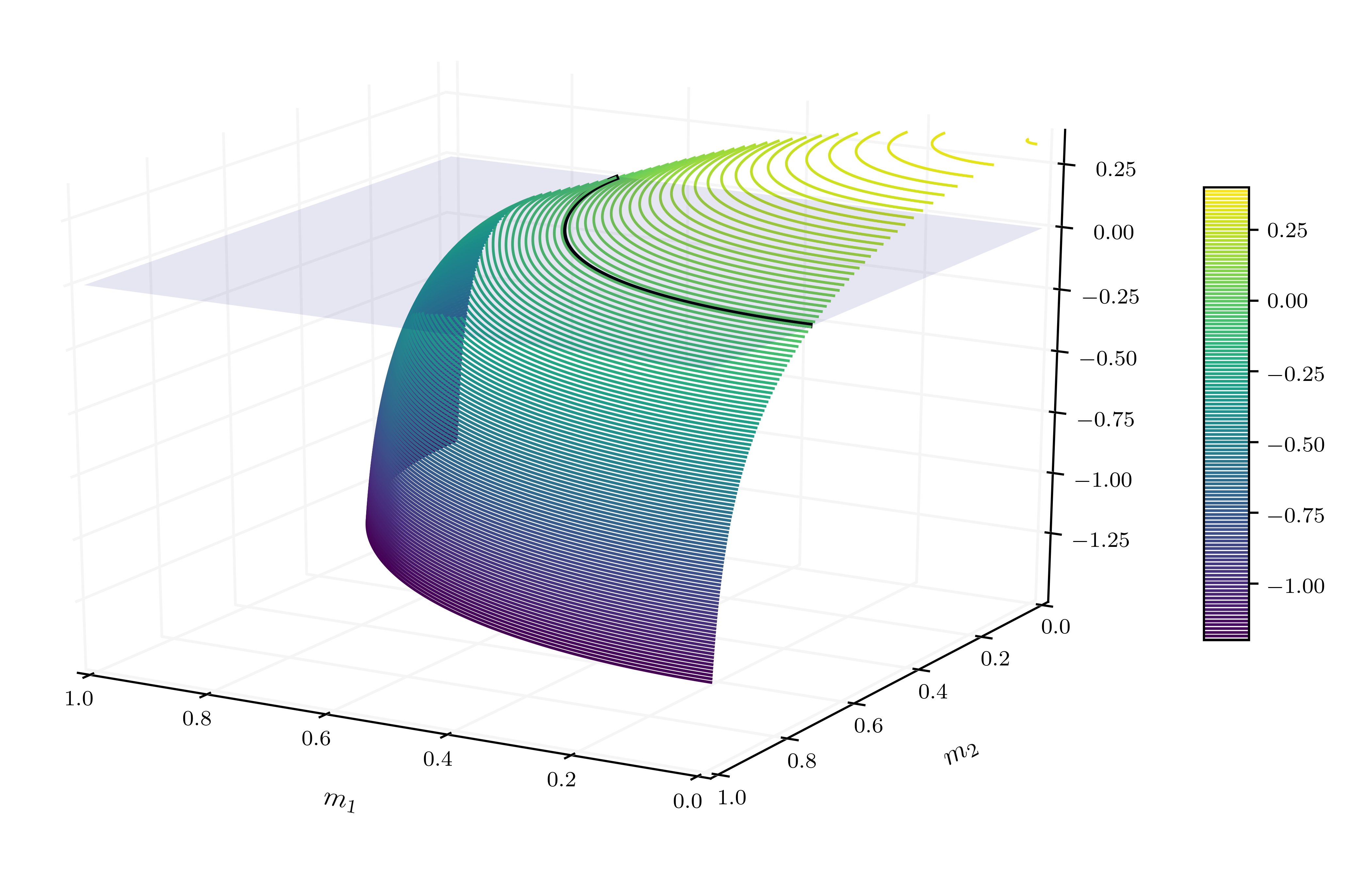}
\caption{\(\lambda_1 = 0.5\) and \(\lambda_2 = 0.2\).}
\label{3D complexity - small}
\end{subfigure}
\vfill
\begin{subfigure}[b]{0.8\textwidth}
\centering
\includegraphics[scale=0.75]{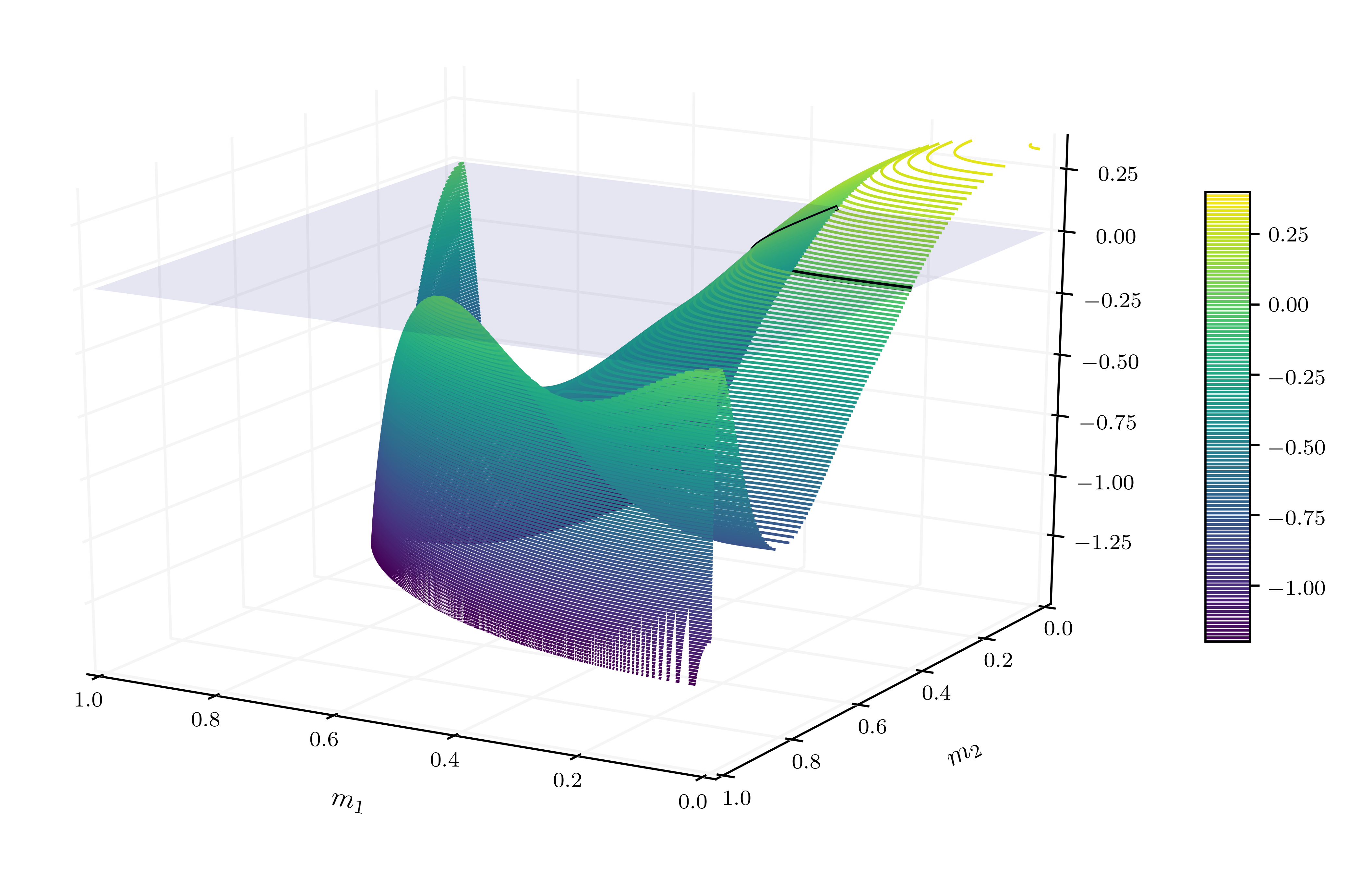}
\caption{\(\lambda_1 = 3\) and \(\lambda_2 = 2\).}
\label{3D complexity - large}
\end{subfigure}
\caption{3D plots of \(\Sigma^{\textnormal{tot}}(\boldsymbol{m})\) for \(r=2\) and \(p=k=3\), with different values of \(\lambda_1\) and \(\lambda_2\). The black curve delineates the boundary where \(\Sigma^{\textnormal{tot}}(\boldsymbol{m})\) transitions from positive to nonpositive complexity. In (a), the complexity is positive for small \(m_1\) and \(m_2\), zero along the black curve, and \emph{negative} for larger values. Thus, critical points with small scalar products are exponentially numerous, while the black curve marks the transition to a region where the expected number is exponentially small. In (b), the region of positive complexity shrinks, and the function becomes \emph{nonpositive} over a larger domain. In particular, \(\Sigma^{\textnormal{tot}}(\boldsymbol{m})\) now vanishes in three additional regions: near \((m_1,m_2) \approx (1,0)\), near \((m_1,m_2) \approx (0,1)\), and for large values where both \(m_1>0\) and \(m_2>0\). This suggests the emergence, at the annealed exponential scale, of regions corresponding to critical points highly correlated either with a single spike or with multiple spikes simultaneously.}
\label{3D complexity}
\end{figure}

We begin by analyzing the variational formula for the total complexity. In Section~\ref{analysis complexity function}, and in particular Proposition~\ref{tot complexity}, we characterize the behavior of the function \(\Sigma^{\textnormal{tot}}(\boldsymbol{m})\) by identifying the regions in \(\boldsymbol{m}\) where it is positive, zero, or negative. Additionally, we study how these regions evolve as the parameters \(\lambda_1, \ldots, \lambda_r >0\) are varied. Broadly speaking, when \(\lambda_1, \ldots, \lambda_r\) are small, the complexity function \(\Sigma^{\textnormal{tot}}(\boldsymbol{m})\) is positive over a relatively large neighborhood around zero and is negative elsewhere (see Figure~\ref{3D complexity - small}). In contrast, when \(\lambda_1, \ldots, \lambda_r\) are sufficiently large and exceed a certain critical threshold, the positive region shrinks to a smaller neighborhood near zero. Outside this neighborhood, the function is nonpositive; moreover, additional zero-complexity regions appear in some regions with large correlation with one or more spikes (see Figure~\ref{3D complexity - large}). 

The following theorem focuses on the regime where the entries of \(\boldsymbol{m}\) are sufficiently large---a condition made precise below---and rigorously formalizes the associated phase transition. 

\begin{thm}\label{thm variational problem crt points}
Assume that \(k_i = k \ge 3\) for every \(1 \leq i \leq r\). For \(\boldsymbol{m} \in \mathcal{D}_r^+\), define 
\[
\tau(\boldsymbol{m}) \coloneqq \frac{k}{p}\sum_{i=1}^r \lambda_i m_i^k, \quad \eta(\boldsymbol{m}) \coloneqq \sum_{i=1}^r \lambda_i^{-\frac{2}{k-2}}\mathbf{1}_{\{m_i \neq 0\}},
\]
and let
\[
\tau_\textnormal{c} (p) \coloneqq \frac{p-2}{\sqrt{2p(p-1)}}, \quad \eta_\textnormal{c}(p,k) \coloneqq \begin{cases} (k-2) \left( \frac{2k^2}{p(k-1)^{k-1}} \right)^{\frac{1}{k-2}} & \textnormal{if} \enspace k \ge p,  \\[1ex] (p-2) \left( \frac{2k^2}{p(p-1)^{k-1}} \right)^{\frac{1}{k-2}} & \textnormal{if} \enspace k \le p. \end{cases}
\]
Then, the complexity function \(\Sigma^{\textnormal{tot}}(\boldsymbol{m})\) is nonpositive for every \(\boldsymbol{m} \in \mathcal{D}_r^+\) such that \(\tau (\boldsymbol{m}) \ge \tau_\textnormal{c}(p)\), and exhibits a phase transition in \(\eta(\boldsymbol{m})\):
\begin{enumerate}
\item If \(\eta(\boldsymbol{m}) > \eta_\textnormal{c}(p,k)\), then \(\Sigma^{\textnormal{tot}}(\boldsymbol{m}) < 0\) for every \(\boldsymbol{m} \in \mathcal{D}_r^+\) such that \(\tau (\boldsymbol{m}) \ge \tau_\textnormal{c}(p)\).
\item If \(\eta(\boldsymbol{m}) \leq \eta_\textnormal{c} (p,k)\), then \(\Sigma^{\textnormal{tot}}(\boldsymbol{m}) \leq 0\) for every \(\boldsymbol{m} \in \mathcal{D}_r^+\) such that \(\tau (\boldsymbol{m}) \ge \tau_\textnormal{c}(p)\), with equality if and only if 
\[
\lambda_i m_i^{k-2} = \lambda_j m_j^{k-2},
\]
for all \(i,j\) such that \(m_i,m_j \neq 0\), and
\[
\tau(\boldsymbol{m}) = \frac{1}{\sqrt{2p}} \frac{\|\boldsymbol{m}\|_2^2}{\sqrt{1 - \|\boldsymbol{m}\|_2^2}}.
\]
\end{enumerate}
\end{thm}

Thus, in the multi-spike setting, the relevant threshold depends not only on the spike strengths but also on which spike directions have nonzero correlation with the critical point. Theorem~\ref{thm variational problem crt points} focuses on the regime \(\tau(\boldsymbol{m}) \geq \tau_\textnormal{c}(p)\), where the complexity function \(\Sigma^{\textnormal{tot}}(\boldsymbol{m})\) is nonpositive and undergoes a phase transition. This regime of larger values of \(\boldsymbol{m}\) is particularly relevant for identifying critical points that are highly correlated with the spike directions. The complementary regime \(\tau(\boldsymbol{m}) < \tau_\textnormal{c}(p)\) is addressed in Section~\ref{analysis complexity function}, where Proposition~\ref{tot complexity} provides a detailed characterization of the regions where \(\Sigma^{\textnormal{tot}}(\boldsymbol{m})\) is positive, zero, or negative. The proof of Theorem~\ref{thm variational problem crt points} is deferred to Section~\ref{analysis complexity function}. 

\begin{figure}[h!] 
\centering
\begin{subfigure}[b]{0.45\textwidth}
\centering
\includegraphics[scale=0.42]{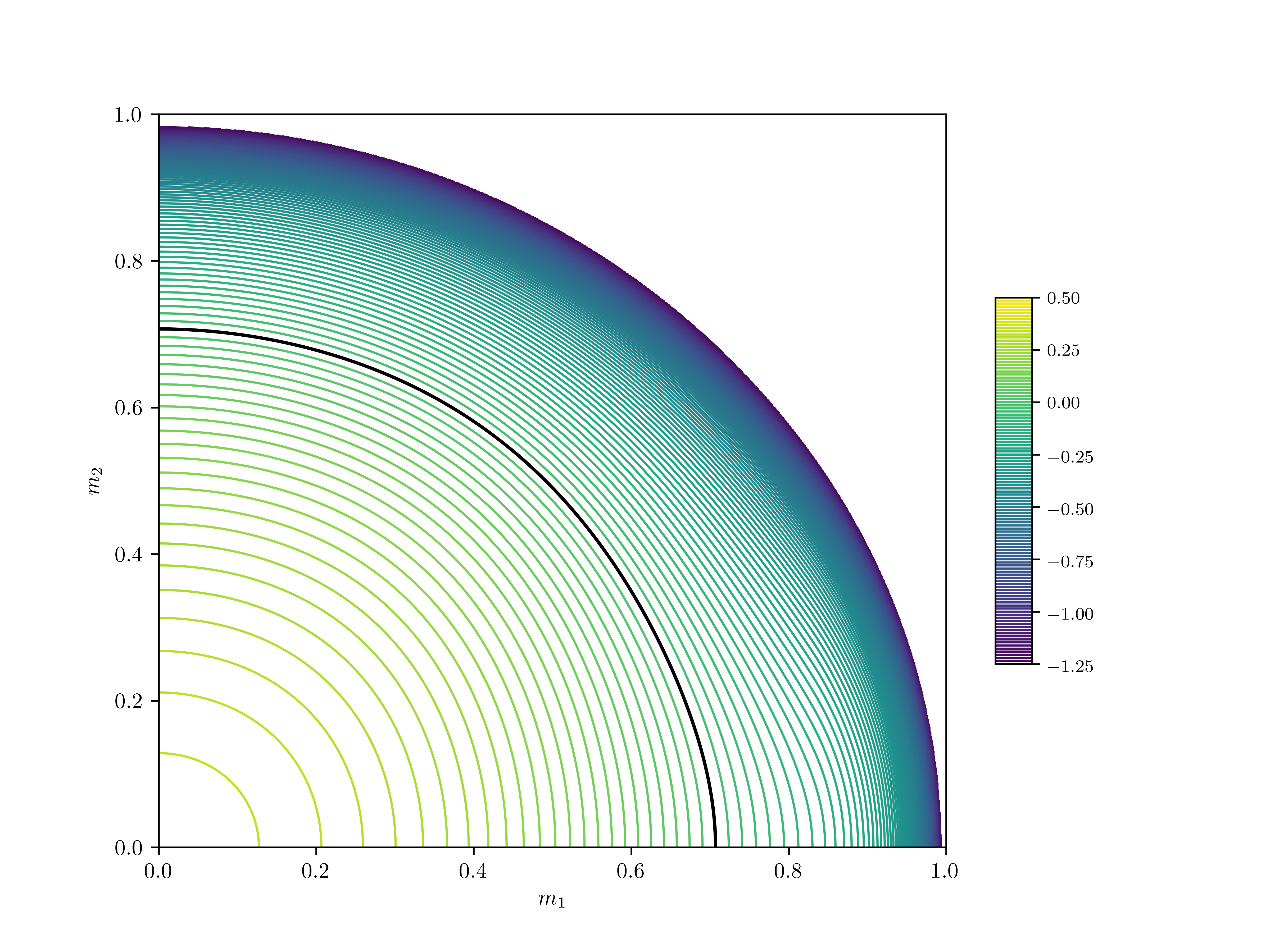}
\caption{\(\lambda_1 = 0.5\) and \(\lambda_2 = 0.2\).}
\label{contour1}
\end{subfigure}
\hspace{0.06\textwidth} 
\begin{subfigure}[b]{0.45\textwidth}
\centering
\includegraphics[scale=0.42]{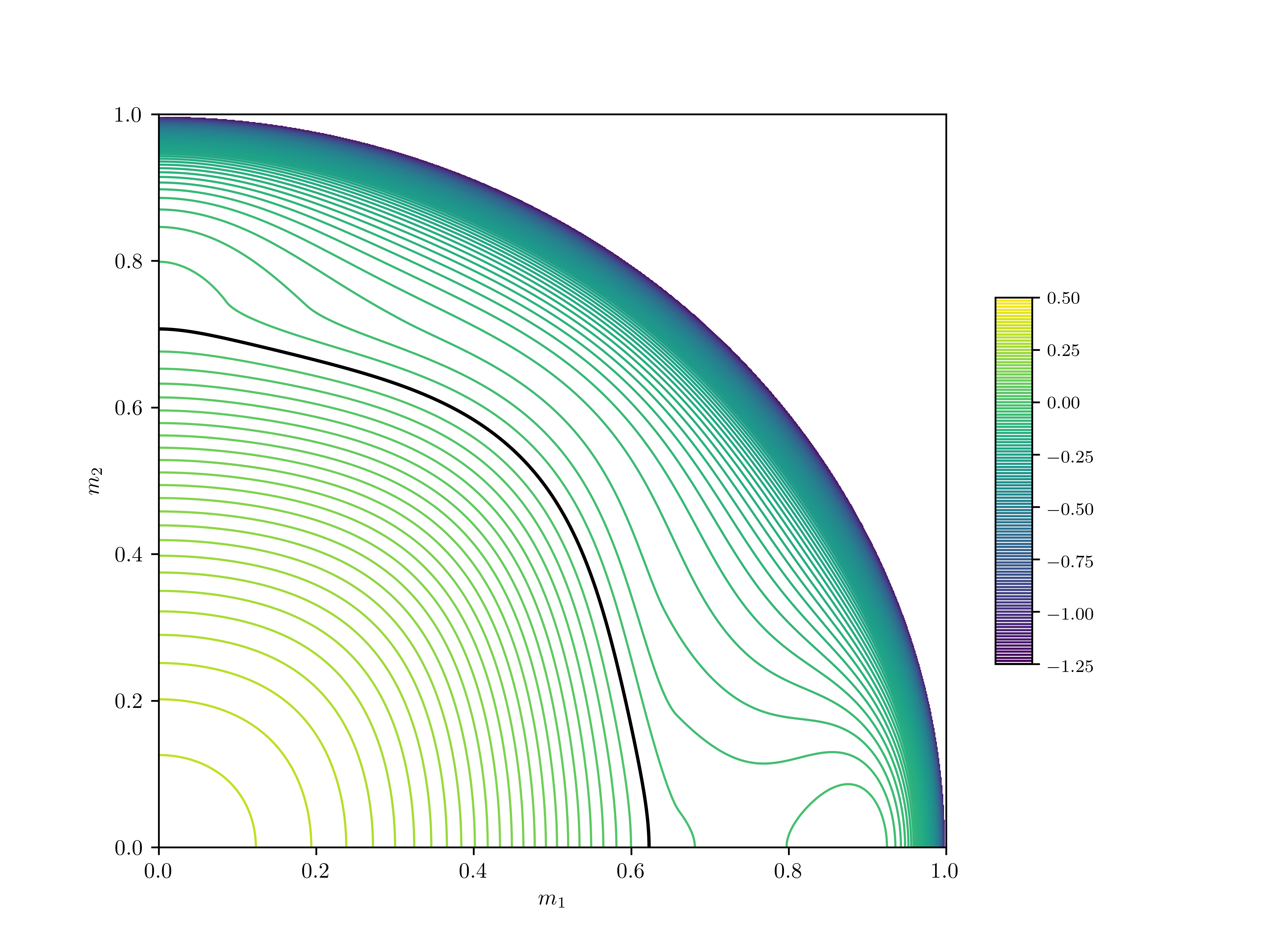}
\caption{\(\lambda_1 = 1\) and \(\lambda_2 = 0.7\).}
\label{contour2}
\end{subfigure}

\vspace{0.5cm} 

\begin{subfigure}[b]{0.45\textwidth}
\centering
\includegraphics[scale=0.42]{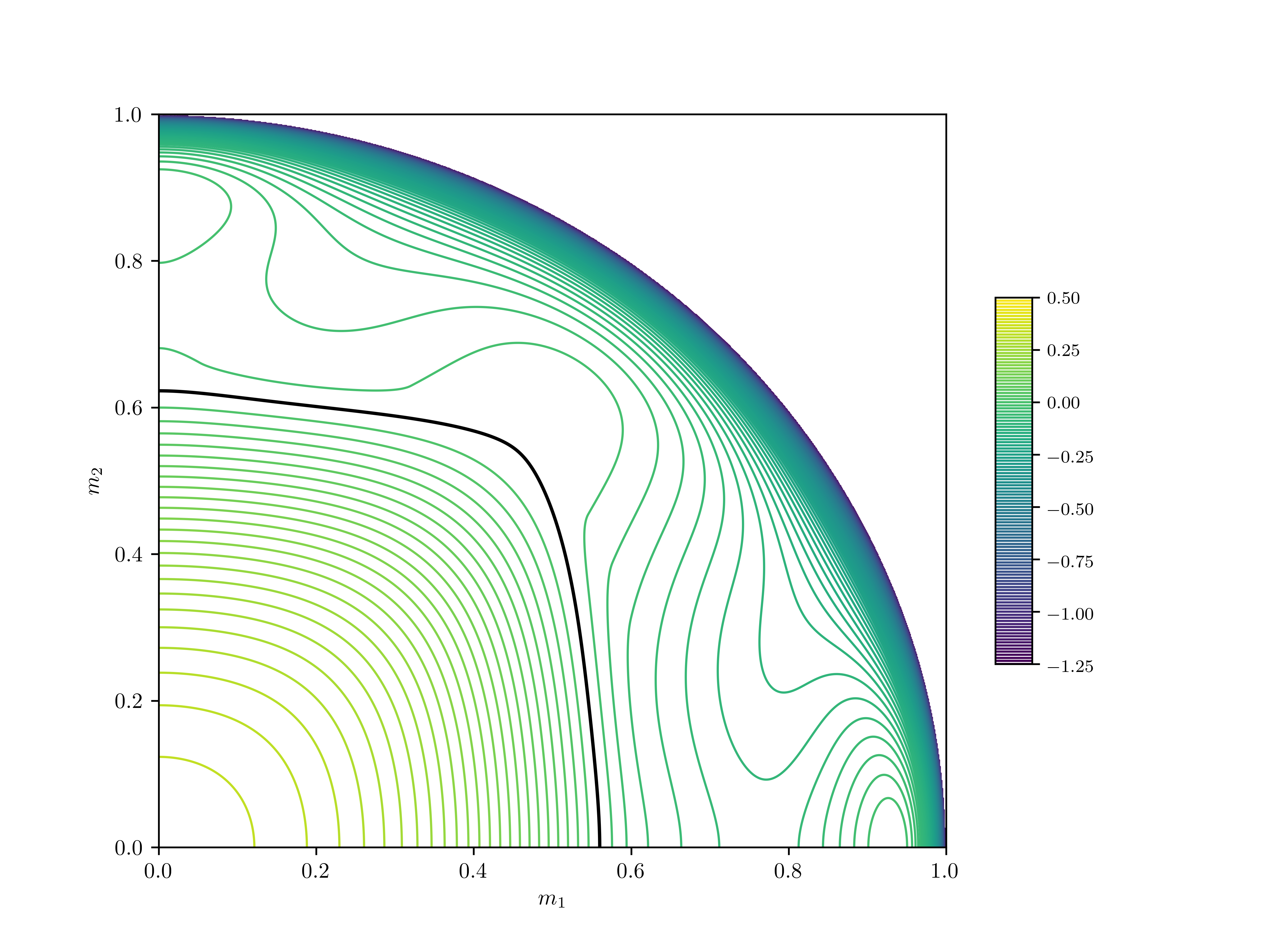}
\caption{\(\lambda_1 = 1.2\) and \(\lambda_2 = 1\).}
\label{contour3}
\end{subfigure}
\hspace{0.06\textwidth} 
\begin{subfigure}[b]{0.45\textwidth}
\centering
\includegraphics[scale=0.42]{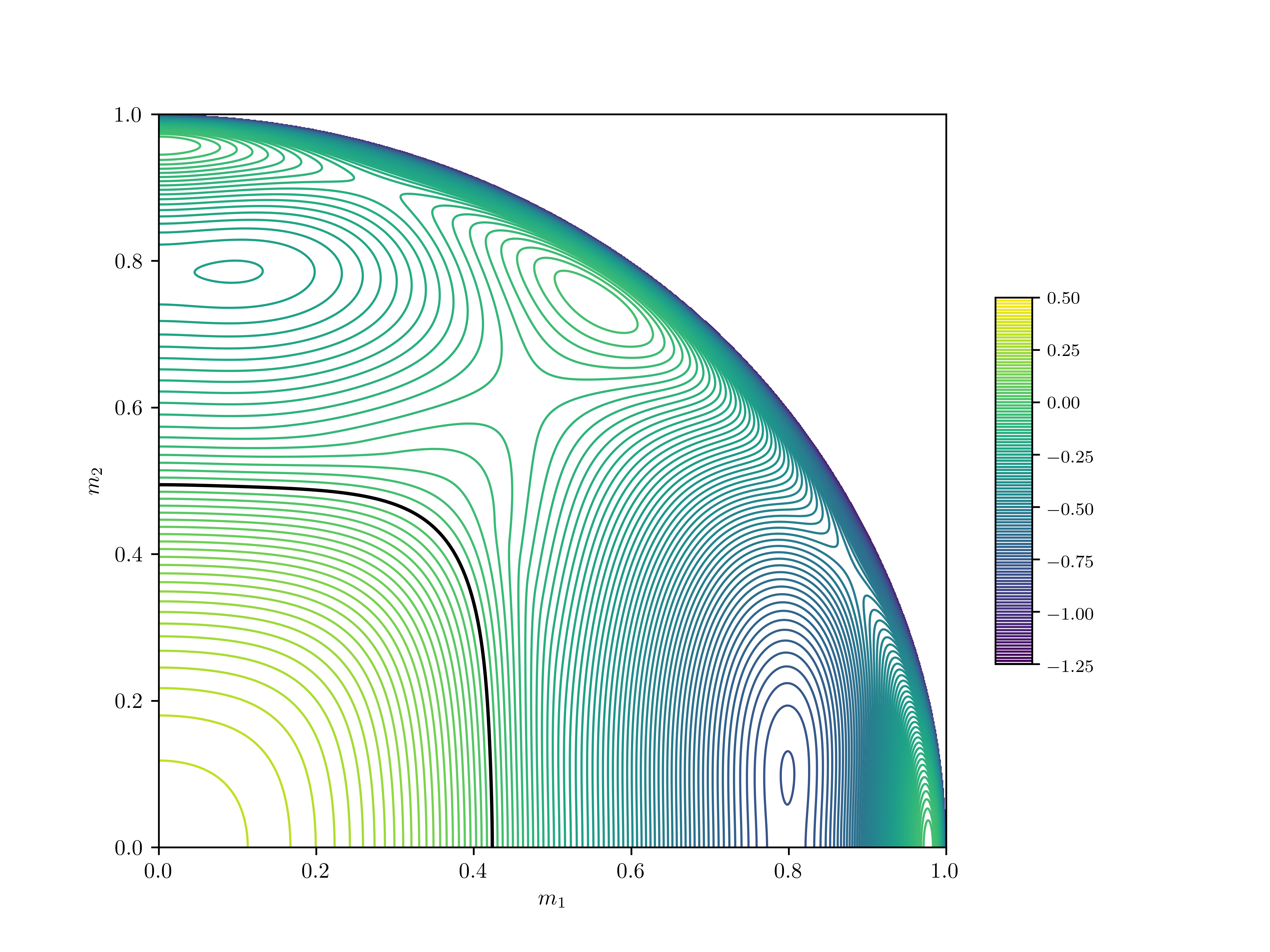}
\caption{\(\lambda_1 = 2\) and \(\lambda_2 = 1.5\).}
\label{contour4}
\end{subfigure}

\caption{Contour plots of \(\Sigma^{\textnormal{tot}}(\boldsymbol{m})\) for \(r=2\) and \(p=k=3\), with different values of \(\lambda_1\) and \(\lambda_2\). The black curve marks the boundary between regions of positive and nonpositive complexity, where \(\Sigma^{\textnormal{tot}}(\boldsymbol{m})\) vanishes. As \(\lambda_1\) and \(\lambda_2\) exceed a critical threshold, additional regions of zero complexity emerge alongside the black curve, indicating the presence of groups of critical points highly correlated with the spikes. Remarkably, in (d), a region of zero complexity emerges, where both \(m_1\) and \(m_2\) are large, signaling the presence of critical points that are simultaneously highly correlated with both spikes.}
\label{contour}
\end{figure}

Theorem~\ref{thm variational problem crt points} extends the single-spike result of~\cite[Proposition 2.4]{MR4011861} to the multi-spike setting (\(r \ge 1\)) and allows for general values of \(k_i \ne p\). In the single-spike case (\(r = 1\)), the condition \(\eta(\boldsymbol{m}) \le \eta_{\mathrm{c}}(p,k)\) of Theorem~\ref{thm variational problem crt points} reduces to the simpler requirement \(\lambda \ge \lambda_{\mathrm{c}}(p,k)\), where 
\begin{equation} \label{eq: crit threshold r=1 1}
\lambda_\text{c} (p,k) = 
\begin{cases}
\sqrt{\frac{p (k-1)^{k-1}}{2k^2 (k-2)^{k-2}}} & \text{if} \enspace k \ge p, \\[1ex]
\sqrt{\frac{p (p-1)^{k-1}}{2k^2 (p-2)^{k-2}}} & \text{if} \enspace k \le p.
\end{cases}
\end{equation}
This threshold delineates a sharp transition in the nature of the critical points:
\begin{enumerate}
\item When \(\lambda < \lambda_{\mathrm{c}}(p,k)\), most critical points are \emph{uninformative}, in the sense that their scalar product with the true signal \(\boldsymbol{u}\) is close to zero.
\item When \(\lambda \ge \lambda_{\mathrm{c}}(p,k)\), \emph{informative} (or \emph{good}) critical points appear, characterized by a large scalar product with \(\boldsymbol{u}\).
\end{enumerate}
In the multi-spike case, the landscape becomes substantially richer, and several distinct regimes arise depending on the values of \(\lambda_1, \ldots, \lambda_r\). Below, we outline the behavior in the case \(r = 2\) and \(\lambda_1 \ge \lambda_2 >0\), though the same qualitative features extend naturally to \(r \geq 3\):
\begin{enumerate}
\item When \(\lambda_\text{c} > \lambda_1 \geq \lambda_2\), there is a band containing a sub-exponential number of critical points with small scalar products \(m_1\) and \(m_2\) (see Figures~\ref{3D complexity - small} and~\ref{contour1}). Most critical points in this regime are therefore uninformative. This band corresponds to the black curve in Figures~\ref{3D complexity} and~\ref{contour}. 
\item When \(\lambda_1 \ge \lambda_\text{c} > \lambda_2\), a new region of zero complexity emerges, characterized by \(m_2 \approx 0\) and large \(m_1\) approaching \(1\) as \(\lambda_1\) increases (see Figure~\ref{contour2}). This indicates the presence of a sub-exponential number of critical points that have a large correlation with \(\boldsymbol{u}_1\). 
\item When \(\lambda_1 \geq \lambda_2 > \lambda_\text{c}\) and \(\lambda_1^{- \frac{2}{k-2}} + \lambda_2^{- \frac{2}{k-2}} > \eta_{\text{c}}\), another zero complexity region appears, in which \(m_1 \approx 0\) and \(m_2\) is large (see Figure~\ref{contour3}). In this regime, two distinct groups of informative critical points can be identified, each highly correlated with one of the spikes \(\boldsymbol{u}_1\) or \(\boldsymbol{u}_2\).
\item When \(\lambda_1 \geq \lambda_2 \ge \lambda_{\text{c}}\) and \(\lambda_1^{- \frac{2}{k-2}} + \lambda_2^{- \frac{2}{k-2}} \leq \eta_{\text{c}}\), yet another zero complexity region emerges, where critical points exhibit large correlation with both \(\boldsymbol{u}_1\) and \(\boldsymbol{u}_2\) simultaneously (see Figures~\ref{3D complexity - large} and~\ref{contour4}). 
\end{enumerate}
As \(\lambda_1\) and \(\lambda_2\) increase beyond the critical thresholds, three distinct groups of informative critical points emerge, corresponding to cases 2--4 above. The most distinctive phenomenon occurs in the last regime: a zero-complexity region appears in which the relevant critical points have nontrivial correlation with several spike directions simultaneously. 

\begin{figure}[ht]
\centering
\begin{subfigure}[b]{0.45\textwidth}
\centering
\includegraphics[scale=0.42]{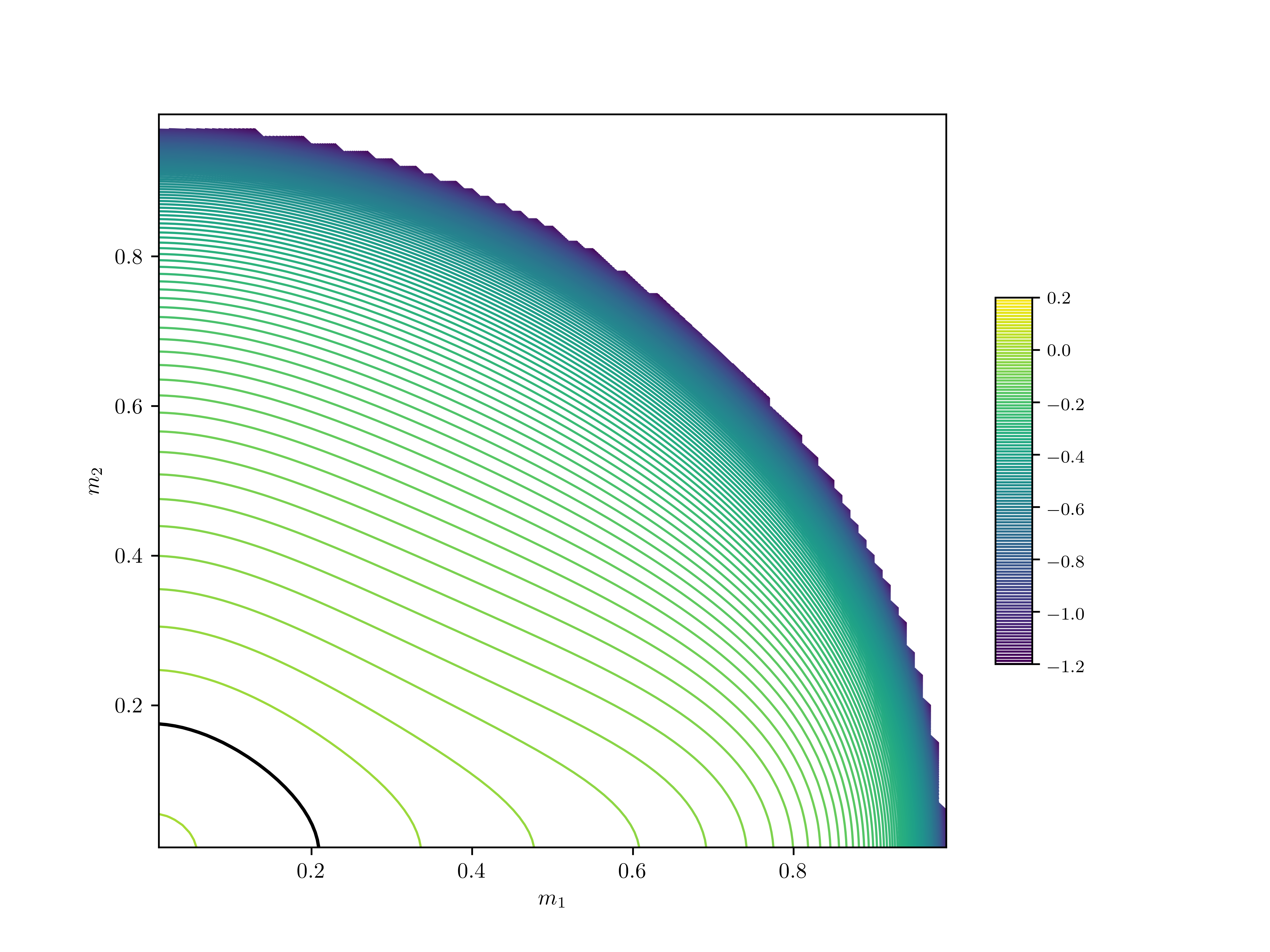}
\caption{\(\lambda_1 = 0.5\) and \(\lambda_2 = 0.2\).}
\label{max1}
\end{subfigure}
\hspace{0.06\textwidth} 
\begin{subfigure}[b]{0.45\textwidth}
\centering
\includegraphics[scale=0.42]{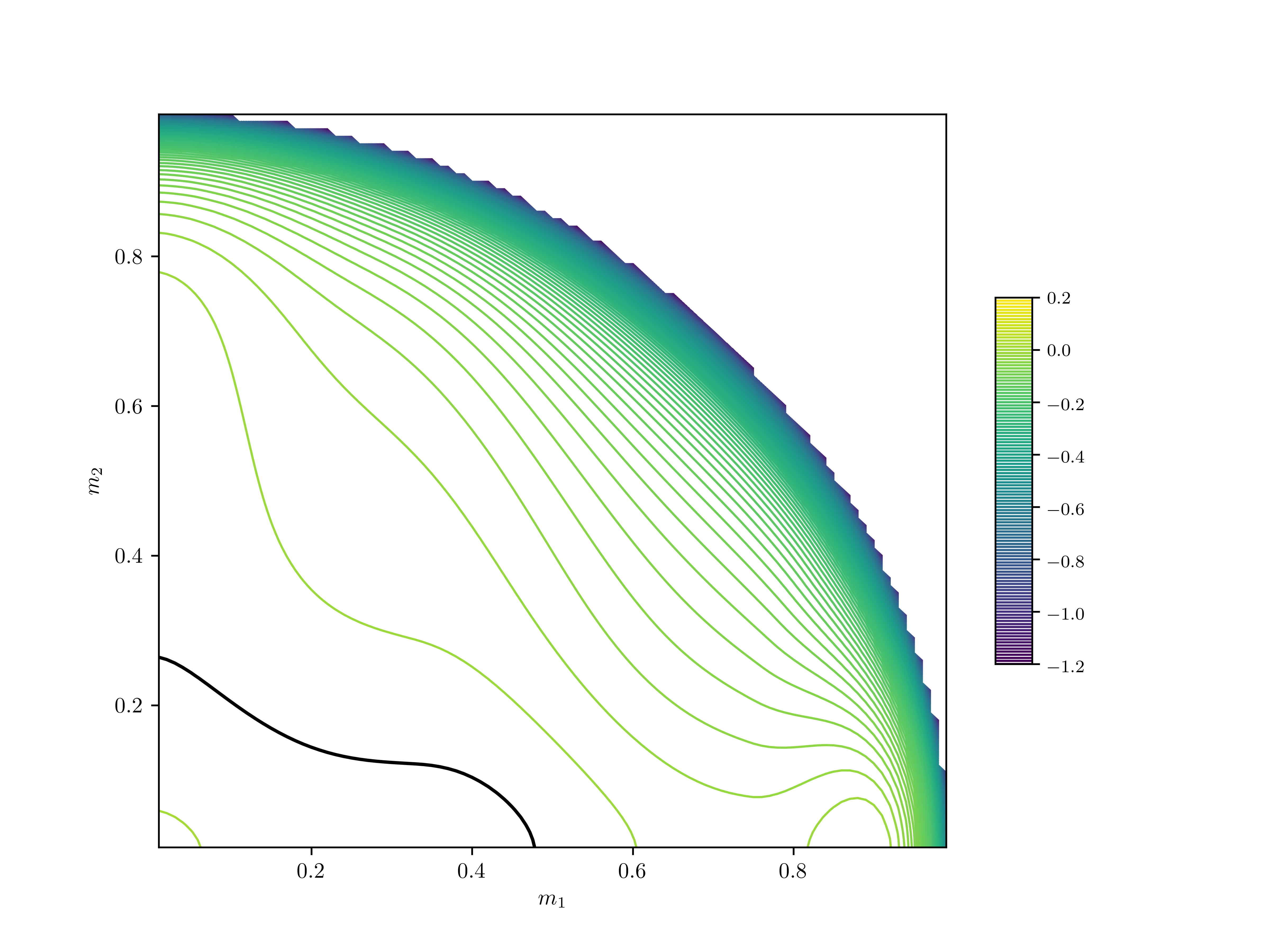}
\caption{\(\lambda_1 = 1\) and \(\lambda_2 = 0.7\).}
\label{max2}
\end{subfigure}

\vspace{0.5cm} 

\begin{subfigure}[b]{0.45\textwidth}
\centering
\includegraphics[scale=0.42]{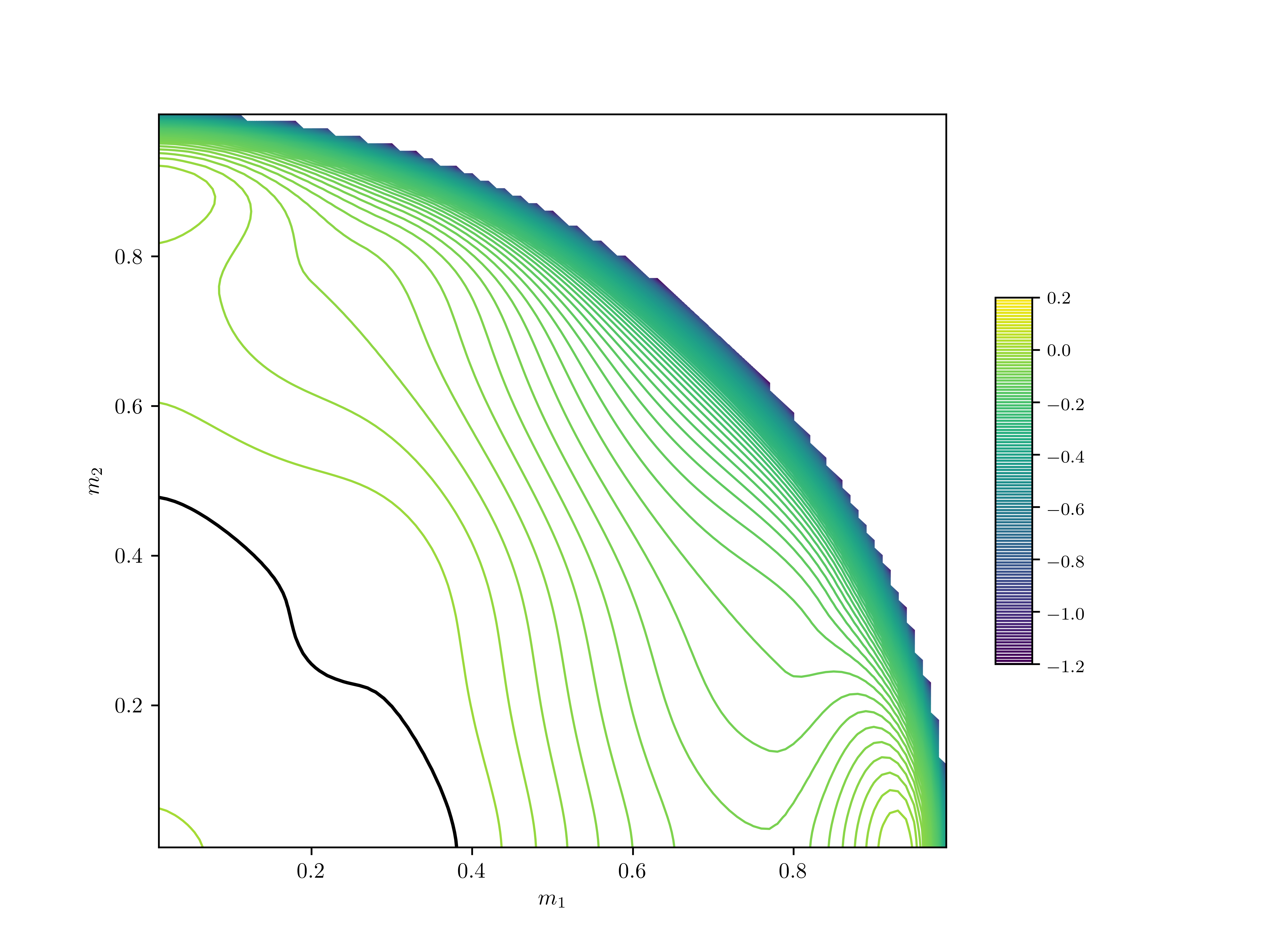}
\caption{\(\lambda_1 = 1.2\) and \(\lambda_2 = 1\).}
\label{max3}
\end{subfigure}
\hspace{0.06\textwidth} 
\begin{subfigure}[b]{0.45\textwidth}
\centering
\includegraphics[scale=0.42]{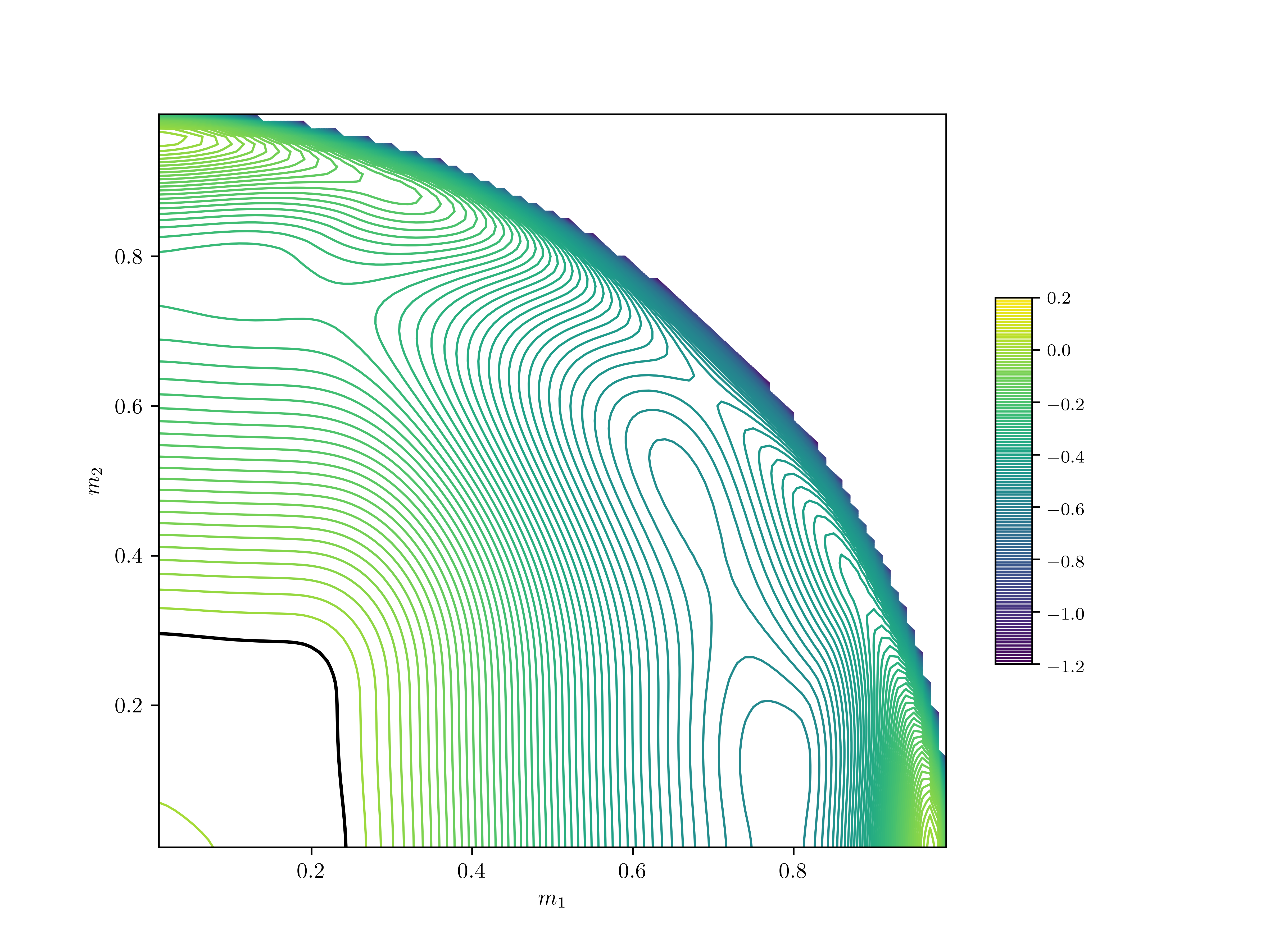}
\caption{\(\lambda_1 = 1.8\) and \(\lambda_2 = 1.5\).}
\label{max4}
\end{subfigure}
\caption{Contour plots of \(\Sigma^{\textnormal{max}}(\boldsymbol{m})\) for \(r=2\) and \(p=k=3\), with different values of \(\lambda_1\) and \(\lambda_2\). The black curve delineates the regions of positive and nonpositive complexity, where \(\Sigma^{\textnormal{max}}(\boldsymbol{m})\) vanishes. Along the curve, local maxima are uninformative, as their scalar product with either \(\boldsymbol{u}_1\) or \(\boldsymbol{u}_2\) is small. As \(\lambda_1\) and \(\lambda_2\) exceed a critical threshold, two new maxima appear where the complexity vanishes, indicating the presence of sub-exponentially many local maxima that are strongly correlated with either \(\boldsymbol{u}_1\) or \(\boldsymbol{u}_2\). This suggests a phase transition in the complexity landscape, where sufficiently strong signals lead to the emergence of groups of local maxima highly correlated with one spike.}
\label{figmax}
\end{figure}

We now turn to the variational problem of Theorem~\ref{thm local maxima} and analyze the complexity function \(\Sigma^{\textnormal{max}}(\boldsymbol{m})\) for \(\boldsymbol{m} \in \mathcal{D}_r^+\). For local maxima, we do not obtain an analogue of Theorem~\ref{thm variational problem crt points} in full generality. Instead, we analyze the variational formula numerically in the case \(r=2\) and \(\lambda_1 \ge \lambda_2 >0\). In this setting, the two ordered eigenvalues \(\gamma_1(\boldsymbol{m})\) and \(\gamma_2(\boldsymbol{m})\) of the matrix \(\boldsymbol{G}(\boldsymbol{m})^{1/2} \boldsymbol{\Theta}(\boldsymbol{m}) \boldsymbol{G}(\boldsymbol{m})^{1/2}\) introduced in Definition~\ref{def: perturbation matrix} are given by
\[
\begin{split}
\gamma_{1,2}(\boldsymbol{m}) & = \frac{1}{2} \left( \theta_1(m_1) + \theta_2(m_2) \right) \\
& \quad \pm \frac{1}{2} \sqrt{ \left( \theta_1(m_1) - \theta_2(m_2) \right)^2 + 4 \theta_1(m_1) \theta_2(m_2) G_{12}(\boldsymbol{m})^2},
\end{split}
\]
where \(\theta_i\) is defined in~\eqref{function theta} and \(G_{12}(\boldsymbol{m})\) in~\eqref{overlap v}. Since we restrict our analysis to \(m_1,m_2 \in [0,1]\), \(\theta_i(m_i) \ge0\) and so \(\boldsymbol{G}(\boldsymbol{m})^{1/2} \boldsymbol{\Theta}(\boldsymbol{m}) \boldsymbol{G}(\boldsymbol{m})^{1/2}\) is positive semi-definite. Consequently, its eigenvalues satisfy \(\gamma_1(\boldsymbol{m}) \geq \gamma_2(\boldsymbol{m}) \geq 0\). Numerical results, illustrated in Figure~\ref{figmax}, suggest the existence of a threshold \(\lambda_{\text{s}} = \lambda_{\text{s}}(p,k)\) with the following qualitative behavior:
\begin{enumerate}
\item When \(\lambda_{\text{s}} > \lambda_1 \ge \lambda_2\), a band of local maxima appears (depicted by the black curve in Figure~\ref{max1}), but these are uninformative.
\item When \(\lambda_1 \ge \lambda_{\text{s}} > \lambda_2\), a new region of zero complexity emerges, characterized by \(m_2 \approx 0\) and with \(m_1\) increasing toward \(1\) as \(\lambda_1\) grows. This indicates the presence of local maxima with large correlation with \(\boldsymbol{u}_1\) (see Figure~\ref{max2}). 
\item When \(\lambda_2 > \lambda_{\text{s}}\), another region with \(\Sigma^{\textnormal{max}} (m_1,m_2)= 0\) arises near \(m_1 \approx 0\) and large \(m_2\), indicating the presence of local maxima highly correlated with \(\boldsymbol{u}_2\) (see Figure~\ref{max3}).
\end{enumerate}
In contrast to the case of total critical points, we do not observe any region where \(\Sigma^{\textnormal{max}}(\boldsymbol{m})\) vanishes for simultaneously large values of both \(m_1\) and \(m_2\) (unlike the fourth regime described above for the total complexity). This suggests that, at least at the annealed level captured by \(\Sigma^{\mathrm{max}}\), critical points strongly correlated with multiple spike directions do not contribute to the local-maxima complexity and are more likely to be saddles.

The threshold \(\lambda_{\text{s}}(p,k)\) numerically coincides with the one identified in the single-spike setting~\cite[Subsection 2.3]{MR4011861}. Indeed, since the function \(\Sigma^{\textnormal{max}}(m_1,m_2)\) vanishes in neighborhoods where either \(m_1 \approx 0\) or \(m_2 \approx 0\), the variational problem effectively reduces to the rank-one case. Since the function \(L\) given in Definition~\ref{def: function L} and appearing in the definition of $\Sigma^{\textnormal{max}}$ is nonnegative, it follows that \(\lambda_{\mathrm{s}}(p,k) \ge \lambda_{\mathrm{c}}(p,k)\). In Section~\ref{analysis complexity function} (specifically, Lemma~\ref{lem: threshold for max}), we show that for \(r=1\) and \(k =  p\), this threshold \(\lambda_{\text{s}}(k,k)\) actually coincides with the critical threshold \(\lambda_{\text{c}}(k,k)\). Numerical simulations further suggest that the same holds for \(k \neq p\). 

However, although the value \(\lambda_{\textnormal{s}}(p,k)\) captures the phase transition for \(\Sigma^{\textnormal{max}}(\boldsymbol{m})\), this threshold does not provide any information about the location of the global maximum. This reflects a limitation of our approach: we analyze expected counts of critical points rather than their typical behavior. As discussed in Subsection~\ref{related work}, the \emph{typical} landscape structure is described by the \emph{quenched complexity}. In the single-spike setting, the authors of~\cite{ros2019complex} computed the quenched complexity and showed that, when \(\lambda \ge \lambda_{\text{c}}(p,k)\), the annealed complexity \(\Sigma^{\textnormal{max}}(m)\) vanishes for large \(m\), but the global maximum still remains close to the equator (\(m \approx 0\)). They further identified a second threshold above which the global maximum moves to regions with larger correlation with the spike. A complete understanding of the phase transition in the location of the global maximum for \(\Sigma^{\textnormal{max}}(\boldsymbol{m})\) in the multi-spike case would therefore require computing the quenched complexity in this setting. This remains an important direction for future work.

\section{Background} \label{background}

In this section, we state all the preliminary results that will be used throughout the paper. In particular, in Subsection~\ref{landscape complexity}, we characterize the complexity of the landscape via the Kac--Rice formula, while in Subsection~\ref{LDP}, we provide the large deviation principle for the largest eigenvalue of a finite-rank perturbation of a Gaussian Wigner matrix. We recall that an \(N\times N\) matrix \(\boldsymbol{W}_N\) in the Gaussian Orthogonal Ensemble (GOE), denoted by \(\boldsymbol{W}_N \sim \text{GOE}(N)\), is a real symmetric matrix whose upper-triangular entries are independent Gaussian random variables with mean zero and variance \(\E \left [ (\boldsymbol{W}_N)_{ij}^2 \right ] = \frac{1+\delta_{ij}}{N}\) for \(1\le i\le j\le N\).

\subsection{Kac--Rice formula and landscape complexity} \label{landscape complexity}

The Kac--Rice formula is the classical tool for computing the expected number of critical points of a smooth Gaussian random function. We apply Theorem 12.1.1 and Corollary 12.1.2 of~\cite{MR2319516} to the random function \(f_N\) defined in~\eqref{eq: function f}. We refer the reader to~\cite{MR2319516,MR2478201} for a broader introduction to the Kac--Rice method. 

\begin{lem}[Kac--Rice formula~\cite{MR2319516}] \label{thm: kac rice formula}
For all Borel sets \(M_1,\ldots, M_r \subset [-1,1]\) and \(B \subset \R\), and for all \(N\), it holds that
\begin{equation} \label{Kac--Rice tot}
\begin{split}
\E \left [ \textnormal{Crt}_N^{\textnormal{tot}} (M,B) \right ] & = \int_{\{\boldsymbol{\sigma} \in \mathbb{S}^{N-1} \colon \left ( \langle \boldsymbol{\sigma}, \boldsymbol{u}_i \rangle \right)_{i=1}^r \in M\}} \textnormal{d}\boldsymbol{\sigma} \, \varphi_{\nabla f_N(\boldsymbol{\sigma})}(\boldsymbol{0})  \\
& \quad \times \E \left [ \left |\det  \nabla^2 f_N(\boldsymbol{\sigma})  \right | \mathbf{1}_{\{f_N(\boldsymbol{\sigma})\in B\}} \Big | \nabla f_N(\boldsymbol{\sigma})=\boldsymbol{0} \right ] ,
\end{split}
\end{equation}
where $M=M_1\times\cdots\times M_r$, \(\textnormal{d}\boldsymbol{\sigma}\) denotes the surface measure on \(\mathbb{S}^{N-1}\), and \(\varphi_{\nabla f_N(\boldsymbol{\sigma})}\) denotes the density of the Riemannian gradient \(\nabla f_N(\boldsymbol{\sigma})\). More generally, let \( \textnormal{Crt}_N^{\ell} (M, B)\) denote the number of critical points of index \(\ell \in \{0,\ldots, N-1\}\), where the index is the number of negative eigenvalues of the Riemannian Hessian \(\nabla^2 f_N\). Then 
\begin{equation} \label{Kac--Rice index}
\begin{split}
\E \left [ \textnormal{Crt}_N^{\ell} (M, B) \right ] & = \int_{\{\boldsymbol{\sigma} \in \mathbb{S}^{N-1} \colon \left ( \langle \boldsymbol{\sigma}, \boldsymbol{u}_i \rangle \right)_{i=1}^r \in M\}}  \textnormal{d}\boldsymbol{\sigma} \, \varphi_{\nabla f_N(\boldsymbol{\sigma})}(\boldsymbol{0})\\
& \quad \times \E\left [ \left |\det \nabla^2  f_N(\boldsymbol{\sigma})  \right |  \mathbf{1}_{ \left \{i(\nabla^2 f_N(\boldsymbol{\sigma}))=\ell, \, f_N(\boldsymbol{\sigma})\in B \right \}} \Big |\nabla f_N(\boldsymbol{\sigma})=\boldsymbol{0} \right ] .
\end{split}
\end{equation}
\end{lem}

Recall the admissible domain \(\mathcal{D}_r = \left \{ \boldsymbol{m} \in [-1,1]^r \colon \|\boldsymbol{m}\|_2 < 1 \right\}\) from Definition~\ref{def: function Sigma tot}. Define \(S\colon [-1,1]^r \times\R\to\R\cup\{-\infty\}\) as follows. For \((\boldsymbol{m},x)  \in \mathcal{D}_r^{\mathrm{c}} \times \R\), we set \(S(\boldsymbol{m},x) = - \infty\). For \((\boldsymbol{m},x) \in \mathcal{D}_r \times \R\), we set
\begin{equation} \label{function S}
\begin{split}
S(\boldsymbol{m},x) & = \frac{1}{2} (\log(p-1) +1) + \frac{1}{2} \log \left (1 - \|\boldsymbol{m}\|_2^2 \right ) - \frac{1}{p} \sum_{i=1}^r \lambda_i^2 k_i^2 m_i^{2k_i-2}  \\
& \quad + \frac{1}{p} \left( \sum_{i=1}^r \lambda_i  k_i m_i^{k_i} \right)^2 - \left ( x - \sum_{i=1}^r \lambda_i m_i^{k_i}\right )^2.
\end{split}
\end{equation}

We next introduce the finite-rank matrix that appears in the conditional law of the Hessian. For \(\boldsymbol{m}\in\mathcal{D}_r\) and \(N>r\), let \(\boldsymbol{v}_1(\boldsymbol{m}),\ldots, \boldsymbol{v}_r(\boldsymbol{m}) \in\mathbb{S}^{N-2}\) be vectors whose Gram matrix is \(\boldsymbol{G}(\boldsymbol{m})\) (see~\eqref{overlap v}), that is, 
\[ 
\left\langle \boldsymbol{v}_i(\boldsymbol{m}), \boldsymbol{v}_j(\boldsymbol{m}) \right\rangle = G_{ij}(\boldsymbol{m}) = \frac{\delta_{ij}-m_im_j} {\sqrt{1-m_i^2}\sqrt{1-m_j^2}},
\] 
for all \(1\leq i,j\leq r\). Such a family exists because \(\boldsymbol{G}(\boldsymbol{m})\) is positive definite on \(\mathcal{D}_r\) of rank \(r\) and \(r \le N-1\). Define 
\begin{equation}\label{def: P matrix} 
\boldsymbol{P}_{N-1}(\boldsymbol{m}) \coloneqq \sum_{i=1}^r \theta_i(m_i) \boldsymbol{v}_i(\boldsymbol{m}) \boldsymbol{v}_i(\boldsymbol{m})^\top, 
\end{equation} 
where \(\theta_i\) is defined in~\eqref{function theta}. 

The following proposition gives an exact finite-$N$ representation of the expected critical-point counts.

\begin{prop} \label{formula complexity}
For all Borel sets \(M_1,\ldots, M_r \subset [-1,1]\) and \(B \subset \R\) and for all \(N >r\), it holds that
\begin{equation}\label{Kac--Rice tot 2}
\begin{split}
\E \left [ \textnormal{Crt}_N^{\textnormal{tot}} (M, B) \right] & = C_N \int_B \textnormal{d}x \int_{M \cap \mathcal{D}_r} \textnormal{d} \boldsymbol{m}  \,  \left (1- \| \boldsymbol{m}\|_2^2\right )^{-\frac{r+2}{2}}  \\
& \quad \times  \exp \{N S(\boldsymbol{m},x)\} \E \left [ |\det \boldsymbol{H}_{N-1}(\boldsymbol{m},x) |\right ] ,
\end{split}
\end{equation}
and 
\begin{equation} \label{Kac--Rice index 2}
\begin{split}
\E \left [ \textnormal{Crt}_N^{\textnormal{max}} (M, B) \right ] & = C_N \int_B \textnormal{d}x \int_{M \cap \mathcal{D}_r} \textnormal{d} \boldsymbol{m}  \,  \left (1- \| \boldsymbol{m}\|_2^2\right )^{-\frac{r+2}{2}} \\
&\quad  \times \exp \{N S(\boldsymbol{m},x)\}  \E \left [ |\det \boldsymbol{H}_{N-1}(\boldsymbol{m},x) | \mathbf{1}_{\{\boldsymbol{H}_{N-1}(\boldsymbol{m},x) \prec 0\}}\right ] .
\end{split}
\end{equation}
Here, \(\boldsymbol{H}_{N-1}(\boldsymbol{m},x)\) is an \((N-1)\)-dimensional square random matrix distributed as
\begin{equation} \label{matrix H_N}
\boldsymbol{H}_{N-1}(\boldsymbol{m},x) \stackrel{d}{=} \boldsymbol{W}_{N-1} + \sqrt{\frac{N}{N-1}} \boldsymbol{P}_{N-1}(\boldsymbol{m}) -\sqrt{\frac{N}{N-1}} t(\boldsymbol{m},x) \boldsymbol{I}_{N-1}, 
\end{equation}
where \(\boldsymbol{W}_{N-1} \sim \textnormal{GOE}(N-1)\), \(\boldsymbol{P}_{N-1} (\boldsymbol{m})\) is given by~\eqref{def: P matrix}, and \(t(\boldsymbol{m},x)\) by~\eqref{function t}. Furthermore, the constant \(C_N = C (N,r,p)\) is exponentially trivial, i.e., \(\lim_{N \to \infty} \frac{1}{N} \log C_N=0\). 
\end{prop}

\begin{rmk} \label{rmk: KR S0}
For later use, we rewrite the Kac--Rice integrand of Proposition~\ref{formula complexity} in a slightly different form. Set $q(\boldsymbol{m}) \coloneqq 1- \| \boldsymbol{m} \|_2^2$.
On \(\mathcal{D}_r\times\R\), define
\[
S_0 (\boldsymbol{m},x) \coloneqq S(\boldsymbol{m},x)-\frac{1}{2} \log q(\boldsymbol{m}).
\]
Then, for \((\boldsymbol{m},x)\in \mathcal{D}_r\times\R\),
\[
q(\boldsymbol{m})^{-\frac{r+2}{2}} \exp\{N S(\boldsymbol m,x)\} = \exp\left\{ N S_0(\boldsymbol m,x) + \frac{N-r-2}{2}\log q(\boldsymbol m) \right\}.
\]
\end{rmk}

We prove Proposition~\ref{formula complexity} using Lemma~\ref{thm: kac rice formula} and the following result which provides the joint law of the Gaussian random variables \((f_N(\boldsymbol{\sigma}), \nabla f_N(\boldsymbol{\sigma}), \nabla^2 f_N(\boldsymbol{\sigma}))\).

\begin{lem} \label{prop}
For all Borel sets \(M_1,\ldots, M_r \subset [-1,1]\) and \(B \subset \R\) and for all \(N > r\), it holds that
\begin{equation} \label{Kac--Rice tot 3}
\begin{split}
\E \left [ \textnormal{Crt}_N^{\textnormal{tot}} (M, B) \right ]& =  \omega_{N-r} \int_{M \cap \mathcal{D}_r} \textnormal{d} \boldsymbol{m}  \,  \varphi_{\boldsymbol{g}_{\boldsymbol{m}}}(\boldsymbol{0}) \left (1- \| \boldsymbol{m}\|_2^2\right )^{\frac{N-r-2}{2}}  \E\left [ |\det \boldsymbol{H}_{\boldsymbol{m}} | \mathbf{1}_{\{f_{\boldsymbol{m}} \in B\}}\right ], 
\end{split}
\end{equation}
and
\begin{equation}  \label{Kac--Rice index 3}
\begin{split} 
\E \left [ \textnormal{Crt}_N^{\textnormal{max}} (M, B) \right ]&= \omega_{N-r} \int_{M \cap \mathcal{D}_r} \textnormal{d} \boldsymbol{m}  \, \varphi_{\boldsymbol{g}_{\boldsymbol{m}}}(\boldsymbol{0})  \left (1- \| \boldsymbol{m}\|_2^2\right )^{\frac{N-r-2}{2}} \E\left [  |\det \boldsymbol{H}_{\boldsymbol{m}} | \mathbf{1}_{\left \{ \boldsymbol{H}_{\boldsymbol{m}} \prec 0 \right \}} \mathbf{1}_{\{f_{\boldsymbol{m}} \in B\}} \right ], 
\end{split}
\end{equation}
where \(\omega_N = \frac{2\pi^{N/2}}{\Gamma(N/2)}\) is the volume of the \((N-1)\)-dimensional unit sphere. Furthermore, the joint distribution of \(f_{\boldsymbol{m}} \in \R\), \(\boldsymbol{g}_{\boldsymbol{m}} \in \R^{N-1}\), and \(\boldsymbol{H}_{\boldsymbol{m}} \in \R^{(N-1) \times (N-1)}\) is given by
\begin{equation} \label{joint distr}
\begin{split}
f_{\boldsymbol{m}} &\stackrel{d}{=} \sum_{i=1}^r\lambda_i m_i^{k_i} + \frac{1}{\sqrt{2N}}Z, \\
\boldsymbol{g}_{\boldsymbol{m}} &\stackrel{d}{=}  \sum_{i=1}^r \lambda_i  k_i m_i^{k_i-1} \sqrt{1-m_i^2} \boldsymbol{v}_i(\boldsymbol{m}) + \sqrt{\frac{p}{2N}} \boldsymbol{g}_{N-1},\\
\boldsymbol{H}_{\boldsymbol{m}} &\stackrel{d}{=} \sum_{i=1}^r \lambda_i k_i (k_i-1) m_i^{k_i-2} (1-m_i^2) \boldsymbol{v}_i(\boldsymbol{m}) \boldsymbol{v}_i(\boldsymbol{m})^\top \\
& \quad  + \sqrt{\frac{p(p-1)(N-1)}{2N}} \boldsymbol{W}_{N-1} - \left ( \sum_{i=1}^r \lambda_i k_i m_i^{k_i} + \frac{p}{\sqrt{2N}} Z\right ) \boldsymbol{I}_{N-1},
\end{split}
\end{equation}
where \(Z \sim \mathcal{N}(0,1)\), \(\boldsymbol{g}_{N-1} \sim \mathcal{N}(0, \boldsymbol{I}_{N-1})\), and \(\boldsymbol{W}_{N-1} \sim \textnormal{GOE}(N-1)\) are independent. Moreover, the vectors \(\boldsymbol{v}_1 (\boldsymbol{m}), \ldots, \boldsymbol{v}_r (\boldsymbol{m}) \in \mathbb{S}^{N-2}\) have Gram matrix \(\boldsymbol{G}(\boldsymbol{m})\) defined in~\eqref{overlap v}. 
\end{lem}

Having Lemma~\ref{prop} at hand, we turn to the proof of Proposition~\ref{formula complexity}.

\begin{proof}[Proof of Proposition~\ref{formula complexity}]
According to~\eqref{joint distr}, we have that \(f_{\boldsymbol{m}}\) is distributed as a Gaussian random variable with mean \(\sum_{i=1}^r \lambda_i m_i^{k_i}\) and variance \(\frac{1}{2N}\), hence the density function \(p_{f_{\boldsymbol{m}}}\) is given by
\begin{equation} \label{density f_p}
p_{f_{\boldsymbol{m}}}(x) = \sqrt{\frac{N}{\pi}} \exp \left \{ - N \left (x - \sum_{i=1}^r \lambda_i m_i^{k_i} \right )^2 \right \}.
\end{equation}
Therefore, the inner expectation in~\eqref{Kac--Rice index 3} can be expanded as
\begin{equation} \label{inner exp}
\begin{split}
&\E\left [ |\det \boldsymbol{H}_{\boldsymbol{m}}| \mathbf{1}_{\left \{\boldsymbol{H}_{\boldsymbol{m}} \prec 0 \right\}} \mathbf{1}_{\{ f_{\boldsymbol{m}} \in B\}}  \right ]\\
& = \left ( \frac{p(p-1)(N-1)}{2N} \right )^{\frac{N-1}{2}} \int_{B} \E\left [ |\det \boldsymbol{H}_{N-1}(\boldsymbol{m},x)| \mathbf{1}_{\{ \boldsymbol{H}_{N-1}(\boldsymbol{m},x) \prec 0\}}  \right ] p_{f_{\boldsymbol{m}}}(x) \text{d}x,
\end{split}
\end{equation}
where the matrix \(\boldsymbol{H}_{N-1}(\boldsymbol{m},x)\) is distributed as in~\eqref{matrix H_N}. Further, the random vector \(\boldsymbol{g}_{\boldsymbol{m}}\) is Gaussian with mean \(\sum_{i=1}^r \lambda_i  k_i m_i^{k_i-1} \sqrt{1-m_i^2} \boldsymbol{v}_i(\boldsymbol{m})\) and covariance matrix \(\frac{p}{2N} \boldsymbol{I}_{N-1}\). Thus its density function at \(\boldsymbol{0}\) is given by 
\begin{equation} \label{density g}
\begin{split}
&\varphi_{\boldsymbol{g}_{\boldsymbol{m}} } (\boldsymbol{0}) = \left (\frac{N}{\pi p} \right )^{\frac{N-1}{2}} \exp \left \{- \frac{N}{p} \sum_{i=1}^r \lambda_i^2 k_i^2 m_i^{2k_i-2} + \frac{N}{p} \left (\sum_{i=1}^r \lambda_i  k_i m_i^{k_i} \right)^2\right \},
\end{split}
\end{equation}
where we used~\eqref{overlap v}. Plugging~\eqref{density f_p}-\eqref{density g} into~\eqref{Kac--Rice index 3}, we obtain
\[
\begin{split}
& \E \left [\textnormal{Crt}_N^{\textnormal{max}}(M,B)\right ] \\
& =  \left (\frac{p(p-1)(N-1)}{2N} \right )^{\frac{N-1}{2}}  \omega_{N-r}  \int_{M \cap \mathcal{D}_r} \text{d} \boldsymbol{m} \,  \varphi_{\boldsymbol{g}_{\boldsymbol{m}}}(\boldsymbol{0})  \left (1- \| \boldsymbol{m} \|_2^2  \right )^{\frac{N-r-2}{2}}\\
& \quad \times \int_B \text{d}x  \, \E \left [ |\det \boldsymbol{H}_{N-1}(\boldsymbol{m},x) |  \mathbf{1}_ {\{ \boldsymbol{H}_{N-1}(\boldsymbol{m},x) \prec 0\}} \right ] p_{f_{\boldsymbol{m}}}(x) \\
& = C_N \int_B \text{d}x \int_{M \cap \mathcal{D}_r} \text{d} \boldsymbol{m} \, \left (1-\| \boldsymbol{m} \|_2^2 \right )^{-\frac{r+2}{2}}  \\
& \quad \times \E \left [ |\det \boldsymbol{H}_{N-1}(\boldsymbol{m},x) | \mathbf{1}_{\{\boldsymbol{H}_{N-1}(\boldsymbol{m},x) \prec 0\}} \right ] \exp \{NS(\boldsymbol{m},x)\},
\end{split}
\]
where \(C_N = C(N,r,p)\) is given by
\[
C(N,r,p) = 2 \left (\frac{N-1}{2e}\right)^{\frac{N-1}{2}} \Gamma\left(\frac{N-r}{2}\right)^{-1} \pi^{- \frac{r-1}{2}} \left ( \frac{N}{(p-1) e \pi} \right )^{1/2}.
\]
This proves~\eqref{Kac--Rice index 2}. Identity~\eqref{Kac--Rice tot 2} follows in the same way. It is straightforward to show that \(C_N\) is exponentially trivial by expanding the function \(\Gamma\) using the Stirling's formula.  
\end{proof}

It remains to derive the random matrices appearing in the Kac--Rice formula by proving Lemma~\ref{prop}.

\begin{proof}[Proof of Lemma~\ref{prop}]
Recall that 
\[
f_N(\boldsymbol{\sigma}) = \sum_{i=1}^r \lambda_i \langle \boldsymbol{u}_i , \boldsymbol{\sigma} \rangle^{k_i} + H_{N,p}(\boldsymbol{\sigma}), 
\]
where 
\[
H_{N,p}(\boldsymbol{\sigma}) = \frac{1}{\sqrt{2N}} \sum_{1 \leq i_1, \ldots, i_p \leq N} W_{i_1,\ldots,i_p}\sigma_{i_1} \cdots \sigma_{i_p}.
\]
Recall also that \( W_{i_1,\ldots,i_p} = \frac{1}{p!} \sum_{\pi\in S_p} G_{i_{\pi(1)},\ldots,i_{\pi(p)}}\), so that $W$ is symmetric, namely \(W_{i_1,\ldots,i_p} =W_{i_{\pi(1)},\ldots,i_{\pi(p)}}\) for every $\pi\in S_p$. We note that \(H_{N,p}\) is a centered Gaussian function whose law is invariant under orthogonal transformations, i.e., for every orthogonal matrix \(\boldsymbol{R} \in \mathcal{O}(N)\), 
\[
\bigl(H_{N,p}(\boldsymbol{R}\boldsymbol{\sigma})\bigr)_{\boldsymbol{\sigma}\in\mathbb{S}^{N-1}} \stackrel{d}{=}\bigl(H_{N,p}(\boldsymbol{\sigma})\bigr)_{\boldsymbol{\sigma}\in\mathbb{S}^{N-1}}.
\]
Let \(\hat{\nabla}\) and \(\hat{\nabla}^2\) denote the Euclidean gradient and Hessian, respectively. A direct computation yields
\[
\begin{split}
\hat{\nabla} f_N(\boldsymbol{\sigma}) & = \sum_{i=1}^r \lambda_i k_i \langle \boldsymbol{u}_i, \boldsymbol{\sigma} \rangle^{k_i-1} \boldsymbol{u}_i + \hat{\nabla} H_{N,p}(\boldsymbol{\sigma}),\\
\hat{\nabla}^2 f_N(\boldsymbol{\sigma}) & = \sum_{i=1}^r \lambda_i k_i(k_i-1) \langle \boldsymbol{u}_i, \boldsymbol{\sigma} \rangle^{k_i-2} \boldsymbol{u}_i \boldsymbol{u}_i^\top + \hat{\nabla}^2 H_{N,p}(\boldsymbol{\sigma}),
\end{split}
\]
where for any \(k, \ell \in \{1, \ldots, N\}\),
\[
\begin{split}
\left ( \hat{\nabla} H_{N,p}(\boldsymbol{\sigma}) \right ) _k & = \frac{p}{\sqrt{2N}} \sum_{1 \le i_1, \ldots, i_{p-1} \le N} W_{k, i_1, \ldots, i_{p-1}} \sigma_{i_1} \cdots \sigma_{i_{p-1}},\\
( \hat{\nabla}^2 H_{N,p}(\boldsymbol{\sigma}))_{k\ell} & = \frac{p(p-1)}{\sqrt{2N}}  \sum_{1 \le i_1, \ldots, i_{p-2} \le N} W_{k, \ell, i_1, \ldots, i_{p-2}} \sigma_{i_1} \cdots \sigma_{i_{p-2}}.
\end{split}
\]
The Riemannian gradient and Hessian on \(\mathbb{S}^{N-1}\) are given by
\[
\nabla f_N(\boldsymbol{\sigma}) = \text{P}^\bot_{\boldsymbol{\sigma}} \hat{\nabla} f_N(\boldsymbol{\sigma}),
\]
and
\[
\nabla^2 f_N(\boldsymbol{\sigma}) = \text{P}^\bot_{\boldsymbol{\sigma}} \hat{\nabla}^2 f_N(\boldsymbol{\sigma}) \text{P}^\bot_{\boldsymbol{\sigma}} - \langle \boldsymbol{\sigma},\hat{\nabla}f_N(\boldsymbol{\sigma}) \rangle \text{P}^\bot_{\boldsymbol{\sigma}},
\]
respectively, where \(\text{P}_{\boldsymbol{\sigma}}^\bot = I_N - \boldsymbol{\sigma} \boldsymbol{\sigma}^\top\) denotes the orthogonal projection from \(\R^N\) onto the tangent space \(\text{T}_{\boldsymbol{\sigma}} \mathbb{S}^{N-1}\). Substituting the expressions above, we obtain 
\[
\begin{split}
\nabla f_N(\boldsymbol{\sigma}) & =  \sum_{i=1}^r \lambda_i k_i \langle \boldsymbol{u}_i,\boldsymbol{\sigma} \rangle^{k_i-1} \left ( \text{P}^\bot_{\boldsymbol{\sigma}} \boldsymbol{u}_i \right ) +  \text{P}^\bot_{\boldsymbol{\sigma}} \hat{\nabla} H_{N,p}(\boldsymbol{\sigma}),\\
\nabla^2 f_N(\boldsymbol{\sigma}) & = \sum_{i=1}^r \lambda_i k_i (k_i-1) \langle \boldsymbol{u}_i, \boldsymbol{\sigma} \rangle^{k_i-2}  \left (\text{P}^\bot_{\boldsymbol{\sigma}} \boldsymbol{u}_i \right )  \left ( \text{P}^\bot_{\boldsymbol{\sigma}} \boldsymbol{u}_i \right )^\top + \text{P}^\bot_{\boldsymbol{\sigma}} \hat{\nabla}^2 H_{N,p}(\boldsymbol{\sigma}) \text{P}^\bot_{\boldsymbol{\sigma}}\\
& \quad - \left ( \sum_{i=1}^r \lambda_i k_i \langle \boldsymbol{u}_i ,\boldsymbol{\sigma} \rangle ^{k_i}+ \langle \boldsymbol{\sigma}, \hat{\nabla}H_{N,p} (\boldsymbol{\sigma}) \rangle \right ) \text{P}^\bot_{\boldsymbol{\sigma}}.
\end{split}
\]
By rotational invariance of \(H_{N,p}\), the joint law of \((f_N(\boldsymbol{\sigma}), \nabla f_N(\boldsymbol{\sigma}) , \nabla^2 f_N(\boldsymbol{\sigma}) )\) depends on \(\boldsymbol{\sigma}\) only through the correlation vector 
\[
\boldsymbol{m} = (m_1, \ldots, m_r), \qquad m_i = \langle \boldsymbol{\sigma}, \boldsymbol{u}_i \rangle. 
\]
In particular, the corresponding Kac--Rice integrand depends on $\boldsymbol{\sigma}$ only through $\boldsymbol{m}$. Choose $\boldsymbol{R}\in\mathcal{O}(N)$ such that $\boldsymbol{R}\boldsymbol{\sigma}=\boldsymbol{e}_N$, and rotate the spike directions accordingly. By isotropy, we may therefore work at $\boldsymbol{\sigma}=\boldsymbol{e}_N$. To simplify notation, we continue to denote the rotated spike directions $\boldsymbol{R}\boldsymbol{u}_1,\ldots,\boldsymbol{R}\boldsymbol{u}_r$ by $\boldsymbol{u}_1,\ldots,\boldsymbol{u}_r$. Specifically, we parametrize each \(\boldsymbol{u}_i\) by
\[ 
\boldsymbol{u}_i = m_i \boldsymbol{e}_N + \sqrt{1-m_i^2} \boldsymbol{\tilde{v}}_i(\boldsymbol{m}),
\]
where \(\boldsymbol{\tilde{v}}_i (\boldsymbol{m}) \in \mathbb{S}^{N-1} \cap \boldsymbol{e}_N^\perp\), i.e., it lies in the equatorial subspace orthogonal to \(\boldsymbol{e}_N\), and is chosen so that \(\langle  \boldsymbol{u}_i , \boldsymbol{e}_N\rangle = m_i \) and \(\langle \boldsymbol{u}_i ,\boldsymbol{u}_j \rangle = \delta_{ij}\). This implies that \(\boldsymbol{\tilde{v}}_i (\boldsymbol{m})\) has zero \(N\)th coordinate, and satisfies
\[ 
\langle \boldsymbol{\tilde{v}}_i(\boldsymbol{m}), \boldsymbol{\tilde{v}}_j(\boldsymbol{m}) \rangle  =  \frac{\delta_{ij} - m_i m_j}{\sqrt{1-m_i^2} \sqrt{1-m_j^2}}.
\] 
Identifying \(\text{T}_{\boldsymbol{e}_N} \mathbb{S}^{N-1}\) with \(\R^{N-1}\), and denoting by $\boldsymbol{v}_i(\boldsymbol{m})\in\mathbb{S}^{N-2}$ the vector obtained from $\widetilde{\boldsymbol{v}}_i(\boldsymbol{m})$ by deleting its last coordinate, the following identities hold jointly:
\[ 
\begin{split}
f_N(\boldsymbol{e}_N) &\stackrel{d}{=} \sum_{i=1}^r \lambda_i m_i^{k_i} + \frac{1}{\sqrt{2N}}Z, \\
\nabla f_N(\boldsymbol{e}_N)&\stackrel{d}{=} \sum_{i=1}^r \lambda_i k_im_i^{k_i-1} \sqrt{1-m_i^2} \boldsymbol{v}_i(\boldsymbol{m}) + \sqrt{\frac{p}{2N}} \boldsymbol{g}_{N-1},
\end{split}
\]
and
\[
\begin{split}
\nabla^2 f_N(\boldsymbol{e}_N) & \stackrel{d}{=} \sum_{i=1}^r \lambda_i k_i (k_i-1)m_i^{k_i-2} (1-m_i^2) \boldsymbol{v}_i(\boldsymbol{m}) \boldsymbol{v}_i(\boldsymbol{m})^\top\\
& \quad + \sqrt{\frac{p(p-1)(N-1)}{2N}} \boldsymbol{W}_{N-1} - \left ( \sum_{i=1}^r \lambda_ik_i m_i^{k_i} + \frac{p}{\sqrt{2N}} Z\right )\boldsymbol{I}_{N-1},
\end{split}
\]
where \(Z \sim \mathcal{N}(0,1)\), \(\boldsymbol{g}_{N-1} \sim \mathcal{N}(0, \boldsymbol{I}_{N-1})\) and \(\boldsymbol{W}_{N-1} \sim \text{GOE}(N-1)\) are independent. In particular, we note that $\nabla f_N(\boldsymbol{e}_N)$ is independent of \(\bigl(f_N(\boldsymbol{e}_N), \nabla^2f_N(\boldsymbol{e}_N) \bigr)\). Thus, the conditional expectation appearing in the Kac--Rice formula, conditioned on $\nabla f_N(\boldsymbol{e}_N)=\boldsymbol{0}$, is equal to the corresponding unconditional expectation.

We now apply the Kac--Rice formula from Lemma~\ref{thm: kac rice formula}. Since the integrands in~\eqref{Kac--Rice tot} and~\eqref{Kac--Rice index} depend on $\boldsymbol{\sigma}$ only through $\boldsymbol{m}$, we apply the co-area formula to the map \(\rho\colon\mathbb{S}^{N-1}\to\mathbb{R}^r\) given by
\[ 
\rho(\boldsymbol{\sigma}) = \bigl( \langle\boldsymbol{\sigma},\boldsymbol{u}_1\rangle, \ldots, \langle\boldsymbol{\sigma},\boldsymbol{u}_r\rangle \bigr). 
\]
For \(\boldsymbol{m}\in M\cap\mathcal{D}_r\), the level set \(\rho^{-1}(\boldsymbol{m})\) is a sphere of dimension \(N-r-1\) and radius \(\sqrt{1-\|\boldsymbol{m}\|_2^2}\). Its volume is therefore 
\[
\omega_{N-r} \left (1- \| \boldsymbol{m} \|_2^2  \right )^{\frac{N-r-1}{2}}.
\]
Moreover, since \(\left( \left\langle \nabla \rho_i(\boldsymbol{\sigma}), \nabla \rho_j(\boldsymbol{\sigma}) \right\rangle \right)_{i,j=1}^r =\boldsymbol{I}_r-\boldsymbol{m}\boldsymbol{m}^{\top}\), the Jacobian of \(\rho\) is
\[
J_{\rho}(\boldsymbol{\sigma}) = \sqrt{ \det\left( \boldsymbol{I}_r-\boldsymbol{m}\boldsymbol{m}^{\top} \right) } = \sqrt{1-\|\boldsymbol{m}\|_2^2}.
\]
Thus, for every \(0\leq \ell\leq N-1\), the coarea formula and \eqref{Kac--Rice index} yield
\begin{equation}\label{kac rice index 2}
\begin{split}
\E \left [ \textnormal{Crt}_N^{\ell} (M, B) \right ]& = \omega_{N-r}\int_{M \cap \mathcal{D}_r} \text{d} \boldsymbol{m} \, \varphi_{ \boldsymbol{g}_{\boldsymbol{m}}}(\boldsymbol{0})  \left (1-\| \boldsymbol{m}\|_2^2 \right )^{\frac{N-r-2}{2}} \E\left [  |\det \boldsymbol{H}_{\boldsymbol{m}} | \, \mathbf{1}_{\left \{i(\boldsymbol{H}_{\boldsymbol{m}})=\ell \right\}}\, \mathbf{1}_{\{f_{\boldsymbol{m}} \in B\}}\right ],
\end{split}
\end{equation}
where \(f_{\boldsymbol{m}} \stackrel{d}{=} f_N(\boldsymbol{e}_N)\), \(\boldsymbol{g}_{\boldsymbol{m}} \stackrel{d}{=} \nabla f_N(\boldsymbol{e}_N)\), and \(\boldsymbol{H}_{\boldsymbol{m}} \stackrel{d}{=} \nabla^2 f_N(\boldsymbol{e}_N)\). Taking \(\ell=N-1\) yields~\eqref{Kac--Rice index 3}, and summing~\eqref{kac rice index 2} over \(\ell\) gives~\eqref{Kac--Rice tot 3}.
\end{proof}

\subsection{Large deviations of finite-rank spiked GOE matrices} \label{LDP}

This section is devoted to the large deviations of the largest eigenvalue of a finite-rank spiked GOE matrix. For a positive integer \(k\), let \(\gamma_1 \geq \cdots \geq \gamma_k \ge 0 > \gamma_{k+1} \geq \cdots \geq \gamma_r\). We then define the rank-\(r\) spiked GOE matrix \(\boldsymbol{X}_N^{\boldsymbol{\gamma}}\) by
\[
\boldsymbol{X}_N^{\boldsymbol{\gamma}} = \sum_{i=1}^r \gamma_i \boldsymbol{e}_i \boldsymbol{e}_i^\top + \boldsymbol{W}_N,
\] 
where \(\boldsymbol{W}_N \sim \text{GOE}(N)\) and \((\boldsymbol{e}_1, \dots, \boldsymbol{e}_r)\) is a family of orthonormal vectors distributed according to the Haar measure on the orthogonal group (equivalently, a uniformly random orthonormal frame). Let \(\lambda_1 \geq \dots \geq \lambda_N\) denote the \(N\) eigenvalues of \(\boldsymbol{X}_N^{\boldsymbol{\gamma}} \). The following lemma specifies the joint distribution of the eigenvalues of the spiked matrix model \(\boldsymbol{X}_N^{\boldsymbol{\gamma}}\) (see e.g.~\cite{MR2336602, MR4436026}).

\begin{lem} \label{joint density eig}
The joint density of the ordered eigenvalues of \(\boldsymbol{X}_N^{\boldsymbol{\gamma}}\) is given by
\begin{equation}\label{eq joint density eig}
\mathbb{P}^{\boldsymbol{\gamma}}_N(\textnormal{d} \lambda_1, \dots, \textnormal{d} \lambda_N) = \frac{1}{Z^{\boldsymbol{\gamma}}_N} \prod_{i<j} |\lambda_i - \lambda_j| I_N(\boldsymbol{\gamma}, \boldsymbol{D}_N (\boldsymbol{\lambda}) )  e^{-\frac{N}{4} \sum{i=1}^N \lambda_i^2} \textnormal{d} \lambda_1 \cdots \textnormal{d} \lambda_N,\end{equation}
where \(\boldsymbol{D}_N (\boldsymbol{\lambda}) \coloneqq \textnormal{diag} (\lambda_1, \ldots, \lambda_N)\) denotes the \(N \times N\) diagonal matrix with entries given by \(\lambda_1, \ldots, \lambda_N\), and \(I_N(\boldsymbol{\gamma}, \boldsymbol{D}_N (\boldsymbol{\lambda}))\) is the finite-rank spherical integral given by
\[
I_N(\boldsymbol{\gamma}, \boldsymbol{D}N (\boldsymbol{\lambda}) ) \coloneqq \int_{\mathcal{O}_N} \exp \left \{ \frac{N}{2}  \Tr \left ( \boldsymbol{U}  \boldsymbol{D}N (\boldsymbol{\lambda}) \boldsymbol{U}^\top \sum_{i=1}^r \gamma_i \boldsymbol{e}_i \boldsymbol{e}_i^\top\right) \right \} \textnormal{d} m_N (\boldsymbol{U}),
\]
with \(m_N\) the Haar probability measure on the orthogonal group \(\mathcal{O}_N\) of size \(N\).
\end{lem}

The large deviation principle (LDP) for the largest eigenvalue of a rank-one deformation of a Gaussian Wigner matrix was first established by Ma\"{i}da~\cite{MR2336602}. Guionnet and Husson~\cite[Proposition 2.7]{MR4436026} subsequently proved a LDP for the joint law of extreme eigenvalues of finite-rank deformations. Their proof relies on the factorization of finite-rank spherical integrals into a product of rank-one spherical integrals in the large-\(N\) limit. We also highlight the recent work of Husson and Ko~\cite{MR4970431}, which generalizes the results for finite-rank spherical integrals from~\cite{MR4436026} to spherical integrals of sub-linear rank. For the purposes of this paper, we focus on the result from~\cite{MR4436026} for the \(k\) largest eigenvalues of \(\boldsymbol{X}_N^{\boldsymbol{\gamma}}\).

\begin{lem}[\cite{MR4436026}] \label{LDP eig}
The joint law of \((\lambda_1, \dots, \lambda_k)\) satisfies a LDP with speed \(N\) and good rate function \(I_{\boldsymbol{\gamma}}\) given by
\[
I_{\boldsymbol{\gamma}}(x_1, \ldots, x_k) = 
\begin{cases}
 \sum_{i=1}^k I_{\gamma_i}(x_i) & \textnormal{if} \enspace x_1 \geq \dots \geq x_k \geq 2, \\ 
+ \infty & \textnormal{otherwise} .
\end{cases}
\]
Here, for every \(\gamma > 0\) and \(x \ge 2\),
\[
I_\gamma (x) = I_{\textnormal{GOE}}(x) - \frac{1}{2}J(\mu_{\textnormal{sc}}, \gamma ,x) -  \inf_{y \ge 2} \left \{ I_{\textnormal{GOE}}(y) -\frac{1}{2}J (\mu_{\textnormal{sc}},\gamma ,y) \right \},
\]
where
\[
\begin{split}
I_{\textnormal{GOE}}(x) & = \frac{1}{2}\int_2^x \sqrt{y^2-4} \, \textnormal{d}y ,\\
J(\mu_{\textnormal{sc}}, \gamma ,x) & = 
\begin{cases}
\gamma x - 1 - \log(\gamma) - \Omega(x)  & \textnormal{if} \enspace G_{\mu_{\textnormal{sc}}}(x) \leq \gamma,\\
\frac{1}{2}\gamma^2 & \textnormal{if} \enspace G_{\mu_{\textnormal{sc}}}(x) > \gamma.
\end{cases}
\end{split}
\]
The function \(\Omega(x)\) is given by~\eqref{def: function Phi star}. Moreover, \(\mu_{\textnormal{sc}}\) denotes the semicircle distribution, and \(G_{\mu_{\textnormal{sc}}}(x) = \frac{1}{2}(x-\sqrt{x^2-4})\) its Cauchy--Stieltjes transform.
\end{lem}

In other words, Lemma~\ref{LDP eig} says that for any closed subset \(F \subset \R^k\), 
\[
 \limsup_{N \to \infty} \frac{1}{N} \log \mathbb{P}^{\boldsymbol{\gamma}}_N\left( (\lambda_1, \ldots, \lambda_k) \in F \right ) \leq - \inf_F I_{\boldsymbol{\gamma}},
 \]
and for any open subset \(O \subset \R^k\), 
\[ 
\liminf_{N \to \infty} \frac{1}{N} \log \mathbb{P}^{\boldsymbol{\gamma}}_N \left( (\lambda_1, \ldots, \lambda_k) \in O \right ) \geq - \inf_O I_{\boldsymbol{\gamma}}. 
\]

\begin{proof}
The result follows directly from~\cite[Proposition 2.7]{MR4436026} combined with the contraction principle (e.g., see~\cite[Theorem 4.2.1]{MR2571413}). By~\cite[Proposition 2.7]{MR4436026}, the joint law of the \(k\) largest and \(r-k\) smallest eigenvalues of \(\boldsymbol{X}_N^{\boldsymbol{\gamma}}\) satisfies a LDP with speed \(N\) and good rate function 
\[
\begin{split}
& I (x_1, \ldots, x_r) \\
&= 
\begin{cases}
 \sum_{i=1}^k I_{\gamma_i}(x_i) +  \sum_{i=k+1}^r I_{-\gamma_i}(-x_i)  & \text{if } x_1 \geq \dots \geq x_k \geq 2 \geq -2 \geq x_{k+1} \geq \dots \geq x_r, \\ 
+\infty & \text{otherwise}.
\end{cases}
\end{split}
\]
To derive the LDP for the joint law of the \(k\) largest eigenvalues, we apply the contraction principle, which implies that the joint law of the \(k\) largest eigenvalues satisfies a LDP with speed \(N\) and good rate function
\[
\begin{split}
I_{\boldsymbol{\gamma}} (x_1, \ldots, x_k) & = \inf_{x_{k+1}, \ldots, x_r} I(x_1, \ldots, x_r) \\
& = 
\begin{cases}
 \sum_{i=1}^k I_{\gamma_i}(x_i) + \sum_{i=k+1}^r \inf_{x_i \le -2} I_{-\gamma_i} (-x_i) & \text{if} \enspace x_1 \geq \dots \geq x_k \geq 2, \\ 
+ \infty & \text{otherwise} .
\end{cases}
\end{split}
\]
The statement follows since \(\inf_{x_i \le -2} I_{-\gamma_i} (-x_i)  = 0\) for all \(k+1 \le i \le r\). 
\end{proof}

We now derive the LDP for the largest eigenvalue of \(\boldsymbol{X}_N^{\boldsymbol{\gamma}}\) using Lemma~\ref{LDP eig} and the contraction principle. In particular, when \(r=1\), the following result reduces to~\cite[Theorem 3.2]{MR2336602}.

\begin{lem} \label{LDP max eig}
The law of the largest eigenvalue of \(\boldsymbol{X}^{\boldsymbol{\gamma}}_N\) satisfies a LDP with speed \(N\) and good rate function \(I_{\max}\) defined as follows.
\begin{enumerate}
\item If \(\gamma_1 \geq \gamma_2 \geq \cdots \geq \gamma_k > 1\), then
\begin{equation} \label{GRF 1}
I_{\max}(x)=
\begin{cases}
I_{\gamma_1}(x) & \textnormal{if} \enspace x \geq \gamma_1+\frac{1}{\gamma_1}, \\
\sum_{i=1}^j I_{\gamma_i}(x) & \textnormal{if} \enspace \gamma_{j+1}+ \frac{1}{\gamma_{j+1}}\leq x < \gamma_j+\frac{1}{\gamma_j}, \enspace \textnormal{for} \enspace j=1, \ldots, k-1,\\
\sum_{i=1}^k I_{\gamma_i}(x) & \textnormal{if} \enspace 2 \leq x < \gamma_k+\frac{1}{\gamma_k}, \\
+ \infty & \textnormal{if} \enspace x < 2,
\end{cases}
\end{equation}
where \(I_\gamma(x)\) is given by
\begin{equation}\label{I gamma}
I_{\gamma}(x) = \frac{1}{4} \int^x_{\gamma + \frac{1}{\gamma}} \sqrt{y^2-4} \textnormal{d}y - \frac{1}{2} \gamma \left ( x - \left ( \gamma + \frac{1}{\gamma}\right ) \right ) + \frac{1}{8} \left ( x^2 - \left (\gamma + \frac{1}{\gamma} \right )^2 \right ).
\end{equation}
\item If \(1 \geq \gamma_1  \geq \cdots \geq \gamma_k > 0\), then
\begin{equation} \label{GRF 2}
I_{\max}(x)= 
\begin{cases}
\frac{1}{4} \int_2^x \sqrt{y^2-4} \textnormal{d}y +\frac{1}{8}x^2 - \frac{1}{2} \gamma_1 x + \frac{1}{4} + \frac{1}{2} \log (\gamma_1) + \frac{1}{4} \gamma_1^2 & \textnormal{if} \enspace x >  \gamma_1 + \frac{1}{\gamma_1}, \\
\frac{1}{2} \int_2^x \sqrt{y^2-4} \textnormal{d}y & \textnormal{if} \enspace 2 \leq x \leq \gamma_1 + \frac{1}{\gamma_1}, \\
+ \infty & \textnormal{if} \enspace x < 2.
\end{cases}
\end{equation}
\item If \(\gamma_1 \geq \cdots \geq \gamma_\ell > 1 \ge \gamma_{\ell+1} \geq \cdots \geq \gamma_k>0\) for some \(\ell \in \{1, \ldots, k-1\}\), then
\begin{equation}  \label{GRF 3}
I_{\max}(x)=
\begin{cases}
 I_{\gamma_1}(x) & \textnormal{if} \enspace x \geq \gamma_1+\frac{1}{\gamma_1}, \\
\sum_{i=1}^j I_{\gamma_i}(x) & \textnormal{if} \enspace \gamma_{j+1}+ \frac{1}{\gamma_{j+1}}\leq x < \gamma_j+\frac{1}{\gamma_j}, \enspace \textnormal{for} \enspace j=1, \ldots, \ell-1,\\
\sum_{i=1}^\ell I_{\gamma_i}(x) & \textnormal{if} \enspace  2 \leq x < \gamma_\ell+\frac{1}{\gamma_\ell}, \\
+ \infty & \textnormal{if} \enspace x < 2,
\end{cases}
\end{equation}
where \(I_{\gamma_i}\) is given by~\eqref{I gamma}.
\end{enumerate}
\end{lem}

\begin{proof}
As the \(\max\) function is continuous, the contraction principle (e.g., see~\cite[Theorem 4.2.1]{MR2571413}) implies that the largest eigenvalue \(\lambda_1\) satisfies a LDP with speed \(N\) and good rate function \(I_{\max}\). For any \(x \ge 2\), \(I_{\max}\) is given by
\[
\begin{split}
I_{\max}(x) & = \inf \left \{I_{\boldsymbol{\gamma}}(x_1, \ldots, x_k) \colon \max_{1 \leq i \leq k} x_i = x \right \} \\
& = \inf \left \{  \sum_{i=1}^k I_{\gamma_i}(x_i) \colon x_1 \geq \dots \geq x_k \geq 2, x_1 = x \right \}\\
&=  I_{\gamma_1}(x) +  \sum_{i=2}^k \inf_{2 \leq y \leq x} I_{\gamma_i}(y).
\end{split}
\]
In particular, for any \(x \in \R\), we have
\[
I_{\max}(x) = 
\begin{cases}
\displaystyle I_{\gamma_1}(x) +  \sum_{i=2}^k \inf_{2 \leq y \leq x} I_{\gamma_i}(y) & \textnormal{if} \enspace x \ge 2,\\[2mm]
+\infty & \textnormal{if} \enspace x < 2.
\end{cases}
\]
We now specify the function \(I_{\max}\) for the three different cases. We first note that
\begin{equation} \label{integral}
\int_2^x \sqrt{y^2-4} \, \textnormal{d}y = \frac{1}{2} x^2 - 1 - 2 \Omega(x).
\end{equation}
\begin{enumerate}
\item
If \(\gamma_1 \geq \cdots \geq \gamma_k > 1\), then for any \(x \geq 2\), \(G_{\mu_{\textnormal{sc}}}(x) < \gamma_i\) for all \(1 \leq i \leq k\). From Lemma~\ref{LDP eig} and~\eqref{integral}, we have 
\[
\begin{split}
J(\mu_{\textnormal{sc}}, \gamma_i,x) &= \gamma_i x - 1 - \log(\gamma_i) - \Omega(x) \\
&= \frac{1}{2} \int_2^x \sqrt{y^2-4} \textnormal{d}y + \gamma_i x - \frac{1}{4} x^2- \frac{1}{2} - \log(\gamma_i).
\end{split}
\]
As a result, we obtain
\begin{equation} \label{GRF goe}
I_{\textnormal{GOE}}(x) - \frac{1}{2}J(\mu_{\textnormal{sc}}, \gamma_i,x) = \frac{1}{4} \int_2^x \sqrt{y^2-4} \, \textnormal{d}y - \frac{1}{2}\gamma_i x + \frac{1}{8}x^2 + \frac{1}{4} + \frac{1}{2} \log(\gamma_i).
\end{equation}
Differentiating this function on \([2,\infty)\), we see that it is decreasing on \([2, \gamma_i + \frac{1}{\gamma_i})\) and then increasing, and its infimum is reached at \( \gamma_i + \frac{1}{\gamma_i}\). This gives \(I_{\gamma_i}\) as in~\eqref{I gamma} and it is straightforward to see that the good rate function \(I_{\max}(x)\) is given as in~\eqref{GRF 1}. 

\item For the second assertion, if \(1 \geq \gamma_1 \geq \cdots \geq \gamma_k > 0\), then on \([2, \gamma_i + \frac{1}{\gamma_i}]\), \(G_{\mu_{\textnormal{sc}}}(x) \geq \gamma_i\), and on \([ \gamma_i + \frac{1}{\gamma_i}, \infty)\),  \(G_{\mu_{\textnormal{sc}}}(x) \leq \gamma_i\), Therefore, we obtain
\[
I_{\textnormal{GOE}}(x) - \frac{1}{2}J(\mu_{\textnormal{sc}}, \gamma_i,x) =  
\begin{cases}
I_1(x) & \text{if} \enspace 2 \leq x \leq \gamma_i+\frac{1}{\gamma_i},\\
I_2(x) & \text{if} \enspace x > \gamma_i+\frac{1}{\gamma_i},\\
\end{cases}
\]
where 
\[
I_1 (x) = \frac{1}{2} \int_2^x \sqrt{y^2-4} \, \textnormal{d}y -\frac{1}{4} \gamma_i^2,
\]
and \(I_2\) is given as in~\eqref{GRF goe}. Both functions \(I_1\) and \(I_2\) are increasing, and the infimum of \(I_{\text{GOE}}(x) - \frac{1}{2}J(\mu_{\textnormal{sc}}, \gamma_i,x)\) over all \(x\) is reached at \(2\) and is equal to \(-\frac{1}{4} \gamma_i^2\). Therefore, \(I_{\gamma_i}\) equals \( \frac{1}{2} \int_2^x \sqrt{y^2-4} \, \textnormal{d}y\) on \([2, \gamma_i + \frac{1}{\gamma_i}] \) and on \((\gamma_i + \frac{1}{\gamma_i}, \infty )\) is given by
\[
I_{\gamma_i} = \frac{1}{4} \int_2^x \sqrt{y^2-4} \textnormal{d}y +\frac{1}{8}x^2 - \frac{1}{2} \gamma_i x + \frac{1}{4} + \frac{1}{2} \log (\gamma_i) + \frac{1}{4} \gamma_i^2.
\]
 Then, for \(x \geq 2\), the good rate function \(I_{\max}\) is given by
\[
I_{\max}(x) = I_{\gamma_1}(x) +  \sum_{i=2}^k \inf_{2 \leq x_i \leq x} I_{\gamma_i}(x_i) = I_{\gamma_1}(x),
\]
since for all \(i \geq 2\), the infimum of \(I_{\gamma_i}(x_i)\) on \([2, x]\) is zero. 

\item Finally, the third assertion is a combination of the previous two. For \(\ell+1 \leq i \leq k\), the infimum \(I_{\gamma_i}(x_i)\) over \(x_i \in [2,x]\) equals zero, and the argument is the same as for the first assertion.
\end{enumerate}
\end{proof}

Recall the function \(L\) from Definition~\ref{def: function L}. We readily obtain the following corollary. 

\begin{cor}\label{cor LDP max eig}
For every \(t \in \R\), it holds that
\begin{align*}
\limsup_{N \to \infty} \frac{1}{N} \log \mathbb{P}^{\boldsymbol{\gamma}}_N \left( \lambda_1 \leq t \right ) &\leq - L(\boldsymbol{\gamma}, t),\\
\liminf_{N \to \infty} \frac{1}{N} \log \mathbb{P}^{\boldsymbol{\gamma}}_N \left( \lambda_1  < t \right ) &\geq - L(\boldsymbol{\gamma}, t_{-}),
\end{align*}
where \(L(\boldsymbol{\gamma}, t_{-}) = \lim_{t' \to t, t' < t} L(\boldsymbol{\gamma}, t')\). Note that \(L(\boldsymbol{\gamma}, t_{-})\) differs from \(L(\boldsymbol{\gamma}, t)\) only when \(t = 2\).
\end{cor}

\begin{proof}
By Lemma~\ref{LDP max eig}, for every \(\gamma_1 \ge\cdots \ge \gamma_k \ge 0 > \gamma_{k+1} \ge \cdots \ge \gamma_r\) and every \(t \in \R\),
\[
L (\boldsymbol{\gamma}, t) = \inf_{x \le t} I_{\max}(x).
\]
If \(1 \ge \gamma_1 \ge \cdots \ge \gamma_k >0\), then according to~\eqref{GRF 2}, 
\[
L (\boldsymbol{\gamma}, t) = 
\begin{cases}
0 & \text{if} \enspace t \ge  2,\\
+\infty & \text{if} \enspace t < 2.
\end{cases}
\]
Otherwise, there exists some \(\ell \in \{1, \ldots, k\}\) such that \(\gamma_1 \geq \cdots \geq \gamma_\ell > 1 \ge \gamma_{\ell+1} \geq \cdots \geq \gamma_k \). Then, according to~\eqref{GRF 3},
\[
L (\boldsymbol{\gamma}, t) = 
\begin{cases}
0 & \text{if} \enspace t \ge  \gamma_1 + \frac{1}{\gamma_1},\\
I_{\gamma_1}(t) & \text{if} \enspace \gamma_2 + \frac{1}{\gamma_2} \le t <  \gamma_1 + \frac{1}{\gamma_1}, \\
I_{\gamma_1}(t) + I_{\gamma_2}(t) & \text{if} \enspace \gamma_3 + \frac{1}{\gamma_3} \le t <  \gamma_2 + \frac{1}{\gamma_2}, \\
\vdots \\
\sum_{i=1}^{\ell-1} I_{\gamma_i}(t) & \text{if} \enspace \gamma_\ell + \frac{1}{\gamma_\ell } \le t <  \gamma_{\ell-1} + \frac{1}{\gamma_{\ell-1}},\\
\sum_{i=1}^\ell I_{\gamma_i}(t) & \text{if} \enspace 2 \le t <  \gamma_\ell + \frac{1}{\gamma_\ell},\\
+\infty & \text{if} \enspace t < 2.
\end{cases}
\] 
This corresponds to Definition~\ref{def: function L}.
\end{proof}

\section{Proofs of the main results} \label{proof main results}

In this section, we prove Theorems~\ref{thm crit points} and~\ref{thm local maxima}. Having the characterization of the landscape complexity for all integers \(N\) at hand (see Proposition~\ref{formula complexity}), the next important step is to study the exponential asymptotics of the determinant of the random matrix \(\boldsymbol{H}_{N-1}\) given by~\eqref{matrix H_N}.

\subsection{Preliminary remarks} \label{preliminaries}

We begin by identifying the eigenvalues of the finite-rank matrix \(\boldsymbol{P}_{N-1}(\boldsymbol{m})\) appearing in the conditional law of the Hessian in Proposition~\ref{formula complexity}. Recall from~\eqref{def: P matrix} that, for \(\boldsymbol{m}\in\mathcal{D}_r\),
\[
\boldsymbol{P}_{N-1}(\boldsymbol{m}) = \sum_{i=1}^r \theta_{i}(m_i) \boldsymbol{v}_i(\boldsymbol{m}) \boldsymbol{v}_i(\boldsymbol{m})^\top,
\]
where \(\theta_i\) is defined in~\eqref{function theta}, and the vectors \(\boldsymbol{v}_1(\boldsymbol{m}),\ldots,\boldsymbol{v}_r(\boldsymbol{m}) \in\mathbb{S}^{N-2}\) are defined via their Gram matrix \(\boldsymbol{G}(\boldsymbol{m})\) given in~\eqref{overlap v}. Recall also that \( \boldsymbol{\Theta}(\boldsymbol{m}) = \operatorname{diag}\bigl( \theta_1(m_1),\ldots,\theta_r(m_r) \bigr)\), and that \(\gamma_1(\boldsymbol{m}) \geq\cdots\geq \gamma_r(\boldsymbol{m}) \) are the ordered eigenvalues of \(\boldsymbol{G}(\boldsymbol{m})^{1/2} \boldsymbol{\Theta}(\boldsymbol{m}) \boldsymbol{G}(\boldsymbol{m})^{1/2}\), as introduced in Definition~\ref{def: perturbation matrix}. The following lemma identifies these quantities with the eigenvalues of \(\boldsymbol{P}_{N-1}(\boldsymbol{m})\).

\begin{lem} \label{lem: spectrum P}
For every \(\boldsymbol{m}\in\mathcal{D}_r\) and \(N>r\), the spectrum of \(\boldsymbol{P}_{N-1}(\boldsymbol{m})\), counted with algebraic multiplicity, consists of \(\gamma_1(\boldsymbol{m}),\ldots,\gamma_r(\boldsymbol{m})\) and \(N-1-r\) additional zero eigenvalues.
\end{lem}

\begin{proof}
Let \(\boldsymbol{V}(\boldsymbol{m}) = [\boldsymbol{v}_1(\boldsymbol{m}), \ldots, \boldsymbol{v}_r(\boldsymbol{m}) ]\). Then \(\boldsymbol{V}(\boldsymbol{m})^\top \boldsymbol{V}(\boldsymbol{m}) = \boldsymbol{G}(\boldsymbol{m})\) and 
\[
\boldsymbol{P}_{N-1}(\boldsymbol{m}) =
\boldsymbol{V}(\boldsymbol{m})
\boldsymbol{\Theta}(\boldsymbol{m})
\boldsymbol{V}(\boldsymbol{m})^\top.
\]
By Sylvester's determinant identity,
\[
\begin{split}
\det \left( z\boldsymbol{I}_{N-1} - \boldsymbol{P}_{N-1}(\boldsymbol{m}) \right) & = z^{N-1-r} \det\left( z\boldsymbol{I}_r - \boldsymbol{\Theta}(\boldsymbol{m}) \boldsymbol{V}(\boldsymbol{m})^\top \boldsymbol{V}(\boldsymbol{m}) \right) \\
& = z^{N-1-r} \det\left( z\boldsymbol{I}_r - \boldsymbol{\Theta}(\boldsymbol{m}) \boldsymbol{G}(\boldsymbol{m})\right).
\end{split}
\]
Since \(\boldsymbol{G}(\boldsymbol{m})\) is positive definite, \(\boldsymbol{\Theta}(\boldsymbol{m})\boldsymbol{G}(\boldsymbol{m})\) is similar to \(\boldsymbol{G}(\boldsymbol{m})^{1/2} \boldsymbol{\Theta}(\boldsymbol{m}) \boldsymbol{G}(\boldsymbol{m})^{1/2}\). Its eigenvalues are therefore \(\gamma_1(\boldsymbol{m}),\ldots,\gamma_r(\boldsymbol{m})\), which proves the claim.
\end{proof}

Since the entries of \(\boldsymbol{G}(\boldsymbol{m})\) and \(\boldsymbol{\Theta}(\boldsymbol{m})\) depend continuously on \(\boldsymbol{m}\in\mathcal{D}_r\), and since the matrix square-root map is continuous on the cone of positive-definite matrices, the symmetric matrix \(\boldsymbol{G}(\boldsymbol{m})^{1/2} \boldsymbol{\Theta}(\boldsymbol{m}) \boldsymbol{G}(\boldsymbol{m})^{1/2}\) depends continuously on \(\boldsymbol{m}\). The continuity of the ordered eigenvalues of a symmetric matrix therefore implies that \(\boldsymbol{m}\mapsto\gamma_i(\boldsymbol{m})\) are continuous on $\mathcal{D}_r$.

It follows from Lemma~\ref{lem: spectrum P} that, for every \(\boldsymbol{m}\in\mathcal{D}_r\), there exists an orthogonal matrix \(\boldsymbol{Q}_{\boldsymbol{m}}\in\mathcal{O}(N-1)\) such that
\[ 
\boldsymbol{Q}_{\boldsymbol{m}}^\top \boldsymbol{P}_{N-1}(\boldsymbol{m}) \boldsymbol{Q}_{\boldsymbol{m}} = \sum_{i=1}^r \gamma_i(\boldsymbol{m}) \boldsymbol{e}_i\boldsymbol{e}_i^\top,
\]
where \(\boldsymbol{e}_1, \ldots, \boldsymbol{e}_r\) denote the standard basis vectors of \(\mathbb{R}^{N-1}\). By the invariance of the GOE under orthogonal conjugation,
\[ 
\boldsymbol{Q}_{\boldsymbol{m}}^\top \boldsymbol{W}_{N-1} \boldsymbol{Q}_{\boldsymbol{m}} \stackrel{d}{=} \boldsymbol{W}_{N-1}. 
\]
Consequently,~\eqref{matrix H_N} gives
\begin{equation}\label{Hess}
\boldsymbol{H}_{N-1}(\boldsymbol{m},x) \stackrel{d}{=} \boldsymbol{W}_{N-1} + \sqrt{\frac{N}{N-1}}\sum_{i=1}^r \gamma_{i}(\boldsymbol{m}) \boldsymbol{e}_i \boldsymbol{e}_i^\top -\sqrt{\frac{N}{N-1}} t(\boldsymbol{m},x) \boldsymbol{I}_{N-1}.
\end{equation} 
Finally, for any sequence of real numbers \(\boldsymbol{\gamma} = (\gamma_1, \ldots, \gamma_r)\), define the finite-rank deformation
\begin{equation} \label{eq: J_N}
\boldsymbol{X}_{N-1}^{\boldsymbol{\gamma}} \coloneqq \boldsymbol{W}_{N-1} + \sum_{i=1}^r \gamma_{i} \boldsymbol{e}_i \boldsymbol{e}_i^\top.
\end{equation} 
With this notation, 
\[ 
\boldsymbol{H}_{N-1}(\boldsymbol{m},x) \stackrel{d}{=} \boldsymbol{X}_{N-1}^{ \sqrt{N/(N-1)}\,\boldsymbol{\gamma}(\boldsymbol{m}) } - \sqrt{\frac{N}{N-1}} t(\boldsymbol{m},x)\boldsymbol{I}_{N-1}. 
\]

\textbf{Notations.} 
Throughout, we write \(\lVert \cdot \rVert\) for the operator norm on elements of \(\R^{N \times N}\) induced by the \(L^2\)-distance on \(\R^N\). For test functions \(f \colon \R \to \R\), we let  \( \lVert f \rVert_{L^\infty} = \sup_x  |f(x)|\) and \(\lVert f \rVert_{\textnormal{Lip}} = \sup_{x \neq y} \lvert \frac{f(x)-f(y)}{x-y} \rvert\). We let \(\mathcal{P}(\R)\) denote the space of probability measures on \(\R\) and we consider the following two distances on probability measures on \(\R\), called bounded-Lipschitz and Wasserstein-1, respectively: for every \(\mu,\nu \in \mathcal{P}(\R)\), 
\begin{align*}
d_{\text{BL}}(\mu, \nu) & = \sup \left \{ \left |\int_\R f \textnormal{d}(\mu-\nu) \right | \colon \lVert f \rVert_{\textnormal{Lip}} + \lVert f \rVert_{L^\infty} \leq 1 \right \},\\
W_1(\mu, \nu) & = \sup \left \{ \left | \int_\R f \textnormal{d}(\mu-\nu) \right | \colon \lVert f \rVert_{\textnormal{Lip}} \leq 1 \right \}.
\end{align*}
If \((X,d)\) is a metric space, we denote by \(\textnormal{B}(x,r)= \{y \in X \colon d(y,x) < r\}\) the open ball of radius \(r>0 \) around \(x  \in X\). For an \(N\times N\) Hermitian matrix \(\boldsymbol{M}\), we write \(\lambda_1 (\boldsymbol{M}),\ldots, \lambda_N(\boldsymbol{M})\) for its eigenvalues and \(\hat{\mu}_{\boldsymbol{M}}=\frac{1}{N} \sum_{i=1}^N \delta_{\lambda_{i}(\boldsymbol{M})}\) for its empirical spectral measure. We let \(\E \left [\hat{\mu}_{\boldsymbol{M}} \right ]\) denote the mean spectral measure, i.e., \(\int_\R f \textnormal{d} \E \left [\hat{\mu}_{\boldsymbol{M}} \right ] \coloneqq \E\left [ \int_\R f \textnormal{d}\hat{\mu}_{\boldsymbol{M}}\right ]\). We denote the semicircle law as \(\mu_{\textnormal{sc}}(\textnormal{d}x) = \frac{1}{2\pi}\sqrt{4-x^2} \mathbf{1}_{\{|x| \leq 2\}}\textnormal{d}x\).

\subsection{Proof of complexity result for total critical points} \label{subsection crit points}

Here, we provide the proof of Theorem~\ref{thm crit points} on the complexity of total critical points. Our general strategy is to establish upper and lower bounds, along with exponential tightness. The following three results play a crucial role in proving Theorem~\ref{thm crit points}.

\begin{lem}[Good rate function] \label{lem: GRF}
The function \(- \Sigma^{\textnormal{tot}}\) given by Definition~\ref{def: function Sigma tot} is a good rate function. Moreover, the function \(S\) given by~\eqref{function S} is continuous in \(\mathcal{D}_r \times \R\) and the function \(\Omega\) given by~\eqref{def: function Phi star} is continuous in \(\R\).
\end{lem}

\begin{lem}[Exponential tightness] \label{lem: exp tight}
It holds that
\[
\lim_{z \to \infty} \limsup_{N \to \infty} \frac{1}{N} \log \E \left [\textnormal{Crt}_N^{\textnormal{tot}}([-1,1]^r, (-\infty, -z] \cup [z, \infty)) \right ]= - \infty.
\]
\end{lem}

\begin{lem}[Determinant asymptotics] \label{lem: det J_N}
For every compact product set \(\mathcal{U}  = \mathcal{U}_1 \times \cdots \times \mathcal{U}_r \subset \R^r\) and every compact set \(\mathcal{T} \subset \R\), it holds that
\begin{align}
\limsup_{N \to \infty} \sup_{(\boldsymbol{\gamma},t) \in \, \mathcal{U} \times \mathcal{T}} \frac{1}{N} \log  \E \left [ |\det (\boldsymbol{X}_N^{\boldsymbol{\gamma}} - t \boldsymbol{I}_N )| \right] & \leq  \sup_{t \in \,\mathcal{T}}  \Omega(t), \label{eq: upper bound det J_N}  \\ 
\liminf_{N \to \infty} \inf_{(\boldsymbol{\gamma},t) \in \, \mathcal{U} \times \mathcal{T}}  \frac{1}{N} \log  \E \left [ |\det (\boldsymbol{X}_N^{\boldsymbol{\gamma}} - t \boldsymbol{I}_N)| \right] & \geq \inf_{t \in \, \mathcal{T}}  \Omega(t) . \label{eq: lower bound det J_N}
\end{align}
\end{lem}

We postpone the proof of these lemmas towards the end of the subsection. We now turn to the proof of Theorem~\ref{thm crit points}. We first show the upper bound~\eqref{sup crit points}.

\begin{proof}[Proof of Theorem~\ref{thm crit points} (Upper bound)]
Let 
\[
L\coloneqq \sup_{(\boldsymbol{m},x)\in \overline{M} \times \overline{B}} \Sigma^{\textnormal{tot}}(\boldsymbol{m},x). 
\]
Suppose first that \(L=-\infty\). Then either \(B=\emptyset\) or \(M\cap\mathcal{D}_r=\emptyset\). Indeed, otherwise there would exist \((\boldsymbol{m},x)\in(M\cap\mathcal{D}_r)\times B\), for which \(\Sigma^{\textnormal{tot}}(\boldsymbol{m},x)\in\mathbb{R}\), contradicting the definition of \(L\). Hence the integration domain in Proposition~\ref{formula complexity} is empty, and therefore \(\E\left[\textnormal{Crt}_N^{\textnormal{tot}}(M,B)\right]=0\). Consequently,
\[
\limsup_{N\to\infty} \frac{1}{N} \log \E\left[\textnormal{Crt}_N^{\textnormal{tot}}(M,B)\right] =-\infty=L.
\]
Thus the desired upper bound holds in this case. We may therefore assume that \(L>-\infty\). By exponential tightness, see Lemma~\ref{lem: exp tight}, it suffices to prove the upper bound when \(B\) is bounded. In this case \(\overline{M}\times\overline{B}\) is compact.

For every \((\boldsymbol{m}_0,x_0)\in\overline{M}\times\overline{B}\), set
\[
M_{\delta_0} \coloneqq \prod_{i=1}^r \left ((m_0)_i-\delta_0,(m_0)_i+\delta_0 \right),
\qquad
B_{\delta_0}\coloneqq (x_0-\delta_0,x_0+\delta_0).
\]
We claim that
\begin{equation} \label{eq:local-upper}
\limsup_{\delta_0\to0_+} \limsup_{N\to\infty} \frac{1}{N}\log \E\left[\textnormal{Crt}_N^{\textnormal{tot}}(M_{\delta_0},B_{\delta_0}) \right] \le \Sigma^{\textnormal{tot}}(\boldsymbol{m}_0,x_0).
\end{equation}
We first explain how this local estimate implies the global upper bound. Assume that~\eqref{eq:local-upper} holds. Fix \(\varepsilon>0\). For every \((\boldsymbol{m},x)\in\overline{M}\times\overline{B}\), the estimate~\eqref{eq:local-upper} implies that there exists \(\delta_{\boldsymbol{m},x}>0\) such that, with
\[
M_{\boldsymbol{m},x} \coloneqq \prod_{i=1}^r(m_i-\delta_{\boldsymbol{m},x},m_i+\delta_{\boldsymbol{m},x}),
\qquad
B_{\boldsymbol{m},x} \coloneqq (x-\delta_{\boldsymbol{m},x},x+\delta_{\boldsymbol{m},x}),
\]
we have
\begin{equation} \label{eq:local-upper2}
\limsup_{N\to\infty} \frac{1}{N}\log \E \left[ \textnormal{Crt}_N^{\textnormal{tot}} (M_{\boldsymbol{m},x},B_{\boldsymbol{m},x}) \right] \le L+\varepsilon.
\end{equation}
Indeed, if \(\Sigma^{\textnormal{tot}}(\boldsymbol{m},x)>-\infty\), this follows directly from~\eqref{eq:local-upper} and the inequality \(\Sigma^{\textnormal{tot}}(\boldsymbol{m},x)\le L\); if \(\Sigma^{\textnormal{tot}}(\boldsymbol{m},x)=-\infty\), it follows from the definition of the \(\limsup\). The sets \(M_{\boldsymbol{m},x}\times B_{\boldsymbol{m},x}\), \((\boldsymbol{m},x)\in\overline{M}\times\overline{B}\), form an open cover of the compact set
\(\overline{M}\times\overline{B}\). Hence there exist \((\boldsymbol{m}_k,x_k)\in\overline{M}\times\overline{B}\), \(1\le k\le n\), such that
\[
\overline{M}\times\overline{B} \subseteq \bigcup_{k=1}^n \left(M_{\boldsymbol{m}_k,x_k}\times B_{\boldsymbol{m}_k,x_k}\right).
\]
Here the integer \(n\) is independent of \(N\). Therefore, using \eqref{Kac--Rice tot 2},
\[
\E \left[\textnormal{Crt}_N^{\textnormal{tot}}(M,B)\right] \le \E\left[\textnormal{Crt}_N^{\textnormal{tot}}(\overline{M},\overline{B})\right] \le \sum_{k=1}^n \E\left[ \textnormal{Crt}_N^{\textnormal{tot}} (M_{\boldsymbol{m}_k,x_k},B_{\boldsymbol{m}_k,x_k}) \right].
\]
Consequently,
\[
\begin{split}
\limsup_{N\to\infty} \frac{1}{N}\log \E\left[\textnormal{Crt}_N^{\textnormal{tot}}(M,B)\right] &\le \max_{1\le k\le n} \limsup_{N\to\infty} \frac{1}{N}\log \E\left[
\textnormal{Crt}_N^{\textnormal{tot}} (M_{\boldsymbol{m}_k,x_k},B_{\boldsymbol{m}_k,x_k}) \right]  \\ &\le L+\varepsilon,
\end{split}
\]
where the first inequality follows from~\cite[Lemma 1.2.15]{MR2571413} and the second from~\eqref{eq:local-upper2}. Since \(\varepsilon>0\) is arbitrary, this
proves the desired upper bound.

It remains to prove~\eqref{eq:local-upper}. Fix \((\boldsymbol{m}_0,x_0)\in\overline{M}\times\overline{B}\) and \(\delta_0>0\). If \(\overline{M}_{\delta_0}\cap \mathcal{D}_r=\emptyset\), then \(\E\left[\textnormal{Crt}_N^{\textnormal{tot}}(M_{\delta_0},B_{\delta_0})\right]=0\) for every \(N\). Thus the desired estimate is immediate in this case. We may therefore assume that \(\overline{M}_{\delta_0}\cap \mathcal{D}_r\neq\emptyset\). Let \(R_0\coloneqq \sup\{|x| \colon x\in\overline{B}_{\delta_0}\}\), which is a finite constant since \(B_{\delta_0}\) is bounded. By Proposition~\ref{formula complexity} and Remark~\ref{rmk: KR S0},
\[
\begin{split}
\E\left[\textnormal{Crt}_N^{\textnormal{tot}}(M_{\delta_0},B_{\delta_0})\right] & \le C_N \int_{(\boldsymbol{m},x)\in(\overline{M}_{\delta_0}\cap \mathcal{D}_r)\times\overline{B}_{\delta_0}}  \textnormal{d}  \boldsymbol{m} \textnormal{d} x                                      \\
& \quad \times \E\left[|\det \boldsymbol{H}_{N-1}(\boldsymbol{m},x)|\right] \exp\left\{ NS_0 (\boldsymbol{m},x) +\frac{N-r-2}{2}\log q(\boldsymbol{m}) \right\}.
\end{split}
\]
Note that the function \(S_0\) is bounded above on \((\overline{M}_{\delta_0}\cap \mathcal{D}_r)\times\overline{B}_{\delta_0}\). Since the integration domain is contained in \([-1,1]^r\times[-R_0,R_0]\), its volume is at most \(2^{r+1}R_0\). Hence
\begin{equation} \label{eq:upper-1-rigorous}
\begin{split}
\E \left[\textnormal{Crt}_N^{\textnormal{tot}}(M_{\delta_0}, B_{\delta_0})\right] &\le 2^{r+1}R_0C_N \sup_{(\boldsymbol{m},x)\in(\overline{M}_{\delta_0}\cap \mathcal{D}_r)\times\overline{B}_{\delta_0}} \E\left[|\det \boldsymbol{H}_{N-1}(\boldsymbol{m},x)|\right]  \\
&\quad \times \sup_{(\boldsymbol{m},x)\in(\overline{M}_{\delta_0}\cap \mathcal{D}_r)\times\overline{B}_{\delta_0}} \exp \left\{ NS_0 (\boldsymbol{m},x) + \frac{N-r-2}{2}\log q(\boldsymbol{m}) \right\}.
\end{split}
\end{equation}

We now bound the determinant term. Recall that
\[
\boldsymbol{P}_{N-1}(\boldsymbol{m}) = \sum_{i=1}^r \theta_i(m_i)\boldsymbol{v}_i(\boldsymbol{m}) \boldsymbol{v}_i(\boldsymbol{m})^\top.
\]
Since \(\|\boldsymbol{v}_i(\boldsymbol{m})\boldsymbol{v}_i(\boldsymbol{m})^\top\|_{\mathrm{op}}=1\), we have \(\|\boldsymbol{P}_{N-1}(\boldsymbol{m})\|_{\mathrm{op}} \le \sum_{i=1}^r |\theta_i(m_i)|\). Each function \(\theta_i\) is bounded on \([-1,1]\) (see~\eqref{function theta}). Therefore, there exists \(K<\infty\) such that \(|\gamma_i(\boldsymbol{m})|\le K\) for every \(\boldsymbol{m}\in \mathcal{D}_r\) and every \(1\le i\le r\). For \(\delta>0\), define
\[
\overline{\mathcal{U}}_\delta \coloneqq [-K-\delta,K+\delta]^r .
\]
Moreover, since \(t\) is continuous and \(\overline{M}_{\delta_0} \times \overline{B}_{\delta_0}\) is compact, the set 
\[
T_{\delta_0} \coloneqq \left\{ t(\boldsymbol{m},x) \colon  (\boldsymbol{m},x)\in \overline{M}_{\delta_0} \times\overline{B}_{\delta_0} \right\}
\]
is compact. For fixed \(\delta_0>0\) and \(\delta>0\), let
\[
\overline{\mathcal{T}}_\delta^{\delta_0} \coloneqq [\inf T_{\delta_0}-\delta,\sup T_{\delta_0}+\delta].
\]
Since \(\sqrt{N/(N-1)}\to1\), the uniform bound \(\sup_{\boldsymbol{m}\in \mathcal{D}_r} \|\boldsymbol{\gamma}(\boldsymbol{m})\|_\infty\leq K\) and the compactness of \(T_{\delta_0}\) imply that there exists \(N_{\delta,\delta_0}\) such that, for every \(N\geq N_{\delta,\delta_0}\),
\[
\sqrt{\frac N{N-1}}\boldsymbol{\gamma}(\boldsymbol{m}) \in \overline{\mathcal{U}}_\delta,
\qquad
\sqrt{\frac N{N-1}}t(\boldsymbol{m},x) \in \overline{\mathcal{T}}_\delta^{\delta_0},
\]
uniformly for \((\boldsymbol{m},x)\in (\overline{M}_{\delta_0}\cap \mathcal{D}_r)\times\overline{B}_{\delta_0}\). Therefore, by~\eqref{Hess}, for every \(N \ge N_{\delta,\delta_0}\),
\[
\sup_{(\boldsymbol{m},x)\in(\overline{M}_{\delta_0}\cap \mathcal{D}_r)\times\overline{B}_{\delta_0}} \E\left[|\det \boldsymbol{H}_{N-1}(\boldsymbol{m},x)|\right]  \\ \le \sup_{(\boldsymbol{\gamma},t)\in \overline{\mathcal{U}}_\delta\times \overline{\mathcal{T}}_\delta^{\delta_0}} \E\left[ \left| \det \left( \boldsymbol{X}_{N-1}^{\boldsymbol{\gamma}} - t\boldsymbol{I}_{N-1} \right) \right| \right].
\]
By Lemma~\ref{lem: det J_N}, in particular~\eqref{eq: upper bound det J_N},
\[
\limsup_{N\to\infty} \sup_{(\boldsymbol{\gamma},t)\in \overline{\mathcal{U}}_\delta\times\overline{\mathcal{T}}_\delta^{\delta_0}} \frac{1}{N-1} \log \E \left[ \left| \det \left( \boldsymbol{X}_{N-1}^{\boldsymbol{\gamma}} - t \boldsymbol{I}_{N-1} \right) \right| \right] \le \sup_{t\in\overline{\mathcal{T}}_\delta^{\delta_0}}\Omega(t).
\]
Consequently, for every \(\varepsilon>0\), there exists \(N_{\varepsilon, \delta,\delta_0} \ge N_{\delta,\delta_0}\) such that for all \(N \ge N_{\varepsilon, \delta,\delta_0}\),
\begin{equation} \label{eq:upper-2-rigorous}
\sup_{(\boldsymbol{m},x) \in(\overline{M}_{\delta_0}\cap \mathcal{D}_r)\times\overline{B}_{\delta_0}} \E\left[|\det \boldsymbol H_{N-1}(\boldsymbol{m},x)|\right] \le \exp \left \{ (N-1) \left( \sup_{t\in\overline{\mathcal{T}}_\delta^{\delta_0}}\Omega(t) + \varepsilon \right) \right\}.
\end{equation}
Combining~\eqref{eq:upper-1-rigorous} and~\eqref{eq:upper-2-rigorous}, taking \(\frac{1}{N}\log\), and using the fact that \(C_N\) is exponentially trivial, we obtain
\begin{equation} \label{eq:upper-3-rigorous}
\begin{split}
\limsup_{N\to\infty} \frac{1}{N}\log \E \left[\textnormal{Crt}_N^{\textnormal{tot}}(M_{\delta_0}, B_{\delta_0})\right]
& \le \limsup_{N\to\infty} \sup_{(\boldsymbol{m},x)\in(\overline{M}_{\delta_0}\cap \mathcal{D}_r)\times\overline{B}_{\delta_0}} \left[S_0 (\boldsymbol{m},x) + \frac{N-r-2}{2N}\log q(\boldsymbol{m}) \right]  \\
&\quad + \sup_{t\in\overline{\mathcal{T}}_\delta^{\delta_0}}\Omega(t) +\varepsilon .
\end{split}
\end{equation}

We now claim that
\[
\limsup_{N\to\infty} \sup_{(\boldsymbol{m},x)\in (\overline{M}_{\delta_0}\cap \mathcal{D}_r)\times\overline{B}_{\delta_0}} \left[S_0 (\boldsymbol{m},x) +\frac{N-r-2}{2N}\log q(\boldsymbol{m}) \right] \le \sup_{(\boldsymbol{m},x)\in(\overline{M}_{\delta_0}\cap \mathcal{D}_r)\times\overline{B}_{\delta_0}} S(\boldsymbol{m},x),
\]
where we recall that \(S (\boldsymbol{m},x)= S_0 (\boldsymbol{m},x) + \frac{1}{2} \log q(\boldsymbol{m})\). To prove the claim, fix \(\rho\in(0,1)\) and split \((\overline{M}_{\delta_0}\cap \mathcal{D}_r)\times\overline{B}_{\delta_0}\) according to whether \(q(\boldsymbol{m})\geq\rho\) or \(q(\boldsymbol{m})<\rho\). On the first set, we have
\[
S_0 (\boldsymbol{m},x) +\frac{N-r-2}{2N}\log q(\boldsymbol{m}) = S(\boldsymbol{m},x) -\frac{r+2}{2N}\log q(\boldsymbol{m}) \le S(\boldsymbol{m},x) -\frac{r+2}{2N}\log\rho.
\]
Taking \(N\to\infty\), this contribution is bounded above by
\[
\sup_{(\boldsymbol{m},x)\in(\overline{M}_{\delta_0}\cap \mathcal{D}_r)\times\overline{B}_{\delta_0}} S(\boldsymbol{m},x).
\]
On the complementary set \(q(\boldsymbol{m})<\rho\), let \(C<\infty\) be an upper bound for \(S_0\) on \((\overline{M}_{\delta_0}\cap \mathcal{D}_r)\times\overline{B}_{\delta_0}\). Since \((N-r-2)/(2N)\to1/2\), we have, for all sufficiently large \(N\), \((N-r-2)/(2N) \ge 1/4\). Since \(\log q(\boldsymbol{m}) < \log \rho < 0\), it follows that the contribution is bounded above, for all large \(N\), by \(C+\frac{1}{4}\log\rho\). Consequently,
\[
\begin{split}
&\limsup_{N\to\infty} \sup_{(\boldsymbol{m},x)\in (\overline{M}_{\delta_0}\cap \mathcal{D}_r)\times\overline{B}_{\delta_0}} \left[S_0 (\boldsymbol{m},x) +\frac{N-r-2}{2N}\log q(\boldsymbol{m}) \right] \\
&\le \max \left \{  \sup_{(\boldsymbol{m},x)\in(\overline{M}_{\delta_0}\cap \mathcal{D}_r)\times\overline{B}_{\delta_0}} S(\boldsymbol{m},x), C + \frac{1}{4} \log \rho \right \}.
\end{split}
\]
Letting \(\rho\downarrow0\) proves the claim. Therefore, by~\eqref{eq:upper-3-rigorous},
\[
\limsup_{N\to\infty} \frac{1}{N}\log \E\left[\textnormal{Crt}_N^{\textnormal{tot}}(M_{\delta_0}, B_{\delta_0})\right] \le \sup_{(\boldsymbol{m},x)\in(\overline{M}_{\delta_0}\cap \mathcal{D}_r)\times\overline{B}_{\delta_0}} S(\boldsymbol{m},x) + \sup_{t\in\overline{\mathcal{T}}_\delta^{\delta_0}}\Omega(t) +\varepsilon.
\]
We now let first \(\varepsilon\to 0^+\) and then \(\delta\to 0^+\). Since \(T_{\delta_0}\) is a compact interval, we obtain, by continuity of \(\Omega\), 
\[ 
\lim_{\delta\to0^+} \sup_{t\in\overline{\mathcal{T}}_\delta^{\delta_0}} \Omega(t) = \sup_{t\in T_{\delta_0}}\Omega(t) = \sup_{(\boldsymbol{m},x)\in \overline{M}_{\delta_0}\times\overline{B}_{\delta_0}} \Omega(t(\boldsymbol{m},x)). 
\]
Therefore,
\[
\limsup_{N\to\infty} \frac{1}{N}\log \E\left[\textnormal{Crt}_N^{\textnormal{tot}}(M_{\delta_0}, B_{\delta_0})\right] \le \sup_{\overline{M}_{\delta_0}\times\overline{B}_{\delta_0}} S + \sup_{\overline{M}_{\delta_0}\times\overline{B}_{\delta_0}} \Omega(t),
\]
where \(S\) is extended by \(-\infty\) on \(\mathcal{D}_r^{\textnormal{c}}\times\R\).

It remains to let \(\delta_0 \to 0^+\). To this end, set \(G(\boldsymbol{m},x)\coloneqq \Omega(t(\boldsymbol{m},x))\) so that \(\Sigma^{\textnormal{tot}}=S+G\). Since \(G\) is continuous, we define its oscillation on \(\overline{M}_{\delta_0}\times\overline{B}_{\delta_0}\) by
\[
\omega(\delta_0) \coloneqq \sup_{(\boldsymbol{m},x),(\boldsymbol{m}',x')\in \overline{M}_{\delta_0}\times\overline{B}_{\delta_0}} \left|G(\boldsymbol{m},x)-G(\boldsymbol{m}',x' )\right|.
\]
Then
\[
\sup_{\overline{M}_{\delta_0}\times\overline{B}_{\delta_0}} S + \sup_{\overline{M}_{\delta_0}\times\overline{B}_{\delta_0}} G \le \sup_{\overline{M}_{\delta_0}\times\overline{B}_{\delta_0}}\Sigma^{\textnormal{tot}} + \omega(\delta_0).
\]
Indeed, for any two points \((\boldsymbol{m},x),(\boldsymbol{m}',x')\in \overline{M}_{\delta_0}\times\overline{B}_{\delta_0}\),
\[
S(\boldsymbol{m},x)+G(\boldsymbol{m}',x') = \Sigma^{\textnormal{tot}}(\boldsymbol{m},x) + G(\boldsymbol{m}',x')-G(\boldsymbol{m},x) \le \sup_{\overline{M}_{\delta_0} \times\overline{B}_{\delta_0}}\Sigma^{\textnormal{tot}} + \omega(\delta_0),
\]
and taking the supremum gives the result. Hence
\[
\limsup_{N\to\infty} \frac{1}{N}\log \E\left[\textnormal{Crt}_N^{\textnormal{tot}}(M_{\delta_0}, B_{\delta_0})\right] \le \sup_{\overline{M}_{\delta_0}\times\overline{B}_{\delta_0}}\Sigma^{\textnormal{tot}} + \omega(\delta_0).
\]
Since \(G\) is continuous, \(\omega(\delta_0)\to0\) as \(\delta_0\to 0^+\). Moreover, by Lemma~\ref{lem: GRF}, \(-\Sigma^{\textnormal{tot}}\) is lower semicontinuous, and hence \(\Sigma^{\textnormal{tot}}\) is upper semicontinuous. Therefore,
\[
\limsup_{\delta_0 \to 0^+} \sup_{(\boldsymbol{m},x)\in\overline{M}_{\delta_0}\times\overline{B}_{\delta_0}} \Sigma^{\textnormal{tot}}(\boldsymbol{m},x) \le \Sigma^{\textnormal{tot}}(\boldsymbol{m}_0,x_0).
\]
Letting \(\delta_0\to 0^+\) proves~\eqref{eq:local-upper} and completes the proof of the upper bound.
\end{proof}

We next provide the proof of the lower bound~\eqref{inf crit points}.

\begin{proof}[Proof of Theorem~\ref{thm crit points} (Lower bound)]
It suffices to consider Borel sets \(M\subset [-1,1]^r, B \subset \R\) such that 
\[
\sup_{(\boldsymbol{m},x) \in M^\circ \times B^\circ} \Sigma^{\textnormal{tot}}(\boldsymbol{m}, x) > -\infty, 
\]
since otherwise the lower bound~\eqref{inf crit points} holds trivially. Fix \(\varepsilon_0>0\). By the definition of the supremum, there exists \((\boldsymbol{m}_0, x_0) \in M^\circ \times B^\circ\) such that
\begin{equation} \label{lower bound}
\Sigma^{\textnormal{tot}}(\boldsymbol{m}_0,x_0) \geq \sup_{(\boldsymbol{m}, x) \in M^\circ \times B^\circ}\Sigma^{\textnormal{tot}}(\boldsymbol{m}, x)-\varepsilon_0.
\end{equation}
Then \(\Sigma^{\textnormal{tot}}(\boldsymbol{m}_0,x_0)>-\infty\), hence \((\boldsymbol{m}_0,x_0)\in \mathcal{D}_r \times \R\). For \(\delta_0>0\), set
\[
M_{\delta_0} \coloneqq \prod_{i=1}^r ((m_0)_i-\delta_0,(m_0)_i+\delta_0) ,
\qquad
B_{\delta_0} \coloneqq (x_0-\delta_0,x_0+\delta_0).
\]
Since \((\boldsymbol{m}_0,x_0) \in (M^\circ\cap \mathcal{D}_r)\times B^\circ\) and \(M^\circ\cap \mathcal{D}_r\) and \(B^\circ\) are open, we may choose \(\delta_0>0\) sufficiently small that \(\overline{M}_{\delta_0}\subset M^\circ\cap \mathcal{D}_r\) and \(\overline{B}_{\delta_0}\subset B^\circ\). By Proposition~\ref{formula complexity}, 
\begin{equation} \label{lower bound 1}
\begin{split}
&\E \left [ \textnormal{Crt}_N^{\textnormal{tot}}(M, B) \right ] \\
&  \geq \E \left [ \textnormal{Crt}_N^{\textnormal{tot}}(M_{\delta_0}, B_{\delta_0})\right ]  \\
& =C_N \int_{M_{\delta_0} \times B_{\delta_0}} \textnormal{d} \boldsymbol{m}  \, \textnormal{d}x   \left(1-\|\boldsymbol{m}\|_2^2\right)^{-\frac{r+2}{2}}  \E\left [|\det \boldsymbol{H}_{N-1}(\boldsymbol{m},x)| \right ]  \exp \{N S(\boldsymbol{m},x) \} \\
& \geq (2\delta_0)^{r+1} C_N \inf_{(\boldsymbol{m},x) \in M_{\delta_0} \times B_{\delta_0}} \E \left [|\det \boldsymbol{H}_{N-1}(\boldsymbol{m},x)| \right ] \inf_{(\boldsymbol{m},x) \in M_{\delta_0} \times B_{\delta_0}} \exp \left \{N S(\boldsymbol{m},x)  \right \}.
\end{split}
\end{equation}
In the last inequality, we used the fact \(0<q(\boldsymbol{m}) = 1-\|\boldsymbol{m}\|_2^2 \leq1\) on \(\mathcal{D}_r\), and hence \(q(\boldsymbol{m})^{-\frac{r+2}{2}} \geq1\).

We now control the determinant term. Define
\[
U_{i,\delta_0} \coloneqq \{ \gamma_i(\boldsymbol{m}) \colon \boldsymbol{m} \in \overline{M}_{\delta_0} \},
\qquad
T_{\delta_0}  \coloneqq \{ t(\boldsymbol{m},x) \colon (\boldsymbol{m},x) \in \overline{M}_{\delta_0} \times \overline{B}_{\delta_0} \}.
\]
Since \(\overline{M}_{\delta_0}\) and \(\overline{B}_{\delta_0}\) are compact and connected, and since \(\gamma_i\) and \(t\) are continuous, the sets \(U_{i,\delta_0}\) and \(T_{\delta_0}\) are compact intervals. For \(\delta>0\), define
\[
\overline{\mathcal{U}}_{\delta,i}^{\delta_0} \coloneqq [\inf U_{i,\delta_0}  - \delta, \sup U_{i,\delta_0}  + \delta], 
\qquad
\overline{\mathcal{T}}_{\delta}^{\delta_0} \coloneqq [\inf T_{\delta_0} - \delta, \sup T_{\delta_0}  + \delta],
\]
and set \(\overline{\mathcal{U}}_{\delta}^{\delta_0} \coloneqq \prod_{i=1}^r \overline{\mathcal{U}}_{\delta,i}^{\delta_0}\). Arguing as in the upper bound, for every \(\delta>0\), there exists \(N_{\delta,\delta_0}>0\) such that for all \(N \geq N_{\delta,\delta_0}\),
\[
\sqrt{\frac{N}{N-1}}\, \gamma_i (\boldsymbol{m}) \in \overline{\mathcal{U}}_{\delta,i}^{\delta_0}
\quad \text{and} \quad
\sqrt{\frac{N}{N-1}}\,t(\boldsymbol{m},x) \in \overline{\mathcal{T}}_\delta^{\delta_0}
\]
uniformly over \((\boldsymbol{m},x) \in \overline{M}_{\delta_0} \times \overline{B}_{\delta_0}\). By Lemma~\ref{lem: det J_N} (in particular~\eqref{eq: lower bound det J_N}),
\[
\liminf_{N \to \infty} \inf_{(\boldsymbol{\gamma},t) \in \, \overline{\mathcal{U}}_\delta^{\delta_0} \times \overline{\mathcal{T}}_\delta^{\delta_0}} \frac{1}{N-1} \log \E \left [|\det(\boldsymbol{X}_{N-1}^{\boldsymbol{\gamma}}-t\boldsymbol{I}_{N-1})| \right] \geq \inf_{t \in \overline{\mathcal{T}}_\delta^{\delta_0}} \Omega(t).
\]
Hence for every \(\varepsilon>0\), there exists \(N_{\varepsilon,\delta,\delta_0}\ge N_{\delta,\delta_0}\) such that for all \(N \geq N_{\varepsilon,\delta,\delta_0}\),
\begin{equation} \label{lower bound 2}
\begin{split} 
\inf_{(\boldsymbol{m},x) \in M_{\delta_0} \times B_{\delta_0}} \E \left [|\det \boldsymbol{H}_{N-1}(\boldsymbol{m},x)| \right ]
& \geq  \inf_{(\boldsymbol{\gamma}, t) \in \, \overline{\mathcal{U}}_\delta^{\delta_0} \times \overline{\mathcal{T}}_\delta^{\delta_0}} \E \left [ |\det(\boldsymbol{X}_{N-1}^{\boldsymbol{\gamma}}-t\boldsymbol{I}_{N-1})| \right ]\\
& \geq \exp \left \{ (N-1) \left ( \inf_{t \in \overline{\mathcal{T}}^{\delta_0}_\delta} \Omega(t) - \varepsilon \right ) \right \}.
\end{split}
\end{equation}
Combining~\eqref{lower bound 1} and~\eqref{lower bound 2}, we obtain
\[
\begin{split}
\liminf_{N \to \infty} \frac{1}{N} \log \E \left [ \textnormal{Crt}_N^{\textnormal{tot}}(M,B) \right ]& \geq \inf_{(\boldsymbol{m},x) \in M_{\delta_0} \times B_{\delta_0}} S(\boldsymbol{m},x) + \inf_{t \in \overline{\mathcal{T}}^{\delta_0}_\delta} \Omega(t)  -\varepsilon,
\end{split}
\]
where we used the fact that \((2\delta_0)^{r+1}C_N\) is exponentially trivial. By compactness of the interval \(T_{\delta_0}\) and continuity of \(\Omega\), we obtain
\[
\lim_{\delta\to0_+}\inf_{t\in \overline{\mathcal{T}}_\delta^{\delta_0}}\Omega(t) = \inf_{t\in T_{\delta_0}}\Omega(t) = \inf_{(\boldsymbol{m},x)\in \overline{M}_{\delta_0}\times \overline{B}_{\delta_0}}\Omega(t(\boldsymbol{m},x)).
\]
Therefore, letting first \(\varepsilon\to 0_+\) and then \(\delta\to 0_+\) yields
\[
\liminf_{N \to \infty} \frac{1}{N} \log \E \left [ \textnormal{Crt}_N^{\textnormal{tot}}(M,B) \right ] \geq  \inf_{(\boldsymbol{m},x) \in M_{\delta_0} \times B_{\delta_0}}  S(\boldsymbol{m},x) +  \inf_{(\boldsymbol{m},x) \in \overline{M}_{\delta_0} \times \overline{B}_{\delta_0}} \Omega(t(\boldsymbol{m},x)).
\]
Since \(S\) is continuous on \(\mathcal{D}_r\times\R\), \(\Omega\) is continuous on \(\R\), and \(t\) is continuous on \([-1,1]^r \times \R\) (see Lemma~\ref{lem: GRF} and~\eqref{function t}), letting \(\delta_0 \to 0_+\) yields
\[
\begin{split}
\liminf_{N \to \infty} \frac{1}{N} \log \E \left [ \textnormal{Crt}_N^{\textnormal{tot}} (M,B) \right ]& \geq S(\boldsymbol{m}_0,x_0) + \Omega(t(\boldsymbol{m}_0,x_0)) \\
& = \Sigma^{\textnormal{tot}}(\boldsymbol{m}_0,x_0) \\
& \geq \sup_{(\boldsymbol{m}, x) \in M^\circ \times B^\circ} \Sigma^{\textnormal{tot}}(\boldsymbol{m}, x)-\varepsilon_0,
\end{split}
\]
where the last inequality follows by~\eqref{lower bound}. Finally, letting \(\varepsilon_0 \to 0_+\) proves the lower bound~\eqref{inf crit points}. 
\end{proof}

It remains to prove the intermediate results, namely Lemmas~\ref{lem: GRF},~\ref{lem: exp tight}, and~\ref{lem: det J_N}.

\begin{proof}[Proof of Lemma~\ref{lem: GRF}]
The function \(S\), defined in~\eqref{function S}, is continuous on \(\mathcal{D}_r \times \R\) because it is a sum of continuous functions. Moreover, \(\Omega\) is continuous in \(\R\) by~\eqref{def: function Phi star}, and \(t\) is continuous by~\eqref{function t}. Therefore, the function \((\boldsymbol{m},x)\mapsto S(\boldsymbol{m},x) +\Omega(t(\boldsymbol{m},x))\) is continuous on \(\mathcal{D}_r\times\R\).

We first prove that \(- \Sigma^{\textnormal{tot}}\) is lower semicontinuous, or equivalently that \(\Sigma^{\textnormal{tot}}\) is upper semicontinuous. By Definition~\ref{def: function Sigma tot}, 
\[
\Sigma^\textnormal{tot}(\boldsymbol{m},x) = S(\boldsymbol{m},x) + \Omega(t(\boldsymbol{m},x)), 
\]
for \((\boldsymbol{m},x) \in \mathcal{D}_r \times \R\), whereas
\[
\Sigma^\textnormal{tot}(\boldsymbol{m},x)  = - \infty,  
\]
otherwise. Since \(\Sigma^{\textnormal{tot}}\) is continuous on \(\mathcal{D}_r\times\R\), it remains to verify upper semicontinuity at points of \(\mathcal{D}_r^{\textnormal{c}}\times\R\). Let \((\boldsymbol{m}_n)_{n\geq1}\subset \mathcal{D}_r\) and \((x_n)_{n\geq1}\subset\R\) be such that
\[ 
\boldsymbol{m}_n\to \boldsymbol{m}, \qquad x_n\to x, 
\] 
where \(\boldsymbol{m}\in \mathcal{D}_r^{\textnormal{c}}\) and \(x\in\R\). Since the function $q(\boldsymbol{m})=1-\|\boldsymbol{m}\|_{2}^{2}$ is continuous, we have that \(q(\boldsymbol{m}_n)\to q(\boldsymbol{m})\). Moreover, \(q(\boldsymbol{m}_n)>0\) for every \(n\), while \(q(\boldsymbol{m})\leq0\) since \(\boldsymbol{m}\notin \mathcal{D}_r\). Therefore, \(q(\boldsymbol{m})=0\) and \(q(\boldsymbol{m}_n)\to 0\). By the explicit formula for \(S_0\) given in Remark~\ref{rmk: KR S0}, \(S_0\) is bounded above on \([-1,1]^r\times K\) for every compact set \(K\subset\R\). Since \(x_n\to x\), the sequence \((x_n)_{n\geq1}\) is contained in a compact subset of \(\R\), and hence \(\sup_{n\geq1} S_0(\boldsymbol{m}_n,x_n)<\infty\). Therefore, 
\[ 
S(\boldsymbol{m}_n,x_n) =S_0 (\boldsymbol{m}_n,x_n) + \frac{1}{2} \log q(\boldsymbol{m}_n) \to -\infty . 
\] 
On the other hand, by the explicit formula~\eqref{function t}, the sequence \(t(\boldsymbol{m}_n,x_n)\) is bounded, since \(x_n\to x\) and \(\boldsymbol{m}_n\in[-1,1]^r\). Since \(\Omega\) is continuous on \(\R\), it follows that the sequence \(\Omega(t(\boldsymbol{m}_n,x_n))\) is bounded above. Therefore 
\[ 
\limsup_{n\to\infty} \Sigma^{\textnormal{tot}}(\boldsymbol{m}_n,x_n) = -\infty = \Sigma^{\textnormal{tot}}(\boldsymbol{m},x). 
\] 
Thus \(\Sigma^{\textnormal{tot}}\) is upper semicontinuous on \([-1,1]^r\times\R\), and consequently \(-\Sigma^{\textnormal{tot}}\) is lower semicontinuous.

It remains to prove that the sub-level sets of \(-\Sigma^{\textnormal{tot}}\) are compact. Equivalently, we show that, for every \(a < \infty\), the set 
\[
\left \{ (\boldsymbol{m},x) \in [-1,1]^r \times \R \colon \Sigma^{\textnormal{tot}}  (\boldsymbol{m},x) \geq - a \right\} 
\]
is compact. By the lower semicontinuity of \(- \Sigma^{\textnormal{tot}}\), this set is closed. It therefore suffices to show that the sub-level sets are included in a compact subset of \([-1,1]^r \times \R\). According to~\eqref{function S}, the function \(S(\boldsymbol{m},x)\) can be bounded above by
\begin{equation} \label{bound function S}
S(\boldsymbol{m},x) \leq \frac{1}{2} \log(p-1) + \frac{1}{2} - x^2 + 2\left (\sum_{i=1}^r  \lambda_i m_i^{k_i} \right ) x .
\end{equation}
Here, we used the inequalities \(\log(1-\|\boldsymbol{m}\|_2^2) \leq 0\), \(\left ( \sum_{i=1}^r \lambda_i m_i^{k_i} \right)^2 \ge 0\), and 
\[
\left ( \sum_{i=1}^r \lambda_i k_i m_i^{k_i}\right)^2 \leq \left ( \sum_{i=1}^r \lambda_i^2 k_i^2 m_i^{2k_i-2} \right) \left( \sum_{i=1}^r  m_i^2\right) \le  \sum_{i=1}^r \lambda_i^2 k_i^2 m_i^{2k_i-2},
\]
by the Cauchy--Schwarz inequality. Additionally, since \(\Omega(x) \leq \frac{x^2}{4}-\frac{1}{2}\) for all \(x\in \R\) by~\eqref{def: function Phi star}, we have 
\begin{equation} \label{bound function Phi star}
\begin{split}
\Omega(t(\boldsymbol{m},x)) &\leq \frac{p}{2(p-1)} \left ( x - \sum_{i=1}^r \lambda_i \left (1- \frac{k_i}{p} \right )m_i^{k_i} \right )^2 - \frac{1}{2}\\
& \le \frac{p}{2(p-1)} x^2 -  \frac{p}{p-1} \left( \sum_{i=1}^r \lambda_i \left (1- \frac{k_i}{p} \right ) m_i^{k_i}  \right) x + C - \frac{1}{2},
\end{split}
\end{equation}
where \(C = C (p,r, (k_i)_i, (\lambda_i)_i) = r^2 \max_{1 \le i \le r} \{\lambda_i^2 (1 - k_i/p)^2 \} < \infty\). Combining~\eqref{bound function S} and~\eqref{bound function Phi star}, we obtain
\begin{equation} \label{bound sigma tot}
\begin{split}
\Sigma^\textnormal{tot}(\boldsymbol{m},x) &= S(\boldsymbol{m},x) + \Omega(t(\boldsymbol{m},x)) \\
& \le - \frac{p-2}{2(p-1)}x^2 + \left ( \frac{p-2}{p-1} \sum_{i=1}^r \lambda_i m_i^{k_i} + \frac{1}{p-1} \sum_{i=1}^r \lambda_i k_i m_i^{k_i}\right)x + \frac{1}{2} \log(p-1) + C\\
& \le - \frac{p-2}{2(p-1)}x^2  +  C_1 |x| + C_2 ,
\end{split}
\end{equation}
where 
\[ 
C_1 = \frac{p-2}{p-1}\sum_{i=1}^r |\lambda_i| + \frac{1}{p-1}\sum_{i=1}^r |\lambda_i|k_i <\infty, 
\qquad C_2=C+\frac12\log(p-1). 
\] 
Indeed, since \(|m_i|\le1\), we have 
\[ 
\left| \frac{p-2}{p-1}\sum_{i=1}^r \lambda_i m_i^{k_i} + \frac{1}{p-1}\sum_{i=1}^r \lambda_i k_i m_i^{k_i} \right| \le C_1. 
\]
Since \(p>2\), the right-hand side of~\eqref{bound sigma tot} tends to \(-\infty\) as \(|x|\to\infty\). Therefore, for every \(a<\infty\), there exists \(R_a<\infty\) such that
\[ 
\left\{ \Sigma^{\textnormal{tot}}\ge -a \right\} \subset [-1,1]^r\times[-R_a,R_a]. 
\]
Since the set \(\{\Sigma^{\textnormal{tot}}\ge -a\}\) is closed, it is a closed subset of the compact set \([-1,1]^r\times[-R_a,R_a]\). Hence it is compact. This proves that \(-\Sigma^{\textnormal{tot}}\) is a good rate function.
\end{proof}

We next show exponential tightness.

\begin{proof}[Proof of Lemma~\ref{lem: exp tight}]
If \(\boldsymbol{m} \notin \mathcal{D}_r\), then \(S(\boldsymbol{m},x)=-\infty\). Hence, by Proposition~\ref{formula complexity}, 
\[
\E \left [ \textnormal{Crt}_N^{\textnormal{tot}} \left (\mathcal{D}_r^\text{c}, (-\infty,-z]\cup[z,\infty) \right) \right ]= 0. 
\]
It remains to control the contribution from \(\mathcal{D}_r\). By Proposition~\ref{formula complexity} and Remark~\ref{rmk: KR S0}, 
\[ 
\begin{split} 
&\E \left[ \textnormal{Crt}_N^{\textnormal{tot}} \left( \mathcal{D}_r, (-\infty,-z]\cup[z,\infty) \right) \right] \\ 
& \le C_N \int_{\{|x|\ge z\}} \int_{\mathcal{D}_r} \E\left[ |\det \boldsymbol H_{N-1}(\boldsymbol{m},x)| \right] \exp\left\{ NS_0(\boldsymbol{m},x) + \frac{N-r-2}{2}\log q(\boldsymbol{m}) \right\} \,\textnormal{d}\boldsymbol{m}\,\textnormal{d}x . 
\end{split} 
\] 
Since \(0<q(\boldsymbol{m})\le1\), we have \(\log q(\boldsymbol{m})\le0\). Therefore, for all \(N\ge r+3\), \(\frac{N-r-2}{2}\log q(\boldsymbol{m})\le0\). Consequently, 
\[ 
\begin{split} 
&\E\left[ \textnormal{Crt}_N^{\textnormal{tot}} \left( \mathcal{D}_r, (-\infty,-z]\cup[z,\infty) \right) \right] \le C_N \int_{\{|x|\ge z\}} \int_{\mathcal{D}_r} \E\left[ |\det \boldsymbol H_{N-1}(\boldsymbol{m},x)| \right] \exp\left\{ NS_0(\boldsymbol{m},x) \right\} \,\textnormal{d}\boldsymbol{m}\,\textnormal{d}x . 
\end{split} 
\]
By the explicit formula for \(S\), the same estimates as in the proof of Lemma~\ref{lem: GRF}, now applied to \(S_0\), give constants \(C_0,C_1<\infty\) such that, uniformly in \(\boldsymbol{m}\in \mathcal{D}_r\), 
\[ 
S_0(\boldsymbol{m},x) \le C_0-x^2+C_1|x|. 
\] 
Thus, for \(z\ge 2C_1\) and \(|x|\ge z\), it follows that \(C_1|x|\le \frac{x^2}{2}\) and so 
\[
S_0(\boldsymbol{m},x) \le C_0-\frac{x^2}{2}. 
\]
We next bound the determinant term. As shown in the proof of the upper bound of Theorem~\ref{thm crit points}, the functions \(\gamma_i\) are bounded on \(\mathcal{D}_r\). Thus, there exists \(K< \infty\) such that \(|\gamma_i(\boldsymbol{m})|\le K\) for every \(\boldsymbol{m}\in \mathcal{D}_r\) and \(1\le i\le r\). Moreover, by the explicit expression for \(t\) in~\eqref{function t}, there exists \(c<\infty\) such that \(|t(\boldsymbol{m},x)|\le c (1+|x|)\) for every \(\boldsymbol{m}\in \mathcal{D}_r\) and \(x\in\R\). Consequently, using~\eqref{Hess} and the inequality \( |\det \boldsymbol{A}| \leq \|\boldsymbol{A}\|^{N-1}\) for every \((N-1) \times (N-1)\) matrix \(\boldsymbol{A}\), 
\[ 
|\det\boldsymbol{H}_{N-1}(\boldsymbol{m},x)| \leq \left( \|\boldsymbol{W}_{N-1}\|_{\mathrm{op}} + \sqrt{\frac{N}{N-1}} K + \sqrt{\frac{N}{N-1}} c(1+|x|) \right)^{N-1}.
\] 
Since \(\sqrt{\frac{N}{N-1}} \) is uniformly bounded for \(N\geq2\), there exists \(C<\infty\), independent of \(N\), \(\boldsymbol{m}\), and \(x\), such that 
\[ 
|\det \boldsymbol H_{N-1}(\boldsymbol{m},x)| \le \left( \|\boldsymbol W_{N-1}\|+C(1+|x|) \right)^{N-1}.
\] 
In particular, for \(|x|\ge z\ge  \max\{1,2C_1\}\), 
\[ 
|\det \boldsymbol H_{N-1}(\boldsymbol{m},x)| \le \left( \|\boldsymbol W_{N-1}\|+2C|x| \right)^{N-1}. 
\] 
Using \(\operatorname{Vol}(\mathcal{D}_r)\le 2^r\), the symmetry of the integration domain, and the change of variables \(u=|x|\), we obtain 
\[ 
\begin{split} 
&\E \left[ \textnormal{Crt}_N^{\textnormal{tot}} \left( \mathcal{D}_r, (-\infty,-z]\cup[z,\infty) \right) \right] \le K_N \int_z^\infty \E\left[ \left( \|\boldsymbol W_{N-1}\|+2Cu \right)^{N-1} \right] \exp\left\{ -\frac{Nu^2}{2} \right\} \,\textnormal{d}u, 
\end{split} 
\] 
where \(K_N\coloneqq 2^{r+1}C_N e^{NC_0}\). Since \(C_N\) is exponentially trivial by Proposition~\ref{formula complexity}, the sequence \((K_N)_{N\geq1}\) is exponentially finite, i.e.,
\[
\limsup_{N \to \infty} \frac{1}{N} \log K_N  <\infty. 
\]
For \(u\ge1\), there exists \(C'<\infty\) such that 
\[ 
\left( \|\boldsymbol W_{N-1}\|+2Cu \right)^{N-1} \le (C')^N (1+\|\boldsymbol W_{N-1}\|)^{N-1} u^{N-1}. 
\] 
Therefore, 
\[ 
\begin{split} 
&\E\left[ \textnormal{Crt}_N^{\textnormal{tot}} \left( \mathcal{D}_r, (-\infty,-z]\cup[z,\infty) \right) \right] \le K_N(C')^N \E\left[ (1+\|\boldsymbol W_{N-1}\|)^{N-1} \right] \int_z^\infty u^{N-1} \exp\left\{ -\frac{Nu^2}{2} \right\} \,\textnormal{d}u . 
\end{split} 
\] 
Since 
\[ 
\left( 1+\|\boldsymbol{W}_{N-1}\| \right)^{N-1} \leq \exp\left\{ N\|\boldsymbol{W}_{N-1}\| \right\}, 
\] 
the sub-Gaussian tail estimate for the operator norm of a GOE matrix, see for example~\cite[Lemma~6.3]{MR1856194}, implies that there exists \(C''<\infty\) such that  
\[ 
\E\left[ (1+\|\boldsymbol W_{N-1}\|)^{N-1} \right] \le  \E \left [ e^{N \lVert \boldsymbol{W}_{N-1} \rVert}\right ] \le (C'')^N. 
\] 
Hence,
\[ 
\begin{split} 
&\lim_{z\to\infty} \limsup_{N\to\infty} \frac1N \log \E\left[ \textnormal{Crt}_N^{\textnormal{tot}} \left( \mathcal{D}_r, (-\infty,-z]\cup[z,\infty) \right) \right] \\ 
&\le \limsup_{N\to\infty}\frac1N\log K_N +\log(C'C'') + \lim_{z\to\infty} \limsup_{N\to\infty} \frac1N \log \int_z^\infty u^{N-1} \exp\left\{ -\frac{Nu^2}{2} \right\} \textnormal{d}u .
\end{split} 
\] 
By~\cite[Lemma A.3]{MR4011861}, 
\[ 
\lim_{z\to\infty} \limsup_{N\to\infty} \frac1N \log \int_z^\infty u^{N-1} \exp\left\{ -\frac{Nu^2}{2} \right\} \,\textnormal{d}u = -\infty. 
\] 
Therefore, 
\[ 
\lim_{z\to\infty} \limsup_{N\to\infty} \frac1N \log \E\left[ \textnormal{Crt}_N^{\textnormal{tot}} \left( \mathcal{D}_r, (-\infty,-z]\cup[z,\infty) \right) \right] = -\infty. 
\] 
Together with the vanishing contribution from \(\mathcal{D}_r^{\textnormal{c}}\), this proves the exponential tightness estimate.
\end{proof}

It remains to prove the determinant asymptotics in Lemma~\ref{lem: det J_N}. The main approach is to apply Theorem 1.2 from~\cite{MR4552227}. 

\begin{proof}[Proof of Lemma~\ref{lem: det J_N}]
Set \(\boldsymbol{H}_N (\boldsymbol{\gamma},t)=\boldsymbol{X}_N^{\boldsymbol{\gamma}} - t\boldsymbol{I}_N\), where \(\boldsymbol{X}_N^{\boldsymbol{\gamma}} = \boldsymbol{W}_N + \sum_{i=1}^r \gamma_i \boldsymbol{e}_i\boldsymbol{e}_i^\top\) is a spiked GOE matrix. Consider the measure \(\mu_{\textnormal{sc}} (t)\), parametrized by \(t \in \R\), which is the shift of the semicircle law \(\mu_{\textnormal{sc}}\), whose density with respect to the Lebesgue measure is given by
\[
\mu_{\textnormal{sc}}(t) (\textnormal{d}  x)= \frac{1}{2\pi} \sqrt{4 - (x + t)^2} \mathbf{1}_{\{|x+t|\le 2\}} \textnormal{d}x.
\]
First, we make the following claims:
\begin{enumerate}
\item For all compact sets \(\mathcal{U} \subset \R^r\) and \(\mathcal{T} \subset \R\), there exists \(\kappa >0\) such that
\[
\sup_{(\boldsymbol{\gamma},t) \in \, \mathcal{U} \times \mathcal{T}} W_1(\E \left [\hat{\mu}_{\boldsymbol{H}_N (\boldsymbol{\gamma},t)} \right ], \mu_{\textnormal{sc}}(t)) \leq N^{-\kappa}.
\]

\item There exist universal constants \(C,c >0\) such that for every Lipschitz function \(f \colon \R \to \R\), every \(\delta >0\), and every compact \(\mathcal{U} \subset \R^r\) and \(\mathcal{T} \subset \R\),
\[ 
\sup_{(\boldsymbol{\gamma},t) \in \, \mathcal{U} \times \mathcal{T}} \mathbb{P} \left ( \left | \frac{1}{N} \Tr f(\boldsymbol{H}_N(\boldsymbol{\gamma},t)) - \frac{1}{N} \E \left [ \Tr f(\boldsymbol{H}_N(\boldsymbol{\gamma},t)) \right ] \right |> \delta \right )  \leq C \exp \left ( - \frac{c N^2 \delta^2}{\lVert f \rVert_{\textnormal{Lip}}^2}  \right ).
\]

\item For every \(\varepsilon > 0\) and every compact \(\mathcal{U} \subset \R^r\) and \(\mathcal{T} \subset \R\),
\[
\lim_{N \to \infty}  \sup_{(\boldsymbol{\gamma},t) \in \, \mathcal{U} \times \mathcal{T}}  \mathbb{P} \left (\boldsymbol{H}_N(\boldsymbol{\gamma},t) \: \text{has no eigenvalues in} \: [-e^{-N^\epsilon}, e^{-N^\epsilon}] \right ) = 1.
\]
\end{enumerate}
Given claims \(1\) through \(3\), the assumptions of~\cite[Theorem 1.2]{MR4552227} are satisfied uniformly over compact sets of \((\boldsymbol{\gamma},t)\). This implies that the conclusion of~\cite[Theorem 1.2]{MR4552227} is also uniform over compact sets (see also the proof of~\cite[Theorem 4.1]{MR4552227}). Specifically, for all compact sets \(\mathcal{U} \subset \R^r\) and \(\mathcal{T} \subset \R\),
\[
\begin{split}
\limsup_{N \to \infty} \frac{1}{N} \log \sup_{(\boldsymbol{\gamma},t) \in \, \mathcal{U} \times \mathcal{T}} \left \{  \E \left [ |\det \boldsymbol{H}_N(\boldsymbol{\gamma},t) | \right ]  e^{- N \int_\R \log |\lambda|  \mu_{\textnormal{sc}}(t) (\textnormal{d}\lambda) } \right \} & \le 0,\\
\liminf_{N \to \infty} \frac{1}{N} \log \inf_{(\boldsymbol{\gamma},t) \in \, \mathcal{U} \times \mathcal{T}} \left \{  \E \left [ |\det \boldsymbol{H}_N(\boldsymbol{\gamma},t) | \right ]  e^{- N \int_\R \log |\lambda|  \mu_{\textnormal{sc}}(t) (\textnormal{d}\lambda) } \right \} & \ge 0.
\end{split}
\]
Note that \(\int_\R \log |\lambda|  \mu_{\textnormal{sc}}(t) (\textnormal{d}\lambda) = \Omega (t)\). This completes the proof of Lemma~\ref{lem: det J_N}. 

It remains to prove the random-matrix claims. We start with the assumption on the Wasserstein distance. Since \(W_1\) is invariant under translations, \(W_1(\E \left [ \hat{\mu}_{\boldsymbol{H}_N (\boldsymbol{\gamma},t)} \right ],\mu_{\textnormal{sc}}(t)) \) is independent of \(t\). It therefore suffices to consider \(t=0\). Recall that \(\boldsymbol{H}_N (\boldsymbol{\gamma},0) =\boldsymbol{X}_N^{\boldsymbol{\gamma}}\) and that \(\mu_{\textnormal{sc}}(0) = \mu_{\textnormal{sc}}\). By the triangle inequality, we have 
\[
\sup_{\boldsymbol{\gamma} \in \, \mathcal{U}} W_1(\E [ \hat{\mu}_{\boldsymbol{X}_N^{\boldsymbol{\gamma}}}  ], \mu_{\textnormal{sc}}) \leq \sup_{\boldsymbol{\gamma} \in \,\mathcal{U}} W_1(\E  [ \hat{\mu}_{\boldsymbol{X}_N^{\boldsymbol{\gamma}}} ], \E \left [\hat{\mu}_{\boldsymbol{W}_N} \right ]) + W_1(\E \left [\hat{\mu}_{\boldsymbol{W}_N} \right ], \mu_{\textnormal{sc}}).
\]
By~\cite{MR2967966}, 
\[ 
\E\left[W_1(\hat\mu_{\boldsymbol{W}_N},\mu_{\mathrm{sc}})\right] \le C\frac{\sqrt{\log N}}{N}. 
\] 
Hence 
\begin{equation} \label{Dallaporta}
\begin{split}
W_1(\E \left [\hat{\mu}_{\boldsymbol{W}_N}\right ], \mu_{\textnormal{sc}}) & = \sup_{\lVert f \rVert_{\textnormal{Lip}}\leq 1} \left | \E \left[ \int_\R f(x) (\hat{\mu}_{\boldsymbol{W}_N} - \mu_{\textnormal{sc}})(\textnormal{d}x) \right ] \right | \\
& \leq  \E \left [W_1(\hat{\mu}_{\boldsymbol{W}_N}, \mu_{\textnormal{sc}}) \right ]\\
& \leq C \frac{\sqrt{\log N}}{N}.
\end{split}
\end{equation}
We now bound the finite-rank perturbation term 
\[
W_1\left(\E  [\hat{\mu}_{\boldsymbol{X}_N^{\boldsymbol{\gamma}}} ], \E \left [\hat{\mu}_{\boldsymbol{W}_N} \right ] \right )  = \sup_{\lVert f \rVert_{\textnormal{Lip}}\leq 1} \E \left [  \left | \frac{1}{N} \Tr f(\boldsymbol{X}_N^{\boldsymbol{\gamma}}) - \frac{1}{N}\Tr f(\boldsymbol{W}_N) \right | \right ].
\]
For every Lipschitz function \(f \colon \R \to \R\),
\begin{equation} \label{bound trace}
\begin{split}
\left | \frac{1}{N} \Tr f(\boldsymbol{X}_N^{\boldsymbol{\gamma}}) - \frac{1}{N}\Tr f(\boldsymbol{W}_N) \right | & \leq \frac{\lVert f \rVert_{\textnormal{Lip}} }{N} \sum_{i=1}^N \left | \lambda_i(\boldsymbol{X}_N^{\boldsymbol{\gamma}}) - \lambda_i(\boldsymbol{W}_N) \right |\\
& \leq \frac{\lVert f \rVert_{\textnormal{Lip}} }{\sqrt{N}}  \left ( \sum_{i=1}^N \left | \lambda_i(\boldsymbol{X}_N^{\boldsymbol{\gamma}}) - \lambda_i(\boldsymbol{W}_N) \right |^2 \right )^{\frac{1}{2}}\\
& \leq \frac{\lVert f \rVert_{\textnormal{Lip}} }{\sqrt{N}}  \left ( \sum_{i=1}^r \gamma_i^2 \right )^{1/2},
\end{split}
\end{equation}
where the first inequality follows by the fact that \(f\) is Lipschitz, the second by the Cauchy--Schwarz inequality, and the last one by the Hoffman--Wielandt inequality (see e.g.~\cite[Lemma 2.1.19]{MR2760897}). Therefore, 
\[ 
\sup_{\boldsymbol{\gamma} \in\mathcal{U}} W_1\left( \E[\hat\mu_{\boldsymbol{X}_N^{\boldsymbol{\gamma}}}], \E[\hat{\mu}_{\boldsymbol{W}_N}] \right) \le \frac{1}{\sqrt{N}} \sup_{\boldsymbol{\gamma} \in\mathcal{U}} \| \boldsymbol{\gamma}\|_2 . 
\] 
Since \(\mathcal{U}\) is compact, the right-hand side is \(\mathcal{O}(N^{-1/2})\). Combining this with the previous estimate gives claim \(1\), for any \(\kappa<1/2\).

We next verify item 2 on concentration of Lipschitz traces. We note that the laws of the entries of \(\sqrt{N} \boldsymbol{H}_N(\boldsymbol{\gamma},t)\) satisfy the log-Sobolev inequality with a uniform constant, since they are Gaussians. This is true uniformly over \((\boldsymbol{\gamma},t)\) since \(\boldsymbol{\gamma}\) and \(t\) only affect the mean, and translating measures preserves log-Sobolev inequality with the same constant. Therefore, by~\cite[Theorem 2.3.5]{MR2760897}, there exist constants \(C,c>0\), independent of \((\boldsymbol{\gamma},t)\), such that for every Lipschitz function \(f \colon \R \to \R\) and every \(\delta>0\), 
\[ 
\sup_{(\boldsymbol{\gamma},t) \in \mathcal{U} \times \mathcal{T}} \mathbb{P} \left( \left| \frac{1}{N} \Tr f(\boldsymbol{H}_N(\boldsymbol{\gamma},t)) - \frac{1}{N} \E\left[\Tr f(\boldsymbol{H}_N(\boldsymbol{\gamma},t))\right] \right|>\delta \right) \le C\exp\left\{ -\frac{cN^2\delta^2}{\|f\|_{\textnormal{Lip}}^2} \right\}. 
\] 
This proves claim \(2\).

Finally, the small-eigenvalue condition follows uniformly over compact \(\mathcal{U} \times\mathcal{T}\) from~\cite[Theorem 2]{MR3631932}.
\end{proof}

\subsection{Proof of complexity result for local maxima} \label{subsection loc maxima}

This part is devoted to the proof of Theorem~\ref{thm local maxima} on the complexity of local maxima. We adopt a similar strategy as in the previous subsection for analyzing the complexity of total critical points. The following results play a key role in the proof of Theorem~\ref{thm local maxima}.

\begin{lem}[Good rate function] \label{lem: GRF max}
The function \(-\Sigma^{\textnormal{max}}\) given by Definition~\ref{def: function Sigma max} is a good rate function. Moreover, the function \(L\) given by Definition~\ref{def: function L} is lower semicontinuous.
\end{lem}

\begin{lem}[Exponential tightness] \label{lem: exp tight max}
It holds that
\[
\lim_{z \to \infty} \limsup_{N \to \infty} \frac{1}{N} \log \E \left [ \textnormal{Crt}_N^{\textnormal{max}}([-1,1]^r, (-\infty, -z] \cup [z, \infty)) \right ]= - \infty.
\]
\end{lem}

\begin{lem}[Determinant asymptotics] \label{lem: det J_N max}
The following statements hold.
\begin{enumerate}
\item Upper bound: For every compact product set \(\mathcal{U}= \mathcal{U}_1 \times \cdots \times \mathcal{U}_r \subset\R^r\) and every compact set \(\mathcal{T} \subset\R\), it holds that
\begin{equation} \label{UB max}
\begin{split}
& \limsup_{N \to \infty} \sup_{(\boldsymbol{\gamma},t) \in \, \mathcal{U} \times \mathcal{T}} \frac{1}{N} \log \E \left [|\det (\boldsymbol{X}_N^{\boldsymbol{\gamma}} -t \boldsymbol{I}_N)| \mathbf{1}_{\{\boldsymbol{X}_N^{\boldsymbol{\gamma}} -t \boldsymbol{I}_N \prec 0\}}\right] \\
& \leq \sup_{t \in \mathcal{T}} \Omega(t) - \inf_{(\boldsymbol{\gamma},t) \in \, \mathcal{U}   \times \mathcal{T}} L(\boldsymbol{\gamma}, t).
\end{split}
\end{equation}
\item Lower bound: For every \(\eta_0 >0, \boldsymbol{\gamma}_0 \in\R^r\), and \(t_0 > 2\), there exists \(\rho_0>0\) such that \(t_0 - \rho_0 > 2\), and 
\begin{equation}\label{LB max}
\begin{split}
& \liminf_{N \to \infty}  \inf_{(\boldsymbol{\gamma},t) \in \, \mathcal{U}_{\rho_0} \times \mathcal{T}_{\rho_0}} \frac{1}{N} \log \E \left [|\det(\boldsymbol{X}_N^{\boldsymbol{\gamma}} -t\boldsymbol{I}_N)| \mathbf{1}_{\{\boldsymbol{X}_N^{\boldsymbol{\gamma}} -t\boldsymbol{I}_N \prec 0\}}\right] \\
& \geq \Omega(t_0)-L(\boldsymbol{\gamma}_0,t_0) - \eta_0,
\end{split}
\end{equation}
where \(\mathcal{U}_{\rho_0} = \prod_{i=1}^r ((\gamma_0)_i-\rho_0, (\gamma_0)_i +\rho_0) \) and \(\mathcal{T}_{\rho_0} =(t_0 - \rho_0, t_0+\rho_0)\).
\end{enumerate}
\end{lem}

We defer the proof of these lemmas to the end of this subsection. We now prove Theorem~\ref{thm local maxima}. First we prove the upper bound~\eqref{sup loc maxima}. The proof follows similar arguments as the upper bound in Theorem~\ref{thm crit points}.

\begin{proof}[Proof of Theorem~\ref{thm local maxima} (Upper bound)]
If \(\sup_{(\boldsymbol{m},x) \in \overline{M} \times \overline{B}} \Sigma^{\textnormal{max}}(\boldsymbol{m},x) =-\infty\), then the desired upper bound is immediate. Henceforth we assume that \(\sup_{(\boldsymbol{m},x) \in \overline{M} \times \overline{B}} \Sigma^{\textnormal{max}}(\boldsymbol{m},x) > -\infty\). By exponential tightness, see Lemma~\ref{lem: exp tight max}, it suffices to prove the upper bound when \(B\) is bounded. Then \(\overline{M}\times\overline{B}\) is compact. As in the proof of the upper bound in Theorem~\ref{thm crit points}, it is enough to prove the following estimate. For \((\boldsymbol{m}_0,x_0)\in\overline{M}\times\overline{B}\) and \(\delta_0>0\), set
\[
M_{\delta_0} \coloneqq \prod_{i=1}^r \left  ((m_0)_i-\delta_0,(m_0)_i+\delta_0 \right ), \qquad
B_{\delta_0} \coloneqq (x_0-\delta_0,x_0+\delta_0).
\]
We claim that
\begin{equation} \label{eq:local-upper-max}
\limsup_{\delta_0\to0_+} \limsup_{N\to\infty} \frac{1}{N} \log \E\left[ \textnormal{Crt}_N^{\textnormal{max}}(M_{\delta_0},B_{\delta_0}) \right] \le \Sigma^{\textnormal{max}} (\boldsymbol{m}_0,x_0).
\end{equation}
Indeed, once~\eqref{eq:local-upper-max} is proved, the finite-cover argument used in the proof of Theorem~\ref{thm crit points} gives the desired global upper bound. 

We prove~\eqref{eq:local-upper-max}. We first note that if \( \Sigma^{\textnormal{tot}} (\boldsymbol{m}_0,x_0)=-\infty\), then the desired estimate follows from the already proved local upper bound for total critical points and from the inequality \(\textnormal{Crt}_N^{\textnormal{max}} (M_{\delta_0}, B_{\delta_0}) \le \textnormal{Crt}_N^{\textnormal{tot}}(M_{\delta_0}, B_{\delta_0})\). We may therefore assume that \(\Sigma^{\textnormal{tot}}(\boldsymbol{m}_0,x_0)>-\infty\). Then \((\boldsymbol{m}_0,x_0)\in \mathcal{D}_r  \times \R\). Since \(\mathcal{D}_r\) is open, we choose \(\delta_0>0\) sufficiently small such that \(\overline{M}_{\delta_0}\subset \mathcal{D}_r\). By Proposition~\ref{formula complexity}, in particular~\eqref{Kac--Rice index 2}, and Remark~\ref{rmk: KR S0},
\begin{equation} \label{eq:upper_max_1} 
\begin{split}
\E \left[ \textnormal{Crt}_N^{\textnormal{max}}(M_{\delta_0},B_{\delta_0}) \right]  & \le 2^{r+1}R_0 C_N \sup_{(\boldsymbol{m},x)\in \overline{M}_{\delta_0} \times \overline{B}_{\delta_0}} \E\left[ |\det\boldsymbol{H}_{N-1}(\boldsymbol{m},x)| \mathbf{1}_{\{\boldsymbol{H}_{N-1}(\boldsymbol{m},x)\prec0\}} \right]   \\
&\quad \times \sup_{(\boldsymbol{m},x)\in \overline{M}_{\delta_0}\times\overline{B}_{\delta_0}} \exp\left\{ NS_0 (\boldsymbol{m},x) +\frac{N-r-2}{2}\log q(\boldsymbol{m}) \right\},
\end{split}
\end{equation}
where \(R_0\coloneqq \sup\{|x| \colon x\in\overline{B}_{\delta_0}\}<\infty\). Here we used that the integration domain is contained in \([-1,1]^r\times[-R_0,R_0]\).

We now bound the determinant term. Define
\[
U_{i,\delta_0} \coloneqq \{\gamma_i(\boldsymbol{m}) \colon \boldsymbol{m}\in\overline{M}_{\delta_0}\}, \qquad
T_{\delta_0} \coloneqq \{t(\boldsymbol{m},x) \colon (\boldsymbol{m},x)\in \overline{M}_{\delta_0} \times \overline{B}_{\delta_0}\}.
\]
Since \(\overline{M}_{\delta_0}\) and \(\overline{M}_{\delta_0}\times\overline{B}_{\delta_0}\) are compact and connected, and since the functions \(\gamma_i\) and \(t\) are continuous on the respective domains, the sets \(U_{i,\delta_0}\) and \(T_{\delta_0}\) are compact intervals. For \(\delta >0\), define
\[
\overline{\mathcal{U}}^{\delta_0}_\delta \coloneqq \prod_{i=1}^r [\inf U_{i,\delta_0}-\delta,\sup U_{i,\delta_0}+\delta], \qquad
\overline{\mathcal{T}}^{\delta_0}_\delta \coloneqq [\inf T_{\delta_0}-\delta,\sup T_{\delta_0}+\delta].
\]
Since \(\sqrt{N/(N-1)}\to1\), there exists \(N_{\delta,\delta_0}\) such that, for all \(N\ge N_{\delta,\delta_0}\),
\[
\sqrt{\frac N{N-1}} \boldsymbol{\gamma}(\boldsymbol{m}) \in \overline{\mathcal{U}}^{\delta_0}_\delta,
\qquad
\sqrt{\frac N{N-1}}t(\boldsymbol{m},x) \in \overline{\mathcal{T}}^{\delta_0}_\delta,
\]
uniformly over \((\boldsymbol{m},x)\in \overline{M}_{\delta_0}\times\overline{B}_{\delta_0}\). Using the distributional identity~\eqref{Hess}, we obtain
\[
\begin{split}
& \sup_{(\boldsymbol{m},x)\in \overline{M}_{\delta_0}\times\overline{B}_{\delta_0}} \E \left[ |\det\boldsymbol{H}_{N-1}(\boldsymbol{m},x)| \mathbf{1}_{\{\boldsymbol{H}_{N-1}(\boldsymbol{m},x)\prec0\}} \right] \\
& \le \sup_{(\boldsymbol{\gamma},t)\in \overline{\mathcal{U}}^{\delta_0}_\delta \times \overline{\mathcal{T}}^{\delta_0}_\delta} \E\left[ \left| \det(\boldsymbol{X}_{N-1}^{\boldsymbol{\gamma}} -t \boldsymbol{I}_{N-1}) \right| \mathbf{1}_{\{\boldsymbol{X}_{N-1}^{\boldsymbol{\gamma}} -t\boldsymbol{I}_{N-1}\prec0\}} \right].
\end{split}
\]
By Lemma~\ref{lem: det J_N max}, in particular~\eqref{UB max}, for every \(\varepsilon>0\), there exists \(N_{\varepsilon, \delta,\delta_0} \ge  N_{\delta,\delta_0}\) such that for all \(N \ge N_{\varepsilon,\delta,\delta_0}\),
\begin{equation}  \label{eq:upper_max_2} 
\begin{split}
& \sup_{(\boldsymbol{m},x)\in \overline{M}_{\delta_0}\times\overline{B}_{\delta_0}} \E \left[ |\det\boldsymbol{H}_{N-1}(\boldsymbol{m},x)| \mathbf{1}_{\{\boldsymbol{H}_{N-1}(\boldsymbol{m},x)\prec0\}} \right]   \\
& \le \exp \left\{ (N-1) \left[ \sup_{t\in\overline{\mathcal{T}}^{\delta_0}_\delta} \Omega(t) - \inf_{(\boldsymbol{\gamma},t)\in \overline{\mathcal{U}}^{\delta_0}_\delta \times \overline{\mathcal{T}}^{\delta_0}_\delta} L(\boldsymbol{\gamma},t) +\varepsilon \right] \right\}.
\end{split}
\end{equation}
Combining~\eqref{eq:upper_max_1} and~\eqref{eq:upper_max_2}, taking \(\frac{1}{N}\log\), and using that \(C_N\) is exponentially trivial, we get
\[
\begin{split}
 \limsup_{N\to\infty} \frac{1}{N}\log \E\left[ \textnormal{Crt}_N^{\textnormal{max}}(M_{\delta_0},B_{\delta_0}) \right]  & \le \limsup_{N\to\infty} \sup_{(\boldsymbol{m},x)\in \overline{M}_{\delta_0}\times\overline{B}_{\delta_0}} \left[S_0(\boldsymbol{m},x) + \frac{N-r-2}{2N}\log q(\boldsymbol{m}) \right]   \\
 & \quad + \sup_{t\in\overline{\mathcal{T}}^{\delta_0}_\delta} \Omega(t) - \inf_{(\boldsymbol{\gamma},t)\in \overline{\mathcal{U}}^{\delta_0}_\delta \times \overline{\mathcal{T}}^{\delta_0}_\delta} L(\boldsymbol{\gamma},t) +\varepsilon .
\end{split}
\]
Since \(\overline{M}_{\delta_0}\subset \mathcal{D}_r\), the function \(q\) is bounded on \(\overline{M}_{\delta_0}\). Recalling that \(S(\boldsymbol{m},x) = S_0(\boldsymbol{m},x) + \frac{1}{2}\log q(\boldsymbol{m})\), we have
\[
\limsup_{N\to\infty} \sup_{(\boldsymbol{m},x)\in \overline{M}_{\delta_0}\times\overline{B}_{\delta_0}} \left[S_0(\boldsymbol{m},x) + \frac{N-r-2}{2N}\log q(\boldsymbol{m}) \right] \le \sup_{(\boldsymbol{m},x) \in \overline{M}_{\delta_0}\times\overline{B}_{\delta_0}} S(\boldsymbol{m},x) .
\]
Therefore,
\[
\limsup_{N\to\infty} \frac{1}{N}\log \E\left[ \textnormal{Crt}_N^{\textnormal{max}}(M_{\delta_0},B_{\delta_0}) \right] \le \sup_{\overline{M}_{\delta_0}\times\overline{B}_{\delta_0}} S + \sup_{t\in\overline{\mathcal{T}}^{\delta_0}_\delta}\Omega(t) - \inf_{(\boldsymbol{\gamma},t)\in \overline{\mathcal{U}}^{\delta_0}_\delta \times \overline{\mathcal{T}}^{\delta_0}_\delta} L(\boldsymbol{\gamma},t) +\varepsilon .
\]
Letting first \(\varepsilon\to 0^+\) and then \(\delta\to 0^+\), using the continuity of \(\Omega\) and the lower semicontinuity of \(L\), gives
\begin{equation} \label{eq:upper_bound_max_final}
\limsup_{N\to\infty} \frac{1}{N}\log \E\left[ \textnormal{Crt}_N^{\textnormal{max}}(M_{\delta_0},B_{\delta_0}) \right] \le \sup_{\overline{M}_{\delta_0}\times\overline{B}_{\delta_0}} S + \sup_{t\in T_{\delta_0}} \Omega(t) - \inf_{(\boldsymbol{\gamma},t)\in U_{\delta_0} \times T_{\delta_0}} L(\boldsymbol{\gamma},t),
\end{equation}
where \(U_{\delta_0} \coloneqq U_{1,\delta_0}\times\cdots\times U_{r,\delta_0}\). Indeed, the compact sets \(\overline{\mathcal{U}}^{\delta_0}_\delta \times \overline{\mathcal{T}}^{\delta_0}_\delta\) decrease to \(U_{\delta_0}\times T_{\delta_0}\) as \(\delta\to 0^+\). Hence, by the lower semicontinuity of \(L\),
\[
\lim_{\delta\to 0^+} \inf_{(\boldsymbol{\gamma},t)\in \overline{\mathcal{U}}^{\delta_0}_\delta \times \overline{\mathcal{T}}^{\delta_0}_\delta} L(\boldsymbol{\gamma},t) \ge \inf_{(\boldsymbol{\gamma},t)\in U_{\delta_0}\times T_{\delta_0}} L(\boldsymbol{\gamma},t).
\]
Finally we let \(\delta_0\to 0^+\). Since \(\overline{M}_{\delta_0}\times\overline{B}_{\delta_0}\) shrink to \(\{(\boldsymbol{m}_0,x_0)\}\) as \(\delta_0 \to 0^+\), and \(S,\Omega,t,\gamma_i\) are continuous on the relevant domain, we have
\[
\lim_{\delta_0\to 0^+} \sup_{\overline{M}_{\delta_0}\times\overline{B}_{\delta_0}} S = S(\boldsymbol{m}_0,x_0), \quad \textnormal{and} \quad \lim_{\delta_0\to 0^+} \sup_{t\in T_{\delta_0}}\Omega(t) = \Omega(t(\boldsymbol{m}_0,x_0)).
\]
Moreover, since \(L\) is lower semicontinuous and \(U_{\delta_0} \times T_{\delta_0}\) shrink to \(\{(\boldsymbol{\gamma}(\boldsymbol{m}_0),t(\boldsymbol{m}_0,x_0))\}\) as \(\delta_0 \to 0^+\), we have
\[
\lim_{\delta_0 \to 0^+} \inf_{(\boldsymbol{\gamma},t)\in U_{\delta_0}  \times T_{\delta_0}} L(\boldsymbol{\gamma},t) \ge L (\boldsymbol{\gamma} (\boldsymbol{m}_0), t(\boldsymbol{m}_0,x_0)).
\]
Taking the \(\limsup_{\delta_0\to0^+}\) in~\eqref{eq:upper_bound_max_final} therefore yields
\[
\begin{split}
\limsup_{\delta_0 \to 0^+} \limsup_{N\to\infty} \frac{1}{N}\log \E\left[ \textnormal{Crt}_N^{\textnormal{max}}(M_{\delta_0},B_{\delta_0}) \right] & \le S(\boldsymbol{m}_0,x_0) +
\Omega(t(\boldsymbol{m}_0,x_0)) - L(\boldsymbol{\gamma}(\boldsymbol{m}_0),t(\boldsymbol{m}_0,x_0))  \\
& = \Sigma^{\textnormal{max}}(\boldsymbol{m}_0,x_0).
\end{split}
\]
This proves~\eqref{eq:local-upper-max}, and therefore the upper bound for local maxima.
\end{proof}

We next show the lower bound~\eqref{inf loc maxima}. 

\begin{proof}[Proof of Theorem~\ref{thm local maxima} (Lower bound)] 
It suffices to consider Borel sets \(M\subset[-1,1]^r\) and \(B\subset\R\) such that \(\sup_{(\boldsymbol{m},x) \in M^\circ \times B^\circ} \Sigma^{\textnormal{max}}(\boldsymbol{m}, x) > -\infty\), since otherwise~\eqref{inf loc maxima} is trivial. Fix \(\varepsilon_0>0\). By definition of the supremum, there exists \((\boldsymbol{m}_0,x_0)\in M^\circ\times B^\circ\) such that
\begin{equation} \label{Sigma min ineq}
\Sigma^{\textnormal{max}} (\boldsymbol{m}_0,x_0) \ge \sup_{(\boldsymbol{m},x)\in M^\circ \times B^\circ} \Sigma^{\textnormal{max}} (\boldsymbol{m} , x) - \varepsilon_0.
\end{equation}
In particular, \(\Sigma^{\textnormal{max}} (\boldsymbol{m}_0,x_0)  > - \infty\), hence \((\boldsymbol{m}_0,x_0) \in \mathcal{D}_r \times \R\) and \(t_0 \coloneqq t(\boldsymbol{m}_0,x_0) > 2\). Set \(\boldsymbol{\gamma}_0 \coloneqq \boldsymbol{\gamma} (\boldsymbol{m}_0)\). 

Fix \(\eta_0>0\). By Lemma~\ref{lem: det J_N max}, choose \(\rho_0>0\) such that \(t_0-\rho_0>2\) and~\eqref{LB max} holds. For \(\delta_0>0\) define
\[
M_{\delta_0} \coloneqq \prod_{i=1}^r ((m_0)_i -\delta_0, (m_0)_i + \delta_0) , 
\qquad
B_{\delta_0}  \coloneqq (x_0-\delta_0, x_0+\delta_0) .
\]
Since \(M^\circ\cap \mathcal{D}_r\) and \(B^\circ\) are open, and since \((\boldsymbol{m}_0,x_0)\in (M^\circ\cap \mathcal{D}_r)\times B^\circ\), we may choose \(\delta_0>0\) sufficiently small such that \(\overline{M}_{\delta_0}\subset M^\circ\cap \mathcal{D}_r\) and \(\overline{B}_{\delta_0}\subset B^\circ\). In particular, \(S\) and \(t\) are continuous on \(\overline{M}_{\delta_0}\times\overline{B}_{\delta_0}\), while \(\boldsymbol{\gamma}\) is continuous on \(\overline{M}_{\delta_0}\). Next, since \((\boldsymbol{m},x)\mapsto (\boldsymbol{\gamma}(\boldsymbol{m}),t(\boldsymbol{m},x))\) is continuous at \((\boldsymbol{m}_0,x_0)\), there exists a neighborhood \(U\) of \((\boldsymbol{m}_0,x_0)\) such that
\[
(\boldsymbol{\gamma}(\boldsymbol{m}),t(\boldsymbol{m},x)) \in \prod_{i=1}^r \left ((\gamma_0)_i - \frac{\rho_0}{2}, (\gamma_0)_i + \frac{\rho_0}{2}\right ) \times \left (t_0- \frac{\rho_0}{2},\, t_0 + \frac{\rho_0}{2} \right )
\]
for all \((\boldsymbol{m},x)\in U\). Shrinking \(\delta_0>0\) if necessary so that \( \overline{M}_{\delta_0} \times \overline{B}_{\delta_0} \subset U\), we obtain
\[
\gamma_i (\boldsymbol{m}) \in \left ((\gamma_0)_i - \frac{\rho_0}{2},\, (\gamma_0)_i + \frac{\rho_0}{2}\right ), 
\qquad
t (\boldsymbol{m},x) \in \left (t_0- \frac{\rho_0}{2},\, t_0 + \frac{\rho_0}{2}\right ),
\]
uniformly for all \((\boldsymbol{m},x)\in  \overline{M}_{\delta_0}\times \overline{B}_{\delta_0}\). Since \(\sqrt{N/(N-1)}\to1\) and each \(\gamma_i\), as well as \(t\), is bounded on the compact set \(\overline{M}_{\delta_0}\times \overline{B}_{\delta_0}\), it follows that there exists \(N_{\rho_0, \delta_0}\) such that for all \(N \ge N_{\rho_0,\delta_0}\),
\[
\sqrt{\frac N{N-1}} \gamma_i (\boldsymbol{m})\in  \left ((\gamma_0)_i - \rho_0,\, (\gamma_0)_i +\rho_0 \right ),
\qquad
\sqrt{\frac N{N-1}} t(\boldsymbol{m},x)\in  \left (t_0- \rho_0,\, t_0 + \rho_0\right ),
\]
for all \((\boldsymbol{m},x)\in \overline{M}_{\delta_0}\times \overline{B}_{\delta_0}\) and all \(i=1,\dots,r\). Equivalently,
\[
\left ( \sqrt{\frac N{N-1}} \boldsymbol{\gamma} (\boldsymbol{m}) , \sqrt{\frac N{N-1}} t(\boldsymbol{m},x)\right) \in \mathcal{U}_{\rho_0} \times \mathcal{T}_{\rho_0},
\]
uniformly over \((\boldsymbol{m},x)\in \overline{M}_{\delta_0}\times \overline{B}_{\delta_0}\), where \(\mathcal{U}_{\rho_0} \coloneqq \prod_{i=1}^r ((\gamma_0)_i-\rho_0,(\gamma_0)_i+\rho_0)\) and \(\mathcal{T}_{\rho_0} \coloneqq (t_0-\rho_0,t_0+\rho_0)\).

By Proposition~\ref{formula complexity},
\[
\begin{split}
& \E \left [ \textnormal{Crt}_N^{\textnormal{max}}(M, B) \right ] \\
& \geq  \E  \left [  \textnormal{Crt}_N^{\textnormal{max}}(M_{\delta_0}, B_{\delta_0})  \right ]  \\
& = C_N \int_{M_{\delta_0}\times B_{\delta_0}} q(\boldsymbol{m})^{-\frac{r+2}{2}} \exp\left\{ NS(\boldsymbol{m},x) \right\} \E\left[ \left| \det\boldsymbol{H}_{N-1}(\boldsymbol{m},x) \right| \mathbf{1}_{\{ \boldsymbol{H}_{N-1}(\boldsymbol{m},x)\prec0 \}} \right] \,\textnormal{d}\boldsymbol{m}\textnormal{d}x.
\end{split}
\]
Since \(0<q(\boldsymbol{m})\leq1\) on \(\mathcal{D}_r\), we have \( q(\boldsymbol{m})^{-\frac{r+2}{2}}\geq1\). Therefore, 
\[ 
\begin{split} 
& \E\left[ \textnormal{Crt}_N^{\textnormal{max}}(M,B) \right] \\
&\geq C_N \inf_{(\boldsymbol{m},x)\in M_{\delta_0} \times B_{\delta_0}} \exp \left \{NS(\boldsymbol{m},x)\right \} \int_{M_{\delta_0}\times B_{\delta_0}} \E\left[ \left| \det\boldsymbol{H}_{N-1}(\boldsymbol{m},x) \right| \mathbf{1}_{\{ \boldsymbol{H}_{N-1}(\boldsymbol{m},x)\prec0 \}} \right] \,\textnormal{d}\boldsymbol{m} \textnormal{d}x. 
\end{split} 
\]
Therefore, for all \(N \ge N_{\rho_0,\delta_0}\),
\[
\begin{split}
&  \int_{M_{\delta_0} \times B_{\delta_0}} \E \left [|\det \boldsymbol{H}_{N-1}(\boldsymbol{m},x)| \mathbf{1}_{\{\boldsymbol{H}_{N-1}(\boldsymbol{m},x) \prec 0\}}\right] \textnormal{d} \boldsymbol{m}  \textnormal{d}x \\
&\ge (2\delta_0)^{r+1} \inf_{(\boldsymbol{m},x)\in M_{\delta_0} \times B_{\delta_0}} \E \left [|\det \boldsymbol{H}_{N-1}(\boldsymbol{m},x)| \mathbf{1}_{\{\boldsymbol{H}_{N-1}(\boldsymbol{m},x) \prec 0\}} \right] \\
&\ge (2\delta_0)^{r+1}\inf_{(\boldsymbol{\gamma},t) \in \mathcal{U}_{\rho_0}\times \mathcal{T}_{\rho_0}} \E \left [ |\det(\boldsymbol{X}_{N-1}^{\boldsymbol{\gamma}} - t \boldsymbol{I}_{N-1})| \mathbf{1}_{\{\boldsymbol{X}_{N-1}^{\boldsymbol{\gamma}} - t \boldsymbol{I}_{N-1}\prec 0\}} \right ].
\end{split}
\]
Since the factor \((2\delta_0)^{r+1}\) is exponentially trivial, it follows by Lemma~\ref{lem: det J_N max} (and in particular by~\eqref{LB max}) that
\[
\begin{split}
& \liminf_{N\to\infty} \frac{1}{N} \log \int_{M_{\delta_0} \times B_{\delta_0}} \E \left [|\det \boldsymbol{H}_{N-1}(\boldsymbol{m},x)| \mathbf{1}_{\{\boldsymbol{H}_{N-1}(\boldsymbol{m},x) \prec 0\}}\right]  \textnormal{d} \boldsymbol{m}  \textnormal{d}x  \\
& \ge \Omega(t_0) - L(\boldsymbol{\gamma}_0,t_0) -\eta_0.
\end{split}
\]
As a result,
\[
\liminf_{N\to\infty} \frac{1}{N} \log\E \left [ \textnormal{Crt}_N^{\textnormal{max}}(M, B) \right ]  \ge \inf_{(\boldsymbol{m},x)\in M_{\delta_0} \times B_{\delta_0}} S(\boldsymbol{m},x) +  \Omega(t_0) - L(\boldsymbol{\gamma}_0,t_0) -\eta_0,
\]
where we used that the pre-factor \(C_N\) is exponentially trivial. For fixed \(\eta_0>0\), the preceding estimate holds for every sufficiently small \(\delta_0>0\). Since \(S\) is continuous on \(\mathcal{D}_r\times\R\), letting \(\delta_0 \to0_+\) gives
\[
\begin{split}
\liminf_{N \to \infty} \frac{1}{N} \log \E \left [ \textnormal{Crt}_N^{\textnormal{max}}(M,B) \right ] & \geq  S(\boldsymbol{m}_0,x_0)+  \Omega(t _0) - L(\boldsymbol{\gamma}_0,t_0) - \eta_0 \\
& = \Sigma^{\textnormal{max}}(\boldsymbol{m}_0,x_0) - \eta_0\\
& \geq \sup_{(\boldsymbol{m},x) \in M^\circ \times B^\circ} \Sigma^{\textnormal{max}}(\boldsymbol{m}, x) - \varepsilon_0 - \eta_0,
\end{split}
\]
where the last inequality follows by~\eqref{Sigma min ineq}. Finally, letting \(\eta_0\to 0^+\) and then \(\varepsilon_0\to 0^+\) proves the result.
\end{proof}

It remains to prove the three intermediate results, namely Lemmas~\ref{lem: GRF max},~\ref{lem: exp tight max}, and~\ref{lem: det J_N max}.

\begin{proof}[Proof of Lemma~\ref{lem: GRF max}]
We first show that 
\[
\Sigma^\textnormal{max}(\boldsymbol{m},x) = \Sigma^\textnormal{tot}(\boldsymbol{m},x) - L(\boldsymbol{\gamma}(\boldsymbol{m}), t(\boldsymbol{m},x))
\]
is upper semicontinuous. By Lemma~\ref{lem: GRF}, the function \(\Sigma^\textnormal{tot}\) is upper semicontinuous. It therefore suffices to show that the map \((\boldsymbol{m},x) \mapsto L(\boldsymbol{\gamma}(\boldsymbol{m}), t(\boldsymbol{m},x))\) is lower semicontinuous. 

We first verify that \(L\) is lower semicontinuous on \(\R^r\times\R\). According to Definition~\ref{def: function L}, \(L(\boldsymbol{\gamma},t)=+\infty\) for \(t < 2\). Thus lower semicontinuity is immediate at points with \(t<2\), since every sequence converging to such a point eventually satisfies \(t_n<2\). Moreover, sequences approaching a point with \(t=2\) through values \(t_n<2\) cannot violate lower semicontinuity, since the corresponding values of \(L\) are \(+\infty\). It remains to consider the case \(t\ge2\). Since \(L\) is symmetric in the coordinates of \(\boldsymbol{\gamma}\), we may arrange them in decreasing order \(\gamma_1\ge\cdots\ge\gamma_r\). When verifying lower semicontinuity along a convergent sequence, we may pass to a subsequence for which the same ordering holds. If \(\gamma_1<1\), then all coordinates remain below \(1\) in a neighborhood of \(\boldsymbol{\gamma}\), and hence \(L(\boldsymbol{\gamma},t)=0\) there. Otherwise, there exists \(\ell\in\{1,\ldots,r\}\) such that 
\[ 
\gamma_1\geq\cdots\geq\gamma_\ell>1 \geq\gamma_{\ell+1}\geq\cdots\geq\gamma_r, 
\] 
with the case of equality at \(1\) treated below. Away from the boundaries \(\gamma_i=1\) and \(t=\gamma_i+\gamma_i^{-1}\), the function \(L\) is continuous, since it is a finite sum of terms \(I_{\gamma_i}(t)\), where
\[ 
I_{\gamma}(t) = \frac{1}{4}\int_{\gamma+\gamma^{-1}}^t \sqrt{y^2-4}\,\textnormal{d}y -\frac{1}{2}\gamma\left(t-(\gamma+\gamma^{-1})\right) +\frac{1}{8} \left(t^2-\left(\gamma+\gamma^{-1}\right)^2\right). 
\] 
The function \(I_\gamma(t)\) is continuous in \((\gamma,t)\) for \(\gamma>1\) and \(t\geq2\). Moreover, \(I_\gamma\left(\gamma+\gamma^{-1}\right)=0\), so no discontinuity occurs when \(t\) crosses the threshold \(\gamma+\gamma^{-1}\). It remains to consider the boundary \(\gamma=1\). Let \(\gamma_n>1\) and \(t_n\geq2\) satisfy \( \gamma_n\to 1\) and \(t_n\to t\). If \(t>2\), then \(\gamma_n+\gamma_n^{-1}\to2<t\), and hence the corresponding contribution to \(L\) vanishes for all sufficiently large \(n\). If \(t=2\), the explicit formula above implies that \( I_{\gamma_n}(t_n)\to0\) whenever this contribution is present. Thus each contribution is lower semicontinuous across the boundary \(\gamma_i=1\). Consequently, \(L\) is lower semicontinuous on \(\R^r\times\R\). Since \(t(\boldsymbol{m},x)\) is continuous in \((\boldsymbol{m},x)\) by \eqref{function t}, and each \(\gamma_i(\boldsymbol{m})\) is continuous by Subsection~\ref{preliminaries}, the composition 
\[ 
(\boldsymbol{m},x) \mapsto L(\boldsymbol{\gamma}(\boldsymbol{m}),t(\boldsymbol{m},x)) 
\] 
is lower semicontinuous.

Finally, we show that the sub-level sets of \(\Sigma^\textnormal{max}\), namely \(\{\Sigma^\textnormal{max} \geq -a\}\), are compact for all \(a < \infty\). Observe that 
\[
\{\Sigma^\textnormal{max} \geq -a\}  \subseteq \{\Sigma^\textnormal{tot} \geq -a\},
\]
since 
\[
- a \leq \Sigma^\textnormal{max}(\boldsymbol{m},x) = \Sigma^\textnormal{tot}(\boldsymbol{m},x) - L(\boldsymbol{\gamma}(\boldsymbol{m}), t(\boldsymbol{m},x)) \leq \Sigma^\textnormal{tot}(\boldsymbol{m},x), 
\]
where we used that \(L(\boldsymbol{\gamma}(\boldsymbol{m}), t(\boldsymbol{m},x))\ge 0\). By Lemma~\ref{lem: GRF}, the set \(\{\Sigma^\textnormal{tot} \geq -a\}\) is compact. Since \(\Sigma^\textnormal{max}\) is upper semicontinuous, the set \(\{\Sigma^\textnormal{max} \geq -a\}\) is closed. Therefore, it is a closed subset of a compact set, hence compact. This shows that \(-\Sigma^\textnormal{max}\) is a good rate function.
\end{proof}

\begin{proof}[Proof of Lemma~\ref{lem: exp tight max}]
Lemma~\ref{lem: exp tight max} follows from the exponential tightness of the expected number of critical points, see Lemma~\ref{lem: exp tight}. 
\end{proof}

The proofs of the upper and lower bounds in Lemma~\ref{lem: det J_N max} follow the same general strategy as the proof of~\cite[Proposition 4.10]{MR4011861}. 

\begin{proof}[Proof of Lemma~\ref{lem: det J_N max} (Upper bound)]
We equip \(\mathcal{P}(\R)\) with the bounded-Lipschitz distance \(d_{\text{BL}}\). For \(\delta >0\), define
\[
\textnormal{B}(\mu_{\textnormal{sc}}, \delta) \coloneqq  \{\mu \in \mathcal{P}(\R) \colon d_{\text{BL}}(\mu,\mu_{\textnormal{sc}}) < \delta\}.
\]
Choose \(M>\max\left\{\sup_{t\in\mathcal{T}}|t|,2\right\}\). Define
\[
\textnormal{B}_M(\mu_{\textnormal{sc}}, \delta) \coloneqq \{ \mu \in \textnormal{B}(\mu_{\textnormal{sc}}, \delta) \colon \mathrm{supp}(\mu) \subset [-M,M]\}.
\]
Fix \((\boldsymbol{\gamma},t) \in \, \mathcal{U} \times \mathcal{T}\). We decompose 
\[
\E \left [|\det (\boldsymbol{X}_N^{\boldsymbol{\gamma}} - t \boldsymbol{I}_N)|\mathbf{1}_{\{ \boldsymbol{X}_N^{\boldsymbol{\gamma}} - t\boldsymbol{I}_N \prec 0\}}\right] \le E_1(\boldsymbol{\gamma},t)  + E_2 (\boldsymbol{\gamma},t)  + E_3 (\boldsymbol{\gamma},t) ,
\]
where
\[
\begin{split}
E_1(\boldsymbol{\gamma},t) & \coloneqq  \E\left [|\det (\boldsymbol{X}_N^{\boldsymbol{\gamma}} - t\boldsymbol{I}_N)| \mathbf{1}_{\{\boldsymbol{X}_N^{\boldsymbol{\gamma}} - t \boldsymbol{I}_N\prec 0, \, \hat{\mu}_{\boldsymbol{X}_N^{\boldsymbol{\gamma}}} \in \, \textnormal{B}(\mu_{\textnormal{sc}}, \delta), \, \max_{i=1}^N |\lambda_i  (\boldsymbol{X}_N^{\boldsymbol{\gamma}} )  | \le M\}} \right]  ,\\
E_2(\boldsymbol{\gamma},t) & \coloneqq \E\left [|\det (\boldsymbol{X}_N^{\boldsymbol{\gamma}}- t\boldsymbol{I}_N)| \mathbf{1}_{\{\hat{\mu}_{\boldsymbol{X}_N^{\boldsymbol{\gamma}}} \notin \, \textnormal{B}(\mu_{\textnormal{sc}}, \delta)\}} \right], \\
E_3(\boldsymbol{\gamma},t) & \coloneqq \E\left [|\det (\boldsymbol{X}_N^{\boldsymbol{\gamma}}- t\boldsymbol{I}_N)| \mathbf{1}_{\{\max_{i=1}^N |\lambda_i (\boldsymbol{X}_N^{\boldsymbol{\gamma}} )| > M\}} \right].
\end{split}
\]
We first consider \( E_2 (\boldsymbol{\gamma}, t)\). Using~\eqref{eq joint density eig}, we write
\[
\begin{split}
E_2 (\boldsymbol{\gamma}, t) & = \int_{\R^N} \prod_{i=1}^N |\lambda_i -t| \mathbf{1}_{\{\hat{\mu}_{\boldsymbol{\lambda}} \notin \, \textnormal{B}(\mu_{\textnormal{sc}}, \delta)\}}  \textnormal{d} \mathbb{P}^{\boldsymbol{\gamma}}_N (\lambda_1, \ldots, \lambda_N) \\
& = \frac{Z_N^0}{Z_N^{\boldsymbol{\gamma}}}  \int_{\R^N} \prod_{i=1}^N |\lambda_i-t| \mathbf{1}_{\{\hat{\mu}_{\boldsymbol{\lambda}} \notin \, \textnormal{B}(\mu_{\textnormal{sc}}, \delta)\}}  I_N (\boldsymbol{\gamma}, \boldsymbol{D}_N (\boldsymbol{\lambda}) ) \textnormal{d} \mathbb{P}^0_N(\lambda_1, \ldots, \lambda_N)  ,
\end{split}
\]
where \(\hat{\mu}_{\boldsymbol{\lambda}} \coloneqq \frac{1}{N} \sum_{i=1}^N \delta_{\lambda_i}\). Since
\[
Z^{\boldsymbol{\gamma}}_N = Z^0_N \int_{\R^N} I_N(\boldsymbol{\gamma}, \boldsymbol{D}_N (\boldsymbol{\lambda})) \textnormal{d}\mathbb{P}^0_N(\lambda_1, \ldots, \lambda_N), 
\]
it follows that
\[
E_2 (\boldsymbol{\gamma}, t) = \frac{\int_{\R^N} \prod_{i=1}^N |\lambda_i-t| \mathbf{1}_{\{\hat{\mu}_{\boldsymbol{\lambda}} \notin \, \textnormal{B}(\mu_{\textnormal{sc}}, \delta)\}}  I_N (\boldsymbol{\gamma}, \boldsymbol{D}_N (\boldsymbol{\lambda}) ) \textnormal{d} \mathbb{P}^0_N (\lambda_1, \ldots, \lambda_N)}{ \int_{\R^N} I_N(\boldsymbol{\gamma}, \boldsymbol{D}_N (\boldsymbol{\lambda})) \textnormal{d}\mathbb{P}^0_N (\lambda_1, \ldots, \lambda_N)} \eqqcolon \frac{I_1(\boldsymbol{\gamma},t)}{I_2(\boldsymbol{\gamma})}.
\]
By~\cite[Lemma A.4]{MR4011861}, the numerator is exponentially vanishing on compact sets of \( (\boldsymbol{\gamma},t)\) and the denominator is exponentially finite on compact sets of \(\boldsymbol{\gamma}\), i.e., 
\[
\lim_{N \to \infty} \sup_{(\boldsymbol{\gamma},t) \in \, \mathcal{U} \times \mathcal{T}}  \frac{1}{N} \log I_1 (\boldsymbol{\gamma},t) = - \infty, \qquad \limsup_{N \to \infty} \sup_{\boldsymbol{\gamma} \in \, \mathcal{U}}  \left | \frac{1}{N} \log I_2 (\boldsymbol{\gamma})\right| <  \infty.
\]
Hence 
\[
\lim_{N \to \infty} \sup_{(\boldsymbol{\gamma},t) \in \, \mathcal{U} \times \mathcal{T}} \frac{1}{N} \log E_2 (\boldsymbol{\gamma}, t) = - \infty.
\]
Similarly, rewriting \(E_3\) under \(\mathbb{P}_N^0\) as above and using the same uniform exponential lower bound on the denominator,~\cite[Lemma A.4]{MR4011861} gives
\[
\lim_{M \to \infty} \limsup_{N \to \infty} \sup_{(\boldsymbol{\gamma},t) \in \, \mathcal{U} \times \mathcal{T}} \frac{1}{N} \log E_3 (\boldsymbol{\gamma}, t) = - \infty.
\]

It therefore remains to estimate \(E_1 (\boldsymbol{\gamma},t)\). Since \(\boldsymbol{X}_N^{\boldsymbol{\gamma}}-t \boldsymbol{I}_N\prec 0\) is equivalent to \(\lambda_{\max}(\boldsymbol{X}_N^{\boldsymbol{\gamma}})<t\), we have
\[
E_1 (\boldsymbol{\gamma},t) =  \int_{\R^N} \prod_{i=1}^N |\lambda_i  -t| \mathbf{1}_{\{\lambda_{\max} < t, \, \hat{\mu}_{\boldsymbol{\lambda}} \in \, \textnormal{B}(\mu_{\textnormal{sc}}, \delta), \, \max_{i=1}^N |\lambda_i| \le M\}} \,\textnormal{d}\mathbb{P}^{\boldsymbol{\gamma}}_N (\lambda_1, \ldots, \lambda_N).
\]
On the event \(\left \{ \hat{\mu}_{\boldsymbol{\lambda}} \in \, \textnormal{B}(\mu_{\textnormal{sc}}, \delta),  \max_{i=1}^N |\lambda_i| \le M\right \}\), we have that \(\textnormal{supp}(\hat{\mu}_{\boldsymbol{\lambda}} ) \subset [-M,M]\), hence \(\hat{\mu}_{\boldsymbol{\lambda}} \in \, \textnormal{B}_M(\mu_{\textnormal{sc}}, \delta)\). Therefore
\begin{equation} \label{eq:boundE1}
\begin{split} 
E_1 (\boldsymbol{\gamma},t) & \le \int_{\R^N} \exp \{ N \Omega_{\hat{\mu}_{\boldsymbol{\lambda}}} (t)\} \mathbf{1}_{\{\hat{\mu}_{\boldsymbol{\lambda}} \in \, \textnormal{B}_M (\mu_{\textnormal{sc}}, \delta), \, \lambda_{\max} < t\}} \,\textnormal{d}\mathbb{P}^{\boldsymbol{\gamma}}_N (\lambda_1, \ldots, \lambda_N)\\
& \leq   \exp \left \{ N \sup_{\mu \in \, \textnormal{B}_M (\mu_{\textnormal{sc}},\delta)} \Omega_\mu(t) \right \} \mathbb{P}^{\boldsymbol{\gamma}}_N \left (\lambda_{\max} \le t \right),
\end{split}
\end{equation}
where \(\Omega_\mu(x)\) denotes the logarithmic potential and is introduced in Definition~\ref{def: log-potential}. Hence
\[
\frac{1}{N} \log E_1 (\boldsymbol{\gamma},t) \le  \sup_{\mu \in \, \textnormal{B}_M (\mu_{\textnormal{sc}},\delta)} \Omega_\mu(t) + \frac{1}{N} \log  \mathbb{P}^{\boldsymbol{\gamma}}_N \left (\lambda_{\max} \leq t \right) .
\]
For \(1\leq i\leq r\), set \(U_i\coloneqq\min\mathcal{U}_i\), \(\boldsymbol{U}\coloneqq(U_1,\ldots,U_r)\), and \(T\coloneqq\max\mathcal{T}\). The map \((\boldsymbol{\gamma},s) \mapsto \mathbb{P}_N^{\boldsymbol{\gamma}} (\lambda_{\max}\leq s) \) is nonincreasing in each coordinate of \(\boldsymbol{\gamma}\) and nondecreasing in \(s\). Indeed, monotonicity in \(s\) is straightforward. Moreover, if \(\boldsymbol{\gamma}\le \boldsymbol{\gamma}'\) coordinatewise, then
\[
\boldsymbol{X}_N^{\boldsymbol{\gamma}'} -\boldsymbol{X}_N^{\boldsymbol{\gamma}} = \sum_{i=1}^r(\gamma_i'-\gamma_i)e_i e_i^\top \succeq 0.
\]
Thus, by Weyl's monotonicity theorem, \(\lambda_{\max}(\boldsymbol{X}_N^{\boldsymbol{\gamma}'}) \ge \lambda_{\max}(\boldsymbol{X}_N^{\boldsymbol{\gamma}})\), and therefore
\[
\mathbb{P}_N^{\boldsymbol{\gamma}'}(\lambda_{\max}\le s) \le \mathbb{P}_N^{\boldsymbol{\gamma}}(\lambda_{\max}\le s).
\]
Since \(\boldsymbol{U}\leq\boldsymbol{\gamma}\) coordinatewise and \(t\leq T\) for every \((\boldsymbol{\gamma},t)\in\mathcal{U}\times\mathcal{T}\), we obtain 
\[ 
\mathbb{P}_N^{\boldsymbol{\gamma}} (\lambda_{\max}\leq t) \leq \mathbb{P}_N^{\boldsymbol{U}} (\lambda_{\max}\leq T). 
\] 
By Corollary~\ref{cor LDP max eig}, we obtain
\begin{equation} \label{eq:UBLDP}
\begin{split}
\limsup_{N \to \infty}  \sup_{(\boldsymbol{\gamma},t) \in \, \mathcal{U} \times \mathcal{T}}  \frac{1}{N} \log \mathbb{P}^{\boldsymbol{\gamma}}_N \left ( 
\lambda_{\max}  \leq t \right) & \le \limsup_{N \to \infty}  \frac{1}{N} \log \mathbb{P}^{\boldsymbol{U}}_N \left ( \lambda_{\max}  \leq T \right) \\
& \le - L (\boldsymbol{U}, T) \\
& \le - \inf_{(\boldsymbol{\gamma},t) \in \,  \mathcal{U} \times \mathcal{T}}  L(\boldsymbol{\gamma},t),
\end{split}
\end{equation}
where the last inequality follows from \((\boldsymbol{U},T)\in\mathcal{U}\times\mathcal{T}\). Combining~\eqref{eq:boundE1} and~\eqref{eq:UBLDP} gives, for every fixed \(M > \max\{\sup_{t \in \mathcal{T} }|t|,2\}\) and \(\delta > 0\),
\[
\limsup_{N \to \infty}  \sup_{ (\boldsymbol{\gamma},t) \in \, \mathcal{U} \times \mathcal{T}}  \frac{1}{N} \log E_1 (\boldsymbol{\gamma},t) \leq   \sup_{ t \in \mathcal{T}}  \sup_{\mu \in \, \textnormal{B}_M (\mu_{\textnormal{sc}},\delta)} \Omega_\mu(t)  - \inf_{(\boldsymbol{\gamma},t) \in \,  \mathcal{U} \times \mathcal{T}}  L(\boldsymbol{\gamma},t).
\]
Next, since \(\mathcal{P}([-M,M]) \times\mathcal{T}\) is compact and \((\mu,x) \mapsto \Omega_\mu(x)\) is upper semicontinuous on this set (see e.g.~\cite{MR1856194}), it follows that
\begin{equation*} 
\limsup_{\delta \to 0} \sup_{ t \in \mathcal{T}} \sup_{\mu \in \, \textnormal{B}_M(\mu_{\textnormal{sc}},\delta)} \Omega_\mu(t) \leq\sup_{t \in \mathcal{T}} \Omega_{\mu_{\text{sc}}}(t)  = \sup_{t \in \mathcal{T}} \Omega(t).
\end{equation*}
Hence, 
\[
\limsup_{\delta \to 0} \limsup_{N \to \infty}  \sup_{ (\boldsymbol{\gamma},t) \in \, \mathcal{U} \times \mathcal{T}}  \frac{1}{N} \log E_1 (\boldsymbol{\gamma},t) \leq  \sup_{t \in \mathcal{T}} \Omega(t) - \inf_{(\boldsymbol{\gamma},t) \in \,  \mathcal{U} \times \mathcal{T}}  L(\boldsymbol{\gamma},t).
\]
Combining the estimate for \(E_1\) with the exponential negligibility of \(E_2\) and \(E_3\), and letting first \(\delta\to0\) and then \(M\to\infty\), we conclude that
\begin{align*}
& \limsup_{N \to \infty}  \sup_{ (\boldsymbol{\gamma},t) \in \, \mathcal{U} \times \mathcal{T} }  \frac{1}{N} \log \E\left [|\det (\boldsymbol{X}_N^{\boldsymbol{\gamma}} -t\boldsymbol{I}_N)| \mathbf{1}_{\{\boldsymbol{X}_N^{\boldsymbol{\gamma}}-t\boldsymbol{I}_N \prec 0\}} \right] \\
& \leq  \sup_{t \in \mathcal{T}} \Omega(t) - \inf_{(\boldsymbol{\gamma},t) \in \,  \mathcal{U} \times \mathcal{T}}  L(\boldsymbol{\gamma},t).
\end{align*}
This proves~\eqref{UB max}.
\end{proof}

\begin{proof}[Proof of Lemma~\ref{lem: det J_N max} (Lower bound)]
Fix \(\eta_0>0\), \(\boldsymbol{\gamma}_0\in\R^r\), and \(t_0>2\). Let \(\rho_0,\varepsilon_0>0\), to be chosen later, be such that \(t_0 - \rho_0 -\varepsilon_0> 2\). Fix also \(\delta >0\), and let \(M>0\) be chosen sufficiently large below. We write 
\begin{equation*}
\begin{split}
&\E \left [|\det(\boldsymbol{X}_N^{\boldsymbol{\gamma}}-t \boldsymbol{I}_N)| \mathbf{1}_{\{\boldsymbol{X}_N^{\boldsymbol{\gamma}}-t\boldsymbol{I}_N \prec 0\}}\right]  \\
& = \int_{\R^N} \prod_{i=1}^N |\lambda_i -t| \mathbf{1}_{ \{ \lambda_{\max}  < t \} } \text{d}\mathbb{P}^{\boldsymbol{\gamma}}_N (\lambda_1, \ldots, \lambda_N)  \\
&\geq  \int_{\R^N} \exp \left \{N \Omega_{\hat{\mu}_{\boldsymbol{\lambda}}} (t) \right \} \mathbf{1}_{\{ \lambda_{\max} < \min\{ t, t_0-\rho_0 - \varepsilon_0 \}, \lambda_{\min} \ge -M , \hat{\mu}_{\boldsymbol{\lambda}}  \in \, \textnormal{B}(\mu_{\textnormal{sc}}, \delta)\}} \text{d}\mathbb{P}^{\boldsymbol{\gamma}}_N (\lambda_1, \ldots, \lambda_N)  \\
& \geq \exp \left \{ N  \inf_{\substack{\mu \in \, \textnormal{B}(\mu_{\textnormal{sc}}, \delta)\\ \text{supp}(\mu) \subset [-M, \, \min\{ t, t_0-\rho_0 - \varepsilon_0 \}]}} \Omega_\mu(t) \right \}  \\
& \quad \times \mathbb{P}^{\boldsymbol{\gamma}}_N \left (  \lambda_{\max} < \min\{ t, t_0-\rho_0 -\varepsilon_0  \}, \lambda_{\min} \ge -M, \hat{\mu}_{\boldsymbol{\lambda}} \in \, \textnormal{B}(\mu_{\textnormal{sc}}, \delta)\right).
\end{split}
\end{equation*}
By the union bound,
\begin{equation} \label{eq:intersection_LB}
\begin{split}
& \mathbb{P}^{\boldsymbol{\gamma}}_N \left (  \lambda_{\max} < \min\{ t, t_0-\rho_0 -\varepsilon_0 \}, \lambda_{\min} \ge -M, \hat{\mu}_{\boldsymbol{\lambda}} \in \, \textnormal{B}(\mu_{\textnormal{sc}}, \delta)\right) \\
& \ge \mathbb{P}^{\boldsymbol{\gamma}}_N \left (  \lambda_{\max} < \min\{ t, t_0-\rho_0  -\varepsilon_0\} \right) -  \mathbb{P}^{\boldsymbol{\gamma}}_N \left (  \lambda_{\min} < -M \right) - \mathbb{P}^{\boldsymbol{\gamma}}_N \left (  \hat{\mu}_{\boldsymbol{\lambda}} \notin \, \textnormal{B}(\mu_{\textnormal{sc}}, \delta)\right).
\end{split}
\end{equation}
By standard large deviation estimates for the empirical spectral measure (see e.g.~\cite[Lemma A.4]{MR4011861}), \(\mathbb{P}^{\boldsymbol{\gamma}}_N \left (\hat{\mu}_{\boldsymbol{\lambda}}  \notin \textnormal{B}(\mu_{\textnormal{sc}}, \delta) \right)\) is exponentially vanishing on compact sets of \(\boldsymbol{\gamma}\), i.e.,
\[
\lim_{N \to \infty} \sup_{\boldsymbol{\gamma} \in \, \mathcal{U}_{\rho_0}} \frac{1}{N} \log \mathbb{P}^{\boldsymbol{\gamma}}_N \left (\hat{\mu}_{\boldsymbol{\lambda}}  \notin \textnormal{B}(\mu_{\textnormal{sc}}, \delta) \right) = -\infty.
\] 
Moreover, the left tail is exponentially vanishing:
\begin{equation} \label{eq:left_tail_estimate}
\lim_{M \to \infty} \limsup_{N \to \infty} \sup_{\boldsymbol{\gamma} \in \, \mathcal{U}_{\rho_0}} \frac{1}{N} \log \mathbb{P}^{\boldsymbol{\gamma}}_N \left ( \lambda_{\min} < -M \right) = -\infty.
\end{equation}
Next, observe that if \((\boldsymbol{\gamma},t)\in \mathcal{U}_{\rho_0} \times \mathcal{T}_{\rho_0}\), then \(t \ge t_0 - \rho_0 > t_0 - \rho_0 - \varepsilon_0\), hence \(\min \{t, t_0 - \rho_0 - \varepsilon_0\} = t_0 - \rho_0 - \varepsilon_0\). Moreover, as shown in the proof of the upper bound, the map \((\boldsymbol{\gamma},s) \mapsto \mathbb{P}^{\boldsymbol{\gamma}}_N \left (\lambda_{\max}  \leq s \right)\) is nonincreasing in each coordinate of \(\boldsymbol{\gamma}\) and nondecreasing in \(s\). Therefore, for every \((\boldsymbol{\gamma},t)\in \mathcal{U}_{\rho_0}\times\mathcal{T}_{\rho_0}\),
\[
\mathbb{P}^{\boldsymbol{\gamma}}_N \left(\lambda_{\max}  < \min\{ t, t_0-\rho_0- \varepsilon_0 \} \right ) \ge \mathbb{P}^{\boldsymbol{\gamma}_0 + \rho_0 \boldsymbol{1}}_N \left(\lambda_{\max}  <t_0-\rho_0- \varepsilon_0  \right ),
\]
where \(\boldsymbol{1} = (1, \ldots, 1) \in \R^r\). Thus, by Corollary~\ref{cor LDP max eig},
\begin{equation*}
\begin{split}
&\liminf_{N \to \infty} \inf_{(\boldsymbol{\gamma}, t) \in \, \mathcal{U}_{\rho_0} \times \mathcal{T}_{\rho_0}} \frac{1}{N}\log \mathbb{P}^{\boldsymbol{\gamma}}_N \left(\lambda_{\max}  < \min\{ t, t_0-\rho_0- \varepsilon_0 \} \right )  \\
&\ge - L \left ( \boldsymbol{\gamma}_0 + \rho_0 \boldsymbol{1},t_0-\rho_0- \varepsilon_0  \right) .
\end{split}
\end{equation*}
Since
\[
\liminf_{N \to \infty} \inf_{\boldsymbol{\gamma}\in\mathcal{U}_{\rho_0}} \frac{1}{N} \log \mathbb{P}^{\boldsymbol{\gamma}}_N \left (\lambda_{\max}<t_0-\rho_0-\varepsilon_0 \right ) \ge -L(\boldsymbol{\gamma}_0+\rho_0\boldsymbol{1},t_0-\rho_0-\varepsilon_0)>-\infty,
\]
there exists \(C<\infty\) such that, for all sufficiently large \(N\),
\[
\inf_{\boldsymbol{\gamma}\in\mathcal{U}_{\rho_0}} \mathbb{P}^{\boldsymbol{\gamma}}_N \left(\lambda_{\max}<t_0-\rho_0-\varepsilon_0\right) \ge e^{-CN}.
\]
Using~\eqref{eq:left_tail_estimate}, choose \(M>2\) sufficiently large that
\[
\limsup_{N\to\infty} \sup_{\boldsymbol{\gamma}\in\mathcal{U}_{\rho_0}} \frac{1}{N} \log  \mathbb{P}_N^{\boldsymbol{\gamma}}(\lambda_{\min}<-M) <-C.
\]
Moreover, for each fixed \(\delta>0\), the empirical-measure error is exponentially vanishing. Consequently, it follows from~\eqref{eq:intersection_LB} and the above estimates that
\[
\begin{split}
& \mathbb{P}^{\boldsymbol{\gamma}}_N \left (  \lambda_{\max} < \min\{ t, t_0-\rho_0- \varepsilon_0 \}, \lambda_{\min} \ge -M, \hat{\mu}_{\boldsymbol{\lambda}} \in \, \textnormal{B}(\mu_{\textnormal{sc}}, \delta)\right) \\
& \ge (1-o(1)) \mathbb{P}^{\boldsymbol{\gamma}}_N \left (  \lambda_{\max} < t_0-\rho_0- \varepsilon_0 \right) ,
\end{split}
\]
uniformly over \((\boldsymbol{\gamma},t) \in \mathcal{U}_{\rho_0} \times \mathcal{T}_{\rho_0}\). In particular, for all sufficiently large \(N\), 
\[ 
\begin{split} 
&\mathbb{P}^{\boldsymbol{\gamma}}_N \left( \lambda_{\max}<\min\{t, t_0-\rho_0- \varepsilon_0\}, \lambda_{\min}\geq-M, \hat{\mu}_{\boldsymbol{\lambda}} \in\textnormal{B}(\mu_{\textnormal{sc}},\delta) \right)  \geq \frac{1}{2} \mathbb{P}^{\boldsymbol{\gamma}}_N \left( \lambda_{\max}<t_0-\rho_0- \varepsilon_0 \right), 
\end{split} 
\] 
uniformly over \((\boldsymbol{\gamma},t)\in \mathcal{U}_{\rho_0}\times\mathcal{T}_{\rho_0}\). This implies that
\begin{equation} \label{eq:LB-exp}
\begin{split}
& \liminf_{N \to \infty} \inf_{(\boldsymbol{\gamma}, t) \in \mathcal{U}_{\rho_0} \times \mathcal{T}_{\rho_0}} \frac{1}{N} \log \E \left [|\det(\boldsymbol{X}_N^{\boldsymbol{\gamma}}-t \boldsymbol{I}_N)| \mathbf{1}_{\{\boldsymbol{X}_N^{\boldsymbol{\gamma}}-t\boldsymbol{I}_N \prec 0\}}\right]   \\
& \ge \lim_{\delta \to 0_+}  \inf_{t \in \mathcal{T}_{\rho_0}} \inf_{\substack{\mu \in \, \textnormal{B}(\mu_{\textnormal{sc}}, \delta)\\ \textnormal{supp}(\mu) \subset [-M,t_0-\rho_0- \varepsilon_0]}} \Omega_\mu(t) - L \left ( \boldsymbol{\gamma}_0 + \rho_0 \boldsymbol{1},t_0-\rho_0- \varepsilon_0  \right) .
\end{split}
\end{equation}

We now claim that
\[
\lim_{\delta \to 0_+} \inf_{t \in \mathcal{T}_{\rho_0}} \inf_{\substack{\mu \in \, \textnormal{B}(\mu_{\textnormal{sc}}, \delta)\\ \textnormal{supp}(\mu) \subset [-M,t_0-\rho_0- \varepsilon_0]}} \Omega_\mu(t) = \inf_{t \in \mathcal{T}_{\rho_0}} \Omega (t).
\]
Since \(\mu_{\text{sc}}\) is supported on \([-2,2] \subset [-M,t_0-\rho_0-\varepsilon_0]\), we have
\[
\inf_{t\in\mathcal{T}_{\rho_0}} \inf_{\substack{\mu \in \textnormal{B}(\mu_{\mathrm{sc}},\delta)\\ \textnormal{supp}(\mu)\subset [-M,t_0-\rho_0-\varepsilon_0]}} \Omega_\mu(t) \le \inf_{t\in\mathcal{T}_{\rho_0}} \Omega(t).
\]
On the other hand, for \(x\in [-M,t_0-\rho_0-\varepsilon_0]\) and \(t\in \overline{\mathcal{T}}_{\rho_0} = [t_0-\rho_0,t_0+\rho_0]\), we have
\[
t-x \ge (t_0-\rho_0)-(t_0-\rho_0-\varepsilon_0)=\varepsilon_0 >0.
\]
Hence the map \((x,t) \mapsto \log (t-x)\) is continuous and bounded on the compact set \([-M,t_0-\rho_0-\varepsilon_0] \times \overline{\mathcal{T}}_{\rho_0}\). It follows that, \((\mu,t) \mapsto \Omega_\mu(t)\) is continuous on \(\mathcal{P} ([-M,t_0-\rho_0-\varepsilon_0])\times \overline{\mathcal{T}}_{\rho_0}\). By uniform continuity on this compact set,
\[ 
\sup_{t\in\overline{\mathcal{T}}_{\rho_0}} \sup_{\substack{ \mu\in\textnormal{B}(\mu_{\textnormal{sc}},\delta)\\ \operatorname{supp}(\mu) \subset[-M,t_0-\rho_0-\varepsilon_0] }} \left| \Omega_\mu(t)-\Omega(t) \right| \to 0 
\] 
as \(\delta \to 0^+\), which proves the claimed convergence. Moreover, since \(\Omega\) is increasing on \((2,+\infty)\), we have
\[
\inf_{t \in \, \mathcal{T}_{\rho_0}} \Omega (t) = \Omega (t_0-\rho_0).
\]
Thus,~\eqref{eq:LB-exp} becomes
\begin{align*}
& \liminf_{N \to \infty} \inf_{(\boldsymbol{\gamma}, t) \in \mathcal{U}_{\rho_0} \times \mathcal{T}_{\rho_0}} \frac{1}{N} \log \E \left [|\det(\boldsymbol{X}_N^{\boldsymbol{\gamma}} - t \boldsymbol{I}_N)| \mathbf{1}_{\{\boldsymbol{X}_N^{\boldsymbol{\gamma}}-t\boldsymbol{I}_N \prec 0\}}\right]  \\
& \geq \Omega(t_0 - \rho_0) - L \left ( \boldsymbol{\gamma}_0 + \rho_0 \boldsymbol{1} ,t_0-\rho_0- \varepsilon_0  \right).
\end{align*}

By continuity of \((\boldsymbol{\gamma},t)\mapsto \Omega(t)-L(\boldsymbol{\gamma},t)\) at \((\boldsymbol{\gamma}_0,t_0)\), choose \(\rho_0,\varepsilon_0>0\) sufficiently small so that \(t_0-\rho_0-\varepsilon_0>2\) and \((\boldsymbol{\gamma}_0+\rho_0\boldsymbol{1},t_0-\rho_0-\varepsilon_0)\) lies in a sufficiently small neighborhood of \((\boldsymbol{\gamma}_0,t_0)\). Then
\[
\Omega(t_0-\rho_0-\varepsilon_0) - L(\boldsymbol{\gamma}_0+\rho_0\boldsymbol{1},t_0-\rho_0-\varepsilon_0) \ge \Omega(t_0)-L(\boldsymbol{\gamma}_0,t_0)-\eta_0.
\]
Since \(\Omega\) is increasing on \((2,\infty)\),
\[
\Omega(t_0-\rho_0) \ge \Omega(t_0-\rho_0-\varepsilon_0).
\]
Therefore
\[
\Omega(t_0-\rho_0) - L(\boldsymbol{\gamma}_0+\rho_0\boldsymbol 1,t_0-\rho_0-\varepsilon_0) \ge \Omega(t_0)-L(\boldsymbol{\gamma}_0,t_0)-\eta_0.
\]
Hence
\[
\begin{split}
& \liminf_{N \to \infty} \inf_{(\boldsymbol{\gamma}, t) \in \mathcal{U}_{\rho_0} \times \mathcal{T}_{\rho_0}} \frac{1}{N} \log \E \left [|\det(\boldsymbol{X}_N^{\boldsymbol{\gamma}} - t \boldsymbol{I}_N)| \mathbf{1}_{\{\boldsymbol{X}_N^{\boldsymbol{\gamma}}-t\boldsymbol{I}_N \prec 0\}}\right]  \\
& \geq \Omega(t_0)-L(\boldsymbol{\gamma}_0,t_0)-\eta_0.
\end{split}
\]
This proves~\eqref{LB max}.
\end{proof}

\section{Analysis of the variational formula for total complexity} \label{analysis complexity function}

In this section, we prove the results stated in Subsection~\ref{analyzing the variational formulas}. We first analyze the variational problem associated with the total complexity and prove Theorem~\ref{thm variational problem crt points}. Throughout, we restrict attention to the nonnegative correlation domain
\[
\mathcal{D}_r^+ = \mathcal{D}_r\cap[0,1]^r = \left \{\boldsymbol{m} \in [0,1]^r \colon \| \boldsymbol{m} \|_2  <1 \right \} . 
\]
For $\boldsymbol{m}\in\mathcal{D}_r^+$, define \(\Sigma^{\textnormal{tot}}(\boldsymbol{m}) \coloneqq \max_{x \in \R} \Sigma^{\textnormal{tot}}(\boldsymbol{m},x)\) which gives the average number of critical points at a given correlation vector \(\boldsymbol{m}\). We also recall the following quantities from Theorem~\ref{thm variational problem crt points}: 
\[
\tau(\boldsymbol{m} )=\frac{1}{p}\sum_{i=1}^r \lambda_i k_i m_i^{k_i} \quad \text{and} \quad \tau_{\textnormal{c}} (p) = \frac{p-2}{\sqrt{2p(p-1)}}. 
\]
Note that \(\tau(\boldsymbol{m}) \geq 0\) for any \(\boldsymbol{m} \in \mathcal{D}_r^+\). 

The complexity function \(\Sigma^{\textnormal{tot}}(\boldsymbol{m})\) has the following explicit expression.

\begin{lem} \label{explicit formula for Sigma tot}
The function \(\Sigma^{\textnormal{tot}}(\boldsymbol{m})\) has the following explicit formula for every \(\boldsymbol{m} \in \mathcal{D}_r^+\):
\begin{equation*} 
\Sigma^{\textnormal{tot}}(\boldsymbol{m}) = 
\begin{cases}
\Sigma^{\textnormal{tot}}_S(\boldsymbol{m}) & \text{if} \enspace \tau(\boldsymbol{m}) < \tau_{\textnormal{c}} (p), \\
\Sigma^{\textnormal{tot}}_L(\boldsymbol{m}) & \text{if} \enspace \tau(\boldsymbol{m}) \geq \tau_{\textnormal{c}} (p),
\end{cases}
\end{equation*}
where the functions \(\Sigma^{\textnormal{tot}}_S(\boldsymbol{m})\) and \(\Sigma^{\textnormal{tot}}_L(\boldsymbol{m})\) are given by
\[
\begin{split}
\Sigma^{\textnormal{tot}}_S(\boldsymbol{m}) & = \frac{1}{2}\log (p-1) + \frac{1}{2}\log \left (1- \| \boldsymbol{m}\|_2^2 \right) - \frac{1}{p} \sum_{i=1}^r \lambda_i^2 k_i^2  m_i^{2 k_i -2}  \\
& \quad + \frac{p-1}{p(p-2)} \left ( \sum_{i=1}^r \lambda_i k_i m_i^{k_i} \right )^2
\end{split}
\]
and
\[
\begin{split}
\Sigma^{\textnormal{tot}}_L(\boldsymbol{m}) & = \frac{1}{2}\log \left (1- \| \boldsymbol{m}\|_2^2 \right)- \frac{1}{p} \sum_{i=1}^r \lambda_i^2 k_i^2 m_i^{2k_i-2} +   \left ( \frac{1}{\sqrt{2p}}  \sum_{i=1}^r \lambda_i k_i m_i^{k_i} \right )^2 \\
& \quad + \sinh^{-1} \left (  \frac{1}{\sqrt{2p}}  \sum_{i=1}^r \lambda_i k_i m_i^{k_i} \right ) \\
& \quad + \left ( \frac{1}{\sqrt{2p}}  \sum_{i=1}^r \lambda_i k_i m_i^{k_i} \right) \cdot \sqrt{1 + \left ( \frac{1}{\sqrt{2p}}  \sum_{i=1}^r \lambda_i k_i m_i^{k_i} \right )^2},
\end{split}
\]
respectively. The two expressions coincide when $\tau(\boldsymbol{m})=\tau_{\textnormal{c}}(p)$.
\end{lem}

We note that \(\Sigma^{\textnormal{tot}}_S\) and \(\Sigma^{\textnormal{tot}}_L\) coincide at \(\tau(\boldsymbol{m})=\tau_{\textnormal{c}}(p)\). Hence \(\Sigma^{\textnormal{tot}}\) is continuous on \(\mathcal{D}_r^+\). We refer to $\Sigma_S^{\textnormal{tot}}$ and $\Sigma_L^{\textnormal{tot}}$ as the small-$\tau$ and large-$\tau$ branches, respectively. Additionally, we remark that Lemma~\ref{explicit formula for Sigma tot} reduces to Proposition 2.3 in~\cite{MR4011861} when \(r=1\) and \(k=p\). This implies that the analysis of the complexity function \(\Sigma^{\textnormal{tot}}(\boldsymbol{m})\) presented here generalizes the results obtained in~\cite[Subsection 2.4]{MR4011861} for the spiked tensor model.

\begin{proof}
Recall that \(\Sigma^{\textnormal{tot}}(\boldsymbol{m},x) = \Sigma (\boldsymbol{m}, y(\boldsymbol{m},x))\), where \( y(\boldsymbol{m},x) = x - \sum_{i=1}^r \lambda_i (1-k_i/p)m_i^{k_i}\). Maximizing \(\Sigma^{\textnormal{tot}}(\boldsymbol{m},x)\) over \(x\) is equivalent to maximizing \(\Sigma (\boldsymbol{m}, y)\) over \(y\). According to Definition~\ref{def: function Sigma tot}, \(\Sigma (\boldsymbol{m}, y)\) is given by
\begin{equation} \label{eq: sigma new}
\begin{split}
\Sigma (\boldsymbol{m}, y) & = \frac{1}{2} \log(p-1) + \frac{1}{2} \log \left (1- \| \boldsymbol{m}\|_2^2 \right) - \frac{1}{p} \sum_{i=1}^r \lambda_i^2 k_i^2 m_i^{2k_i-2} \\
& \quad + \frac{p-1}{p^2}\left( \sum_{i=1}^r \lambda_i k_i m_i^{k_i} \right)^2  - \frac{p-2}{2(p-1)} y^2 + \frac{2}{p} \left( \sum_{i=1}^r \lambda_i k_i m_i^{k_i} \right)y  \\
& \quad - \frac{1}{2} \int_2^{\sqrt{\frac{2p}{p-1}}|y|} \sqrt{t^2-4} \, \textnormal{d}t \cdot \textbf{1}_{\{\sqrt{\frac{2p}{p-1}}|y| \geq 2\}},
\end{split}
\end{equation}
where we used the following identity:
\[
\int_2^x \sqrt{t^2-4} \, \textnormal{d}t = \frac{x}{2} \sqrt{x^2-4} - 2 \log \left (  \frac{x + \sqrt{x^2-4}}{2} \right ).
\]
Thus, maximizing \(\Sigma(\boldsymbol{m} ,y)\) over \(y\) corresponds to maximizing 
\begin{equation} \label{max eq}
\begin{split}
- \frac{p-2}{2(p-1)} y^2 + \frac{2}{p} \left ( \sum_{i=1}^r \lambda_i k_i m_i^{k_i} \right) y -\frac{1}{2} \int_{2}^{\sqrt{\frac{2p}{p-1}}|y|} \sqrt{t^2-4} \, \textnormal{d}t \cdot \mathbf{1}_{\{\sqrt{2p/(p-1)}|y| \geq 2\}}.
\end{split}
\end{equation}
Set \(u \coloneqq \sqrt{\frac{2p}{p-1}} y\). Then, the optimization problem~\eqref{max eq} is equivalent to solving \(-\frac{1}{2} \min_u g(u)\), where the function \(g\) is given by
\[
g(u) = a u^2 - b u + \int_{2}^{|u|} \sqrt{t^2-4} \, \textnormal{d}t \cdot \mathbf{1}_{\{|u| \geq 2\}},
\] 
with coefficients \(a = \frac{p-2}{2p} >0\) and \(b = \frac{4}{p} \sqrt{\frac{p-1}{2p}} \sum_{i=1}^r \lambda_i k_i m_i^{k_i}>0\). Since the derivative of \(g\),
\[
g'(u) = 2au -b + \text{sgn}(u) \sqrt{u^2-4} \cdot  \mathbf{1}_{\{|u| \geq 2\}},
\]
is monotone increasing, \(g\) has a unique minimum, which occurs when
\[
b-2au = \text{sgn}(u) \sqrt{u^2-4} \cdot  \mathbf{1}_{\{|u| \geq 2\}}.
\]
A simple computation shows that if \(-2 \leq \frac{b}{2a} \leq 2\), then the minimum \(u_{\min}\) occurs at \(\frac{b}{2a}\). Otherwise, it occurs at
\[
u_{\min} =  \frac{2ab-\sqrt{b^2+4-16a^2}}{4a^2-1}.
\]
This leads to
\[
-\frac{1}{2} \min_u g(u) = 
 \begin{cases}
\frac{b^2}{8a} & \text{if} \enspace b \leq 4a,\\
 \frac{1}{4}b u_{\min} + \log \left ( \left (\frac{1}{2} -a \right )u_{\min}+\frac{1}{2}b\right ) & \text{if} \enspace b \geq 4a.
\end{cases}
\]
With our notation for \(a\) and \(b\), when \(b \leq 4a\), which corresponds to \(\tau(\boldsymbol{m} ) \leq \tau_{\textnormal{c}}\), then \(-\frac{1}{2} \min_u g(u)\) equals to 
\[
\frac{b^2}{8a} = 2\frac{p-1}{p^2(p-2)} \left (\sum_{i=1}^r \lambda_i k_i m_i^{k_i}\right)^2.
\]
Combining this with the terms in~\eqref{eq: sigma new} that do not depend on \(y\), we obtain the desired expression for \(\Sigma_S^{\textnormal{tot}}\). On the other hand, when \(b \geq 4a\), i.e., \(\tau(\boldsymbol{m} ) \geq \tau_{\textnormal{c}}\), the minimum \(u_{\min}\) is given by
\[
u_{\min}= \frac{p}{\sqrt{p-1}} \sqrt{1 + \left ( \frac{1}{\sqrt{2p}} \sum_{i=1}^r \lambda_i k_i m_i^{k_i} \right)^2} - \frac{p-2}{\sqrt{p-1}}  \left ( \frac{1}{\sqrt{2p}} \sum_{i=1}^r \lambda_i k_i m_i^{k_i} \right),
\]
so that \(-\frac{1}{2} \min_u g(u)\) results in
\[
\begin{split}
 \frac{1}{4}b u_{\min} + \log \left ( \left (\frac{1}{2} -a \right )u_{\min}+\frac{1}{2}b\right ) &= \left ( \frac{1}{\sqrt{2p}}\sum_{i=1}^r \lambda_i k_i m_i^{k_i}\right) \sqrt{ \left (  \frac{1}{\sqrt{2p}} \sum_{i=1}^r \lambda_i k_i m_i^{k_i}\right)^2 +1} \\
 & \quad - \frac{p-2}{p} \left (\frac{1}{\sqrt{2p}}\sum_{i=1}^r \lambda_i k_i m_i^{k_i}\right)^2 \\
& \quad  + \sinh^{-1} \left (  \frac{1}{\sqrt{2p}} \sum_{i=1}^r \lambda_i k_i m_i^{k_i} \right ) - \frac{1}{2} \log (p-1),
\end{split}
\]
where we used the identity \( \sinh^{-1}(x) = \log(x + \sqrt{1+x^2})\). Combining this with the terms in~\eqref{eq: sigma new} that do not depend on \(y\), we obtain the desired expression for \(\Sigma_L^{\textnormal{tot}}\).
\end{proof}

At $\boldsymbol{m}=\boldsymbol{0}$, the explicit formula gives
\[
\Sigma^{\textnormal{tot}}(\boldsymbol{0}) = \frac{1}{2}\log(p-1)>0.
\]
Hence, in the remainder of this analysis, we assume that $\boldsymbol{m}\in\mathcal{D}_r^+ \backslash \{\boldsymbol{0}\}$.
Next, for \(\boldsymbol{m}\in\mathcal{D}_r^+\backslash \{\boldsymbol{0}\}\), we define the parameters \(\alpha(\boldsymbol{m})\) and \(\beta(\boldsymbol{m} )\) by
\begin{equation} \label{alpha}
\alpha(\boldsymbol{m}) \coloneqq  \| \boldsymbol{m}\|_2^2,
\end{equation}
and
\begin{equation} \label{beta}
\beta(\boldsymbol{m} ) \coloneqq \frac{\alpha(\boldsymbol{m})}{\tau(\boldsymbol{m} )^2} \left (\frac{1}{p^2}\sum_{i=1}^r \lambda_i^2 k_i^2m_i^{2 k_i  -2} \right),
\end{equation}
respectively. Since \(\lambda_i > 0\) for every \(1 \le i \le r\), \(\tau(\boldsymbol{m} ) > 0 \) for every \(\boldsymbol{m}\in\mathcal{D}_r^+\backslash \{\boldsymbol{0}\}\). Moreover, \(\alpha(\boldsymbol{m}) \in (0,1)\) and \(\beta(\boldsymbol{m} ) \geq 1\). Indeed, by the Cauchy--Schwarz inequality, it follows that
\[
\left (\frac{1}{p} \sum_{i=1}^r \lambda_i k_i m_i^{k_i-1} m_i\right )^2 \leq \left ( \sum_{i=1}^r m_i^2 \right ) \cdot  \left (\frac{1}{p^2} \sum_{i=1}^r \lambda_i^2 k_i^2 m_i^{2k_i-2} \right ),
\]
which, in our notation, is equivalent to \(\tau(\boldsymbol{m} )^2 \leq \alpha(\boldsymbol{m}) \frac{\beta(\boldsymbol{m}) \tau(\boldsymbol{m} )^2}{\alpha(\boldsymbol{m})}\). Thus, \(\beta(\boldsymbol{m} ) \geq 1\). 

Our goal is now to identify the regions for \(\boldsymbol{m}\) where the complexity function \(\Sigma^{\textnormal{tot}}(\boldsymbol{m})\) is positive, equal to zero, or negative. To this end, we introduce the parameter \(\tau_\ast (\boldsymbol{m})\):
\[
\tau_\ast (\boldsymbol{m}) \coloneqq  \sqrt{- \frac{1}{2p}\frac{\log \left ( (1-\alpha(\boldsymbol{m}))(p-1) \right)}{\frac{p-1}{p-2} - \frac{\beta(\boldsymbol{m})}{\alpha(\boldsymbol{m})}}}.
\]
Note that \(\tau_\ast(\boldsymbol{m})\) is defined in two cases: either when \(\frac{\alpha(\boldsymbol{m})}{\beta(\boldsymbol{m})} < \alpha(\boldsymbol{m}) \le \frac{p-2}{p-1}\) or when \(\alpha (\boldsymbol{m}) \ge \frac{\alpha(\boldsymbol{m})}{\beta(\boldsymbol{m})} > \frac{p-2}{p-1}\). 

\begin{lem} \label{thm function SU}
The function \(\Sigma^{\textnormal{tot}}_S(\boldsymbol{m})\) has the following properties:
\begin{enumerate}
\item In the region \(\alpha(\boldsymbol{m}) < \frac{p-2}{p-1}\), \(\Sigma^{\textnormal{tot}}_S(\boldsymbol{m})\) is positive for \(\tau(\boldsymbol{m}) <\tau_\ast(\boldsymbol{m})\), zero at \(\tau(\boldsymbol{m}) = \tau_\ast(\boldsymbol{m})\), and negative for \(\tau(\boldsymbol{m}) > \tau_\ast(\boldsymbol{m})\).
\item In the region \(\frac{\alpha(\boldsymbol{m})}{\beta(\boldsymbol{m})} \leq \frac{p-2}{p-1} \leq \alpha(\boldsymbol{m})\), \(\Sigma^{\textnormal{tot}}_S(\boldsymbol{m})\) is nonpositive.
\item In the region \(\frac{\alpha(\boldsymbol{m})}{\beta(\boldsymbol{m})} > \frac{p-2}{p-1}\), \(\Sigma^{\textnormal{tot}}_S(\boldsymbol{m})\) is negative for \(\tau(\boldsymbol{m}) < \tau_\ast(\boldsymbol{m})\), zero at \(\tau(\boldsymbol{m}) = \tau_\ast(\boldsymbol{m})\), and positive for \(\tau(\boldsymbol{m}) > \tau_\ast(\boldsymbol{m})\).
\end{enumerate}
\end{lem}

\begin{proof}
With our choice of \(\alpha\) and \(\beta\), as defined in~\eqref{alpha} and~\eqref{beta}, respectively, a direct substitution shows that \(\Sigma^{\textnormal{tot}}_S(\boldsymbol{m}) = g_{\alpha,\beta} \left ( \sqrt{p} \tau(\boldsymbol{m}) \right )\), where \(g_{\alpha,\beta}\) is defined in Lemma~\ref{lem function S_U}. Therefore, the result follows directly from Lemma~\ref{lem function S_U} by noting that \(\alpha(\boldsymbol{m}) \ge \frac{\alpha(\boldsymbol{m})}{\beta(\boldsymbol{m})}\) since \(\beta(\boldsymbol{m}) \ge 1\).
\end{proof}

\begin{lem} \label{thm function SG}
The function \(\Sigma^{\textnormal{tot}}_L(\boldsymbol{m})\) is nonpositive. In addition, \(\Sigma^{\textnormal{tot}}_L(\boldsymbol{m} ) = 0\) if and only if there exists a constant \(\delta>0\) such that \(\frac{k_i}{p}\lambda_i m_i^{k_i-1} = \delta m_i\) for all \(1 \leq i \leq r\), and the zeros are described by
\begin{equation*}
\tau(\boldsymbol{m} ) =  \frac{1}{\sqrt{2p}} \frac{\alpha(\boldsymbol{m})}{ \sqrt{1-\alpha(\boldsymbol{m})}}.
\end{equation*}
\end{lem}

\begin{proof}
With our choice of \(\alpha\) and \(\beta\), as defined in~\eqref{alpha} and~\eqref{beta}, respectively, a direct substitution shows that \(\Sigma^{\textnormal{tot}}_L(\boldsymbol{m}) = f_{\alpha,\beta} \left ( \sqrt{p/2} \, \tau(\boldsymbol{m} ) \right )\), where \(f_{\alpha,\beta}\) is defined in Lemma~\ref{lem function S_G}. According to Lemma~\ref{lem function S_G}, \(\Sigma^{\textnormal{tot}}_L(\boldsymbol{m})\) is nonpositive for all \(\boldsymbol{m} \in \mathcal{D}_r^+\) and equals zero if and only if 
\begin{equation*} 
\beta(\boldsymbol{m} )=1 \quad \text{and} \quad \sqrt{\frac{p}{2}} \tau(\boldsymbol{m} ) = \frac{1}{2} \frac{\alpha(\boldsymbol{m})}{\sqrt{1-\alpha(\boldsymbol{m})}}.
\end{equation*} 
The condition \(\beta(\boldsymbol{m} )=1\) is equivalent to equality in the Cauchy--Schwarz inequality, meaning that there exists \(\delta>0\) such that \(\frac{k_i}{p} \lambda_i m_i^{k_i-1} = \delta m_i\) for all \(1 \le i \le r\). This completes the proof.
\end{proof}

The following result follows directly from Lemmas~\ref{thm function SU} and~\ref{thm function SG}.

\begin{prop} \label{tot complexity}
Within the region \(\{\tau(\boldsymbol{m}) < \tau_{\textnormal{c}} (p)\} \cap \{ \alpha(\boldsymbol{m}) < \frac{p-2}{p-1}\}\), the complexity function \(\Sigma^{\textnormal{tot}}(\boldsymbol{m})\) is positive for \(0 \leq \tau(\boldsymbol{m}) < \tau_\ast(\boldsymbol{m})\), zero at \(\tau(\boldsymbol{m}) = \tau_\ast (\boldsymbol{m})\), and negative for \(\tau(\boldsymbol{m}) > \tau_\ast(\boldsymbol{m})\). Outside this region, i.e., in \(\{\tau(\boldsymbol{m}) < \tau_{\textnormal{c}} (p)\} \cap \{ \alpha(\boldsymbol{m}) \geq \frac{p-2}{p-1}\}\) as well as in \(\{\tau(\boldsymbol{m}) \geq \tau_{\textnormal{c}} (p)\}\), the complexity \(\Sigma^{\textnormal{tot}}(\boldsymbol{m})\) is nonpositive. 
\end{prop}

\begin{proof}
This result follows directly from combining Lemmas~\ref{thm function SU} and~\ref{thm function SG}. First, note that when \(\tau(\boldsymbol{m}) \ge \tau_{\textnormal{c}} (p)\), \(\Sigma^{\textnormal{tot}} (\boldsymbol{m}) = \Sigma^{\textnormal{tot}}_L (\boldsymbol{m}) \le 0\) according to Lemma~\ref{thm function SG}. Now, consider the region \(\{\tau(\boldsymbol{m}) < \tau_{\textnormal{c}} (p)\}\), where \(\Sigma^{\textnormal{tot}} (\boldsymbol{m}) = \Sigma^{\textnormal{tot}}_S (\boldsymbol{m})\). Within \(\{\alpha(\boldsymbol{m}) < \frac{p-2}{p-1}\}\), item 1 of Lemma~\ref{thm function SU} gives the following: \(\Sigma^{\textnormal{tot}}(\boldsymbol{m})\) is positive for \(\tau(\boldsymbol{m}) <\tau_\ast(\boldsymbol{m})\), zero at \(\tau(\boldsymbol{m}) = \tau_\ast(\boldsymbol{m})\), and negative for \(\tau(\boldsymbol{m}) > \tau_\ast(\boldsymbol{m})\). In the region \(\{\frac{\alpha(\boldsymbol{m})}{\beta(\boldsymbol{m})} \le \frac{p-2}{p-1} < \alpha(\boldsymbol{m}) \}\), item 2 of Lemma~\ref{thm function SU} implies that \(\Sigma^{\textnormal{tot}} (\boldsymbol{m}) \le 0\). Finally, consider the region \(\{\frac{p-2}{p-1} < \frac{\alpha(\boldsymbol{m})}{\beta(\boldsymbol{m})} \le \alpha(\boldsymbol{m}) \}\). We show that in this case \(\tau_\ast(\boldsymbol{m}) > \tau_{\textnormal{c}} (p)\). To see this, observe the following:
\[
\begin{split}
-\log ((1-\alpha(\boldsymbol{m}))(p-1)) & \geq 1-(1-\alpha(\boldsymbol{m}))(p-1) \\
& = \alpha(\boldsymbol{m})(p-1)-(p-2) \\
& \geq \alpha(\boldsymbol{m})(p-1)-\beta(\boldsymbol{m})(p-2),
\end{split}
\]
where we used that \(\log(x) \le x-1\) for \(x>0\) in the first inequality and that \(\beta(\boldsymbol{m}) \geq 1\) in the second. From this, we deduce  
\[
\tau_\ast(\boldsymbol{m}) \geq \sqrt{\frac{ \alpha(\boldsymbol{m})(p-2)}{2p}} > \frac{p-2}{\sqrt{2p(p-1)}} = \tau_{\textnormal{c}} (p), 
\]
where the second inequality follows from \(\alpha(\boldsymbol{m}) \geq \frac{\alpha(\boldsymbol{m})}{\beta(\boldsymbol{m})} > \frac{p-2}{p-1}\). Thus, in this region, item 3 of Lemma~\ref{thm function SU} ensures that \(\Sigma^{\textnormal{tot}}(\boldsymbol{m}) \le 0\). Combining the last two cases, we conclude that within \(\{\tau(\boldsymbol{m}) < \tau_{\textnormal{c}} (p)\} \cap \{ \alpha(\boldsymbol{m}) \geq \frac{p-2}{p-1}\}\), \(\Sigma^{\textnormal{tot}}(\boldsymbol{m}) \le 0\). This completes the proof.
\end{proof}

Having established Proposition~\ref{tot complexity}, we now turn to characterizing the regions in which the complexity function \(\Sigma^{\textnormal{tot}}_L\) vanishes as the parameters \(\lambda_1, \ldots, \lambda_r\) are varied. For simplicity, we focus on the case where \(k_i = k\) for all \(1 \leq i \leq r\); that is, the deterministic polynomials in~\eqref{eq: function f} are of the same degree. We recall the following definition from Theorem~\ref{thm variational problem crt points}:
\[
\eta(\boldsymbol{m}) = \sum_{i=1}^r \lambda_i^{-\frac{2}{k-2}} \mathbf{1}_{\{m_i \neq 0\}}.
\]

\begin{prop} \label{thm function SG 2}
Let \(k_i = k\) for all \(1 \leq i \leq r\), and define
\[
 \tilde{\eta}_{\textnormal{c}}(p,k) \coloneqq (k-2) \left( \frac{2k^2}{p(k-1)^{k-1}} \right)^{\frac{1}{k-2}}.
\]
Then, 
\begin{enumerate}
\item If \(\eta(\boldsymbol{m}) > \tilde{\eta}_{\textnormal{c}}(p,k)\), then \(\Sigma^{\textnormal{tot}}_L(\boldsymbol{m}) < 0\).
\item If \(\eta(\boldsymbol{m}) \leq \tilde{\eta}_{\textnormal{c}}(p,k)\), then \(\Sigma^{\textnormal{tot}}_L(\boldsymbol{m}) \leq 0\), with equality if and only if \(\boldsymbol{m}\) satisfies:
\begin{equation} \label{eq: cond zero 1}
\lambda_i m_i^{k-2} = \lambda_j m_j^{k-2} \enspace \text{for all} \enspace 1 \leq i,j \leq r \enspace \text{such that} \enspace m_i, m_j \neq 0,
\end{equation}
and
\begin{equation} \label{eq: cond zero 2}
\tau(\boldsymbol{m}) = \frac{1}{\sqrt{2p}} \frac{\alpha(\boldsymbol{m})}{\sqrt{1 - \alpha(\boldsymbol{m})}}.
\end{equation}
\end{enumerate}
\end{prop}

\begin{proof}
According to Lemma~\ref{thm function SG}, the function \(\Sigma^{\textnormal{tot}}_L\) vanishes if and only if the following two conditions hold:
\begin{enumerate}
\item[(i)] \(\beta(\boldsymbol{m})=1\), 
\item[(ii)] \(\sqrt{\frac{p}{2}} \tau(\boldsymbol{m}) = \frac{1}{2} \frac{\alpha(\boldsymbol{m})}{\sqrt{1-\alpha(\boldsymbol{m})}}\). 
\end{enumerate}
Condition (i) implies that there exists a scalar \(\delta>0\) such that 
\[
\frac{k}{p} \lambda_i m_i^{k-1} = \delta m_i \quad \text{for all } i \text{ with } m_i \neq 0,
\]
which is equivalent to
\[
\lambda_i m_i^{k-2} = \lambda_j m_j^{k-2} \enspace \text{for all} \enspace 1 \leq i,j \leq r \enspace \text{such that} \enspace m_i, m_j \neq 0.
\]
Substituting this into the definitions of \(\alpha(\boldsymbol{m})\) and \(\tau(\boldsymbol{m})\), we obtain
\begin{equation} \label{equations for alpha and tau}
\alpha(\boldsymbol{m}) = \delta^{\frac{2}{k-2}} \left( \frac{p}{k} \right)^{\frac{2}{k-2}} \eta(\boldsymbol{m}), \qquad
\tau(\boldsymbol{m}) = \delta^{\frac{k}{k-2}} \left( \frac{p}{k} \right)^{\frac{2}{k-2}} \eta(\boldsymbol{m}).
\end{equation}
Plugging these expressions into condition (ii), we obtain the following equation for \(\delta\):
\begin{equation} \label{eq for delta}
\delta^{-2} + 2p \left( \frac{p}{k} \right)^{\frac{2}{k-2}} \eta(\boldsymbol{m}) \, \delta^{\frac{2}{k-2}} - 2p = 0.
\end{equation}
Define the function
\[
h(x) = x^{-2} + 2p \left( \frac{p}{k} \right)^{\frac{2}{k-2}} \eta(\boldsymbol{m}) \, x^{\frac{2}{k-2}} - 2p, \qquad x > 0.
\]
The derivative of $h$ vanishes at a unique point 
\begin{equation} \label{eq: min value h}
x_{\min} = \left( \frac{k}{p} \right)^{\frac{1}{k-1}} \left( \frac{k-2}{2p} \cdot \frac{1}{\eta(\boldsymbol{m})} \right)^{\frac{k-2}{2(k-1)}} > 0,
\end{equation}
Moreover, $h'(x)<0$ for $x<x_{\min}$ and $h'(x)>0$ for $x>x_{\min}$. Hence, $h$ attains its unique global minimum at $x_{\min}$, with minimum value
\[
y_{\min} = h(x_{\min}) = (k-1) \left( \frac{p}{k} \right)^{\frac{2}{k-1}} \left( \frac{2p \eta(\boldsymbol{m})}{k-2} \right)^{\frac{k-2}{k-1}} - 2p.
\]
Now we observe that if \(\eta(\boldsymbol{m}) > \tilde{\eta}_{\textnormal{c}} (p,k) \), then \(y_{\min} > 0\), and consequently \(h(x) > 0\) for all \(x >0\). In this case,~\eqref{eq for delta} admits no solution, and thus \(\Sigma^{\textnormal{tot}}_L (\boldsymbol{m}) < 0\). On the other hand, if \(\eta(\boldsymbol{m}) \leq \tilde{\eta}_{\textnormal{c}} (p,k) \), then \(y_{\min}\leq 0\), so \(h(x)\) has at least one root in \((0,\infty)\). Therefore, there exists \(\delta > 0\) such that~\eqref{eq for delta} is satisfied. In this case, \(\Sigma^{\textnormal{tot}}_L\) vanishes if and only if both conditions (i) and (ii) are satisfied.
\end{proof}

We now provide the proof of Theorem~\ref{thm variational problem crt points}.

\begin{proof}[Proof of Theorem~\ref{thm variational problem crt points}]
Theorem~\ref{thm variational problem crt points} follows from Lemma~\ref{explicit formula for Sigma tot}, Lemma~\ref{thm function SG}, and Proposition~\ref{thm function SG 2}. In particular, from Lemma~\ref{explicit formula for Sigma tot} we know that, on the region
\(\{ \boldsymbol{m} \colon \tau(\boldsymbol{m}) \ge \tau_\textnormal{c}(p)\}\),
\[
\Sigma^{\textnormal{tot}} (\boldsymbol{m}) =\Sigma_L^{\textnormal{tot}} (\boldsymbol{m}).
\]
Since Lemma~\ref{thm function SG} implies \(\Sigma_L^{\textnormal{tot}}(\boldsymbol{m})\le 0\) for all \(\boldsymbol{m}\), the first conclusion follows: 
\[
\Sigma^{\textnormal{tot}} (\boldsymbol{m}) \le 0 ,
\]
for all \(\boldsymbol{m}\) such that \(\tau(\boldsymbol{m}) \ge \tau_\textnormal{c}(p)\). 

To analyze the equality case, recall the function \(h(x)\) 
\[
h(x) = x^{-2} + 2p \left( \frac{p}{k} \right)^{\frac{2}{k-2}} \eta(\boldsymbol{m})  x^{\frac{2}{k-2}} - 2p, \qquad x > 0,
\]
introduced in the proof of Proposition~\ref{thm function SG 2}. The function \(h\) attains its global minimum at \(x_{\min}=x_{\min}(p,k,\eta(\boldsymbol{m}))\) (see~\eqref{eq: min value h}). From Proposition~\ref{thm function SG 2}, \(\Sigma^{\textnormal{tot}}_L(\boldsymbol{m}) = 0\) if and only if \(h(\delta) = 0\) for some \(\delta >0\), and in this case 
\[
\tau (\boldsymbol{m}) = \delta^{\frac{k}{k-2}} \left ( \frac{p}{k} \right)^{\frac{2}{k-2}} \eta (\boldsymbol{m}).
\] 
Furthermore, by Proposition~\ref{thm function SG 2} we know that, if \(\eta(\boldsymbol{m}) \le \tilde{\eta}_\textnormal{c}(p,k)\), then \(h(x) = 0\) has one or two positive roots, hence \(\Sigma^{\textnormal{tot}} (\boldsymbol{m}) \le 0\). Denote by \(x_\ast = x_\ast (p,k,\eta(\boldsymbol{m}))\) the largest such root. Thus, a configuration \(\boldsymbol{m}_0\) with \( \tau(\boldsymbol{m}_0) \ge \tau_{\textnormal{c}}(p)\) satisfies \(\Sigma^{\textnormal{tot}} (\boldsymbol{m}_0) = 0\) if and only if 
\begin{itemize}
\item \(\eta (\boldsymbol{m}_0) \le  \tilde{\eta}_\textnormal{c}(p,k)\),
\item \(\tau (\boldsymbol{m}_0)  = \tau_\ast (\eta (\boldsymbol{m}_0)) \ge \tau_{\mathrm{c}}(p)\), where \(\tau_\ast (\eta) \coloneqq x_\ast (p,k,\eta)^{\frac{k}{k-2}} \left ( \frac{p}{k} \right)^{\frac{2}{k-2}} \eta\).
\end{itemize}
This motivates the definition of the critical threshold \(\eta_\textnormal{c}(p,k)\):
\[
\begin{split}
\eta_\textnormal{c}(p,k) & \coloneqq \sup \left \{ \eta(\boldsymbol{m}) > 0 \colon \exists \, \boldsymbol{m} \: \mathrm{with} \: \tau (\boldsymbol{m})  \ge \tau_c(p) \: \mathrm{and} \: \Sigma^{\textnormal{tot}}_L (\boldsymbol{m})= 0 \right \} \\
&=\sup \left\{ \eta \le \tilde{\eta}_\textnormal{c}(p,k) \colon \tau_\ast (\eta)  \ge \tau_\textnormal{c}(p) \right\}.
\end{split}
\]
If \(\eta(\boldsymbol{m})>\eta_{\mathrm c}(p,k)\), then either \(\eta(\boldsymbol{m})>\tilde \eta_{\mathrm c}(p,k)\), in which case \(\Sigma^{\textnormal{tot}}_L (\boldsymbol{m}) < 0\) by Proposition~\ref{thm function SG 2}, or else any \(\boldsymbol{m}\) with \(\Sigma^{\textnormal{tot}}_L(\boldsymbol{m})=0\) satisfies \(\tau(\boldsymbol{m})<\tau_{\mathrm c}(p)\). Hence for all \(\boldsymbol{m}\) with \(\tau(\boldsymbol{m})\ge\tau_{\mathrm c}(p)\) we have \(\Sigma^{\textnormal{tot}}(\boldsymbol{m})<0\). Conversely, if \(\eta(\boldsymbol{m})\le\eta_{\mathrm c}(p,k)\), then \(\Sigma^{\textnormal{tot}}(\boldsymbol{m})=0\) in the region \(\{\tau(\boldsymbol{m})\ge\tau_{\mathrm{c}}(p)\}\) precisely when conditions~\eqref{eq: cond zero 1} and~\eqref{eq: cond zero 2} hold, as established in Proposition~\ref{thm function SG 2}.

It remains to compute \(\eta_\textnormal{c}(p,k)\). When \(k \ge p\), we have
\[
\tau_\ast (\tilde{\eta}_\textnormal{c}(p,k)) = \frac{k - 2}{\sqrt{2p(k - 1)}} \ge \frac{p - 2}{\sqrt{2p(p - 1)}} = \tau_\textnormal{c}(p),
\]
so \(\eta_{\mathrm c}(p,k)=\tilde\eta_{\mathrm c}(p,k)\). When \(k < p\), the above inequality fails and \(\eta_{\mathrm c}(p,k)<\tilde\eta_{\mathrm c}(p,k)\). In this case, substituting
\[
x_\ast (p,k,\eta) =\left(\frac{\tau_{\mathrm c}(p)}{\eta}\right)^{\frac{k-2}{k}}\left(\frac{k}{p}\right)^{\frac{2}{k}}
\]
into the equation \(h(x_\ast)=0\) and solving for \(\eta\) yields a unique positive solution
\[
\eta_{\mathrm c}(p,k)=(p-2)\left(\frac{2k^2}{p(p-1)^{k-1}}\right)^{\frac{1}{k-2}}.
\]
\end{proof}

As a direct consequence of Theorem~\ref{thm variational problem crt points}, we show that in the special case \(r=1\) and \(k = p\), the critical threshold \(\lambda_\text{c}(k,k)\) for \(\Sigma^{\textnormal{tot}}\) coincides with the threshold for \(\Sigma^{\textnormal{max}}\). We emphasize that the restriction \(k=p\) is purely technical: only in this case can the argument be deduced directly from Theorem~\ref{thm variational problem crt points}. The general case \(k \neq p\) would require additional analysis and is not treated here (see also the discussion in Subsection~\ref{analyzing the variational formulas}).

\begin{lem} \label{lem: threshold for max}
Let \(r = 1\) and define 
\[
m_\textnormal{c} (\lambda) \coloneqq \left ( \frac{1}{\lambda} \frac{p}{k} \frac{p-2}{\sqrt{2p(p-1)}}\right)^{\frac{1}{k}}.
\]
Then, 
\[
\Sigma^{\textnormal{max}}(m) \le 0 \qquad \text{for all } m \in [m_\textnormal{c} (\lambda), 1).
\]
Furthermore, let \(\lambda_\textnormal{c}(p,k)\) denote the critical threshold defined in~\eqref{eq: crit threshold r=1 1}. It holds that
\begin{enumerate}
\item If \(\lambda < \lambda_\textnormal{c}(p,k)\), then \(\Sigma^{\textnormal{max}}(m) < 0\) for all \(m \in [m_\textnormal{c} (\lambda), 1)\);
\item If \(k = p\) and \(\lambda \ge \lambda_\textnormal{c}(k,k)\) , then there exists \(m_0 (\lambda)  \in [m_\textnormal{c} (\lambda), 1)\) such that 
\[
\Sigma^{\textnormal{max}}(m_0(\lambda)) = 0,
\]
and \(m_0 (\lambda)\) satisfies 
\[
m^{2k-4} (1 -m^2) = \frac{1}{2k\lambda^2}.
\] 
\end{enumerate}
\end{lem}

\begin{proof}
Recall from Definition~\ref{def: function Sigma max} that 
\[
\Sigma^{\textnormal{max}}(m,x) =  \Sigma^{\textnormal{tot}}(m,x) - L(\theta(m), t(m,x)) \le \Sigma^{\textnormal{tot}}(m,x),
\]
since \(L \ge 0\). Taking the maximum over \(x\) yields
\[
\Sigma^{\textnormal{max}}(m) \le \Sigma^{\textnormal{tot}}(m).
\]
According to Theorem~\ref{thm variational problem crt points} and in particular Lemma~\ref{thm function SG}, \( \Sigma^{\textnormal{tot}}(m) \le 0\) for all \(m\) such that \(\tau(m) \ge \tau_{\mathrm{c}}(p)\), which is equivalent to \(m \ge m_{\mathrm{c}}(\lambda)\). Moreover, if \(\lambda < \lambda_\textnormal{c} (p,k)\), it follows that \(\Sigma^{\textnormal{max}} (m) \le \Sigma^{\textnormal{tot}}(m) < 0\) for all \(m \in [m_\textnormal{c} (\lambda), 1)\). If we write
\[
\lambda_{\mathrm{s}} (p,k) \coloneqq \inf \left \{\lambda >0 \colon \exists m \in [m_{\mathrm{c}}(\lambda),1)  \enspace \mathrm{with} \enspace \Sigma^{\textnormal{max}} (m)=0 \right\},
\]
then the above implies that \(\lambda_{\mathrm{s}} (p,k) \ge \lambda_{\mathrm{c}}(p,k)\). 

Now consider \(\lambda \ge \lambda_\textnormal{c} (p,k)\). By Theorem~\ref{thm variational problem crt points}, if \(\lambda \ge \lambda_\textnormal{c} (p,k)\), then \(\Sigma^{\textnormal{tot}}(m)\) vanishes on \([m_\textnormal{c} (\lambda), 1)\) exactly at the solutions of
\begin{equation} \label{eq: zero eq for tot compl}
m^{2k - 4}(1 - m^2) = \frac{p}{2k^2\lambda^2}.
\end{equation}
The function \(f(m) \coloneqq m^{2k - 4}(1 - m^2)-\frac{p}{2k^2\lambda^2}\) is strictly increasing on \((0,m_{\max})\) and strictly decreasing on \((m_{\max},1)\), where \(m_{\max} = \sqrt{\frac{k-2}{k-1}}\), and satisfies \(f(0)=f(1) <0\). Hence, for any \(p,k \ge 3\) and \(\lambda \ge \lambda_{\mathrm{c}}(p,k)\),~\eqref{eq: zero eq for tot compl} has one or two solutions \(0<m_1(\lambda)\le m_2(\lambda)<1\) with \(0 < m_1(\lambda) \le m_{\max} \le m_2(\lambda) < 1\). When \(k=p\) and \(\lambda \ge \lambda_{\mathrm{c}}(k,k)\), we have \(m_{\mathrm{c}}(\lambda) \le m_{\max}\), so in particular \(m_2(\lambda) \in [m_\textnormal{c} (\lambda), 1)\). We claim that for \(k = p\), \(m_2(\lambda)\) is also a solution of 
\begin{equation}\label{claim r1}
\Sigma^{\textnormal{max}}(m_2(\lambda)) = 0.
\end{equation}
This shows the second item of the statement since \(\lambda_{\mathrm{s}} (k,k) \le  \lambda_{\mathrm{c}}(k,k)\), and from the first part of the proof we already know \(\lambda_{\mathrm{s}} (k,k) \ge  \lambda_{\mathrm{c}}(k,k)\), thus \(\lambda_{\mathrm{s}} (k,k) =  \lambda_{\mathrm{c}}(k,k)\).

To prove~\eqref{claim r1}, we consider the function \(\theta(m)\) defined in~\eqref{function theta} with \(k=p\):
\[
\theta(m) = \sqrt{2k(k-1)} \lambda m^{k-2}(1 - m^2).
\]
Since \(m_2(\lambda)\) satisfies~\eqref{eq: zero eq for tot compl} and \(m_2 (\lambda) \ge \sqrt{\frac{k-2}{k-1}}\), we obtain
\[
\theta (m_2(\lambda)) =  \sqrt{2k(k-1)} (1-m_2^2) \sqrt{\frac{1}{2k (1-m_2^2)}} =\sqrt{k-1}\sqrt{1-m_2^2} \le 1.
\]
According to Definition~\ref{def: function L} and~\eqref{function t}, it follows that
\[
\Sigma^{\textnormal{max}}(m_2(\lambda), x) = 
\begin{cases}
\Sigma^{\textnormal{tot}}(m_2(\lambda), x) &  \text{if} \enspace t(x) \ge 2,\\
- \infty &  \text{if} \enspace t(x) < 2,
\end{cases}
\]
where \(t(x) \coloneqq t(m_2(\lambda),x)=\sqrt{\frac{2k}{k-1}}x\). Taking the maximum over \(x \in \R\), we obtain 
\[
\Sigma^{\textnormal{max}}(m_2(\lambda)) = \max_x \Sigma^{\textnormal{max}}(m_2(\lambda), x) = \max_{x \colon t(x) \ge 2} \Sigma^{\textnormal{tot}}(m_2(\lambda), x).
\]
By Lemma~\ref{explicit formula for Sigma tot}, for each fixed \(m\), the function \(x \mapsto \Sigma^{\textnormal{tot}} (m,x)\) is strictly concave and admits a unique maximizer \(x_\ast(m)\), and for \(k=p\) this maximizer is given explicitly by
\[
x_\ast (m) = 
\begin{cases}
2 \frac{k-1}{k-2} \lambda m^k &  \text{if} \enspace | t(x_\ast (m)) |\le 2,\\
\sqrt{\frac{k}{2}} \sqrt{1 + \frac{k}{2} \lambda^2 m^{2k}} - \frac{k-2}{2} \lambda m^k&  \text{otherwise} .
\end{cases}
\]
For \(m \ge m_{\mathrm{c}}(\lambda)\), the proof of Lemma~\ref{explicit formula for Sigma tot} (see also the proof of \cite[Proposition 2.3]{MR4011861} in the case \(r=1, k=p\)) shows that the global maximizer is always given by the second branch, corresponding to \(t(x_\ast (m)) \ge 2\), and
\[
\max_{x \colon t(x) \ge 2} \Sigma^{\textnormal{tot}}(m, x) = \Sigma^{\textnormal{tot}}_L (m).
\]
In particular, since \(m_2(\lambda) \ge m_{\mathrm{c}}(\lambda)\), 
\[
\max_{x \colon t(x) \ge 2} \Sigma^{\textnormal{tot}}(m_2(\lambda), x) = \Sigma^{\textnormal{tot}}_L (m_2(\lambda)) =0.
\]
This shows that 
\[
\Sigma^{\textnormal{max}}(m_2(\lambda))  = 0,
\]
as desired.
\end{proof}

\appendix

\section{Auxiliary results for Section~\ref{analysis complexity function}}\label{appendix}

\begin{lem}[Auxiliary result for \(\Sigma^{\textnormal{tot}}_S(\boldsymbol{m})\)]  \label{lem function S_U}
For every value of \(a \in (0,1)\) and \(b \geq 1\), we define \(g_{a,b} \colon \R \to \R\) by
\[ 
g_{a,b}(x)= \frac{1}{2} \log (1-a) + \frac{1}{2}\log(p-1) + \left (\frac{p-1}{p-2} - \frac{b}{a} \right) x^2 .
\]
Moreover, for \(\frac{a}{b} < a \le \frac{p-2}{p-1}\) and for \(\frac{a}{b} > \frac{p-2}{p-1}\), we define \(x_\ast\) by
\[
x_\ast = \sqrt{\frac{-\frac{1}{2} (\log(1-a) + \log(p-1))}{\frac{p-1}{p-2} - \frac{b}{a} }}.
\] 
Then, the following holds:
\begin{enumerate}
\item If \(\frac{a}{b} > \frac{p-2}{p-1}\), then \(g_{a,b}\) is negative in \((-x_\ast, x_\ast)\), zero at \(x=\pm x_\ast\), and positive otherwise.
\item If \(\frac{a}{b} < \frac{p-2}{p-1}\) we distinguish three cases. If \(a <  \frac{p-2}{p-1}\), then \(g_{a,b}(x)\) is positive in \((-x_\ast, x_\ast)\), zero at \(x=\pm x_\ast\), and negative otherwise. If  \(a =  \frac{p-2}{p-1}\), then \(g_{a,b}(x)\) has exactly one zero in \(x_\ast\). Otherwise, if \(a > \frac{p-2}{p-1}\), \(g_{a,b}(x)\) is negative for all \(x\).
\item If \(\frac{a}{b}=  \frac{p-2}{p-1}\), then \(g_{a,b}\) is a nonpositive, constant function.
\end{enumerate}
\end{lem}

\begin{proof}
If \(\frac{a}{b}=  \frac{p-2}{p-1}\), then \(g_{a,b} (x) = \frac{1}{2} \log((1-a)(p-1))\) is a nonpositive, constant function. Otherwise, if \(\frac{a}{b} \neq \frac{p-2}{p-1}\), the function \(g_{a,b}\) is a parabola which is symmetric about the \(y\)-axis, and is opening to the top if \(\frac{a}{b} > \frac{p-2}{p-1}\) and to the bottom if \(\frac{a}{b} < \frac{p-2}{p-1}\). If \(\frac{a}{b} > \frac{p-2}{p-1}\), the function has two zeros since the \(y\)-coordinate of vertex \(y_V = g_{a,b}(0)= \frac{1}{2} \log ((1-a)(p-1)) < 0\). If \(\frac{a}{b} < \frac{p-2}{p-1}\), we then have: if \(a > \frac{p-2}{p-1}\), then \(y_V\) is negative and \(g_{a,b}\) is always negative; if \(a = \frac{p-2}{p-1}\), then \(y_V=0\) and \(g_{a,b}\) has exactly one zero; if \(a < \frac{p-2}{p-1}\), then \(y_V\) is positive, so the function \(g_{a,b}\) has two zeros at \(x = \pm x_\ast\).
\end{proof}

\begin{lem}[Auxiliary result for \(\Sigma^{\textnormal{tot}}_L(\boldsymbol{m})\)] \label{lem function S_G}
For every value of \(a \in (0,1)\) and \(b \geq 1\), we define \(f_{a,b} \colon \R \to \R\) by
\[ 
f_{a,b}(x) = \frac{1}{2} \log (1-a) - \frac{2b}{a} x^2 + x^2 + x \sqrt{1 +x^2} + \sinh^{-1}(x).
\]
Then, the following holds:
\begin{enumerate}
\item \(f_{a,b}(x) \leq 0\) for all \(x\);
\item \(f_{a,b}(x) \) has exactly one maximum at \(x_{a,b}=\frac{1}{2} \frac{a}{\sqrt{b (b-a)}}\);
\item \(f_{a,b}(x_{a,b})=0\) if and only if \(b=1\). 
\end{enumerate}
\end{lem}

\begin{proof}
By differentiating \(f_{a,b}\), it can be verified that this function has exactly one maximum at 
\[
x_{a,b} = \frac{a}{2 \sqrt{b (b-a)}}.
\]
Plugging the maximum value \(x_{a,b}\) into the function \(f_{a,b}\) yields
\[
\begin{split}
f_{a,b}(x_{a,b}) & = \frac{1}{2}\log(1-a) - \frac{2b}{a} \frac{a^2}{4b (b-a)} + \frac{a^2}{4b (b-a)} \\
& \quad + \frac{a}{2 \sqrt{b (b-a)}} \sqrt{1 + \frac{a^2}{4b(b-a)} } + \sinh^{-1} \left (\frac{a}{2 \sqrt{b (b-a)}}  \right ) \\
& =\frac{1}{2}\log(1-a) +  \frac{a (a-2b)}{4b (b-a)} - \frac{a (a-2b)}{4b (b-a)} + \log \left (\sqrt{\frac{b}{b-a}} \right) = \frac{1}{2} \log \left (\frac{b(1-a)}{ b-a} \right),
\end{split}
\]
where in the second line we used the identity \(\sinh^{-1}(x) = \log(x+\sqrt{x^2+1})\). By assumption \(b\geq 1\), thus \(f_{a,b}(x_{a,b}) \leq 0\). In particular, \(f_{a,b}(x) \leq 0\) for all \(x\). Moreover, we have that \(f_{a,b}(x_{a,b})\) equals zero if and only if \(b=1\). 
\end{proof}

\printbibliography
\end{document}